\documentclass[12pt, a4paper, oneside
]{book}

\pdfoutput=1
\usepackage[dvipdfmx]{graphicx}
\usepackage{bmpsize}

\usepackage{vuwthesis} 

\usepackage{palatino} 

\usepackage{url} 

\usepackage{jmbarrett}
\usepackage{capt-of}
\usepackage{mathtools}
\usepackage{pdflscape}

\usepackage{tikz}
\usetikzlibrary{matrix}
\usetikzlibrary{calc}
\definecolor{dgreen}{rgb}{0,0.5,0}
\definecolor{dred}{rgb}{0.8,0,0}

\makeatletter
\@removefromreset{theorem}{section}
\@addtoreset{theorem}{chapter}
\makeatother

\renewcommand{\D}[2][]{{\Delta^{#1}_{#2}}}
\newcommand{\Si}[2][]{{\Sigma^{#1}_{#2}}}
\newcommand{\Ppi}[2][]{{\Pi^{#1}_{#2}}}

\newcommand{\Ind}{\mathsf{I}}
\newcommand{\Ln}{\mathsf{L}}
\newcommand{\Bd}{\mathsf{B}}
\newcommand{\Cmp}{\mathsf{C}}
\newcommand{\BCmp}{\mathsf{BC}}

\newcommand{\ISi}[2][]{\Ind\Si[#1]{#2}}
\newcommand{\LSi}[2][]{\Ln\Si[#1]{#2}}
\newcommand{\LPi}[2][]{\Ln\Ppi[#1]{#2}}
\newcommand{\BSi}[2][]{\Bd\Si[#1]{#2}}

\newcommand{\RAD}{\mathsf{RAD}}

\newcommand{\Ls}{\Ll_2}             
\newcommand{\ZL}{\mathsf{ZL}}       

\newcommand{\zerojump}{{\varnothing'}}
\newcommand{\zerodouble}{{\varnothing''}}
\newcommand{\Inf}{\texttt{Inf}}     

\DeclareMathOperator{\Frac}{Frac}           
\newcommand{\loc}[2]{\operatorname{Loc} \left( {#1}, {#2} \right)}
\newcommand{\idquo}[2]{{#1}\,{:}\,{#2}}     
\renewcommand{\div}{\mid}
\DeclareMathOperator{\lcm}{lcm}
\DeclareMathOperator{\Irr}{Irr}             

\newcommand{\Bez}{B\'ezout}
\DeclareMathOperator{\degr}{deg}            
\DeclareMathOperator{\cont}{cont}           
\newcommand{\finthis}{{\color{red} finish this}}

\newenvironment{proofcite}[1]{
    \begin{proof}[Proof\: {\normalfont\cite{#1}}]
}{
    \end{proof}
}

\begin{document}

\frontmatter


\title{Reverse mathematics of rings}
\author{Jordan Mitchell Barrett}

\subject{Mathematics}
\abstract{Using the tools of reverse mathematics in second-order arithmetic, as developed by Friedman, Simpson, and others, we determine the axioms necessary to develop various topics in commutative ring theory. Our main contributions to the field are as follows. We look at fundamental results concerning primary ideals and the radical of an ideal, concepts previously unstudied in reverse mathematics. Then we turn to a fine-grained analysis of four different definitions of Noetherian in the weak base system $\RCA+\ISi{2}$. Finally, we begin a systematic study of various types of integral domains: PIDs, UFDs and \Bez\ and GCD domains.}

\mscthesisonly




\maketitle

\chapter*{Acknowledgments}\label{C:ack}


The document you are now reading would not have been possible without the support of numerous people. 
First and foremost, I'd like to thank my supervisor Dan Turetsky, who provided an invaluable source of support, ideas and knowledge throughout my Masters. I appreciate the generous scholarship from Victoria University of Wellington which made my studies possible. Thanks also to Long Qian and Valentino Vito for friendship, laughs and thoughtful discussions.

I am grateful to Rod, Noam, and especially Martino, for their support over the years. Thanks to the other members of the logic group (Thomas, Diamant, Linus, Andre) and the proof assistants group (Marco, Julian) for keeping me sane. 


Thanks to everyone else in the School of Mathematics and Statistics that made it a lovely place to learn and grow. There are too many to name, but I'd particularly like to mention Joseph, Liam, Amber, Sahas, Steve, Astrid, Lisa, Evelyn and Matthew. A special mention to Caitlin, Alec and the other admin staff who keep things running so smoothly.


Thanks to all my friends, especially Kaspar, my bandmates in Solid Walls of Sound, and the wonderful peeps in Vic Uke. Finally, thanks are owed to my family for their endless love and support throughout. Mum, Jacob, Lottie \& Maddie, Dad \& Mandy, Ani---I love you all.

\tableofcontents


\mainmatter



\chapter{Introduction}
\label{sec:intro}

This thesis concerns the intersection of two distinct areas of mathematics: commutative ring theory, and mathematical logic. The first, commutative ring theory, also known as \textit{commutative algebra}, has its roots in classical algebraic number theory and
algebraic geometry \cite{kleiner_numbers_1998}.

19th century number theory was concerned with problems such as solvability of Diophantine equations, or of polynomial congruences, over the integers $\Z$. A common technique that emerged was to extend $\Z$ by an \textit{algebraic integer}, obtaining new domains such as $\Z[\sqrt{2}i]$, $\Z[i]$, or $\Z[\omega]$ for $\omega$ a primitive root of unity. If these new domains had \textit{unique factorisation}, then one could draw conclusions about the original equations or congruences.

Unfortunately, these domains failed to have unique factorisation in many important cases. Kummer's idea was to further add ``ideal primes'' to the domain to restore unique factorisation \cite{kummer_preupreber_1847,kummer_zur_1847}. While brilliant, Kummer's ideas were vague, and Dedekind later put Kummer's work on a rigorous footing by giving the modern definition of an \textit{ideal} in a ring \cite{dedekind_supplement_1871}. Two particular features of Dedekind's work---a focus on axiomatic methods, as well as an acceptance of nonconstructive procedures---marked a new style of mathematics which would come to dominate the 20th century, ultimately paving the way for mathematical logic.

\newcommand{\Pp}{\mathcal{P}}
On the other hand, algebraic geometry is concerned with \textit{algebraic varieties}. Given a fixed set of polynomials $\Pp \subseteq \R[x_1,\ldots,x_n]$, the variety $V_\Pp$ is the set of points in $\R^n$ satisfying the equations $p(x_1,\ldots,x_n) = 0$ for all $p \in \Pp$. These are the higher-dimensional generalisation of \textit{algebraic curves} (which are the case $n=2$). Now, to any variety $V$, we can assign the set $I(V)$ of polynomials which vanish on $V$; this is an ideal in $\R[x_1,\ldots,x_n]$. Hence, we can study algebraic varieties by studying ideals in polynomial rings. This correspondence was exploited to great effect by Hilbert in his basis theorem \cite{hilbert_ueber_1890} and \textit{Nullstellensatz} \cite{hilbert_ueber_1893}, and later by Lasker \cite{lasker_zur_1905} and Macauley \cite{macaulay_resolution_1913}.

Attempting to create a general theory encompassing all these ideas, Fraenkel gave the first abstract definition of a ring in \cite{fraenkel_preupreber_1915}, and Sono gave the modern definition soon after \cite{sono_congruences_1917}. This opened the door to the pioneering work of Noether, which established abstract ring theory as a subject. Specifically, Noether \cite{noether_idealtheorie_1921} generalised the results of Hilbert, Lasker and Macauley to what are now called \textit{Noetherian rings}, and later recast the work of Dedekind et al in an abstract setting \cite{noether_abstrakter_1927}.

Meanwhile, a mathematical revolution had been brewing. Traditionally, mathematics had been concerned with finite objects and constructive procedures, and grounded in reality \cite{evesIntroductionHistoryMathematics1969, metakidesIntroductionNonrecursiveMethods1982}. The work of Dedekind, Hilbert, Peano, Cantor and others in the late 1800s marked a departure from this, thereby ushering in modern, abstract pure mathematics. This new style of mathematics was distinguished by its focus on abstraction and the axiomatic method, and acceptance of nonconstructive proofs---those which prove the existence of an object without actually constructing an example.

Another feature of this new mathematics was the acceptance of \textit{completed infinity} rather than just \textit{potential infinity}---the idea that infinite sets could be manipulated as mathematical objects in their own right. This idea perhaps appeared first in Dedekind's work on ideals \cite{dedekind_supplement_1871,kleiner_numbers_1998}. %
Cantor was the first to systematically study infinity, founding the field of set theory with his seminal work on cardinals \cite{cantorUeberEigenschaftInbegriffs1874} and ordinals \cite{cantorUeberUnendlicheLineare1883}.
%
As set theory developed, paradoxes arose (most notably Russell's), and the need for a careful and rigorous foundation for mathematics became clear. One such foundation was provided by $\ZFC$ in the 1920s \cite{zermeloUberGrenzzahlenUnd1930a}.

Cantor's work provided new impetus to mathematical logic, a small subfield of mathematics developed by Boole, De Morgan, and Peano in the mid-to-late 1800s \cite{booleInvestigationLawsThought1854, demorganFormalLogic1847, peanoArithmeticesPrincipiaNova1889}. Around this time, the ideas of computation, mathematical truth and mathematical proof were formalised for the first time. By the 1930s, logic was a thriving area of mathematics---highlights included G\"odel's (in)completeness theorems \cite{godelUberVollstandigkeitLogikkalkuls1929, goedelUeberFormalUnentscheidbare1931}, Turing's negative solution to the \textit{Entscheidungsproblem} \cite{turingComputableNumbersApplication1937}, Tarski's development of model theory \cite{vaughtAlfredTarskiWork1986}, and Hilbert's work on proof theory \cite{hilbertGrundlagenMathematik1934} and geometry \cite{hilbertGrundlagenGeometrie1899}.

A later development in logic was reverse mathematics, initiated by Harvey Friedman in the late 1960s \cite{friedmanSubsystemsSetTheory1967, friedmanBarInductionPi1969}. Reverse mathematics asks, for a given theorem of mathematics $\varphi$, ``what axioms are really necessary to prove $\varphi$?'' More broadly, it studies the logical implications between foundational principles of mathematics. An early example was the discovery of non-Euclidean geometries, thereby proving the independence of the parallel postulate from Euclid's other axioms \cite{lobachevskyConciseOutlineFoundations1829, bolyaiAppendixScientiamSpatii1832}. Another early result, more in the style of reverse mathematics, was the demonstration that over $\ZF$, the axiom of choice, Zorn's lemma, and the well-ordering principle are all pairwise equivalent \cite{birkhoffLatticeTheory1940, fraenkelFoundationsSetTheory1958, traylorEquivalenceAxiomChoice1962}.

Traditionally, reverse mathematics is done in second-order arithmetic, in which there are two types of objects: natural numbers $n, m, k, \ldots$, and sets of natural numbers $A, B, C, \ldots$, and quantification is allowed over both types of objects. Restricting oneself to natural numbers may seem unnecessary limiting, but this is not so. In fact, most mathematics deals with countable or ``essentially countable'' objects (such as separable metric spaces), and so can be formalised in second-order arithmetic. This includes virtually all ``classical'' mathematics, or that taught in undergraduate courses \cite[xiv]{simpson_subsystems_2009}.

In practice, reverse mathematics involves attempting to prove a theorem $\varphi$ of ``ordinary'' mathematics in a weak subsystem $\Ss$ of second-order arithmetic. But, supposing we can do this, how do we know we've found the optimal (weakest) system? The empirical phenomenon is thus:
\begin{center}
	\textit{``When the theorem is proved from the right axioms,\\ the axioms can be proved from the theorem.''}
\end{center}\vspace*{-5.7mm}
\begin{flushright}
	\textit{---Harvey Friedman} \cite{friedmanSystemsSecondOrder1974}
\end{flushright}
This is the ``reverse'' part of reverse mathematics. Having proved $\varphi$ from $\Ss$, to show this is optimal, we want to demonstrate a \textit{reversal} of $\varphi$: a proof of $\Ss$ from $\varphi$. This means that $\varphi$ cannot be proved in a weaker system $\Ss'$, because if it could, then $\Ss'$ would also prove $\Ss$ via $\varphi$, meaning $\Ss'$ is not actually a weaker system after all.

The utility of reverse mathematics is abundant. Apart from its obvious use in finding the ``best'' proof of a given statement $\varphi$, it also gives us a way to quantify how nonconstructive or noncomputable $\varphi$ is. The idea is that stronger subsystems correspond to more nonconstructive power, 
so the ``constructiveness'' of $\varphi$ is inversely proportional to the strength of the systems $\Ss$ in which $\varphi$ can be proved \cite{friedman_countable_1983}.
Similarly, many theorems guarantee a solution to a given problem---reverse mathematics then tells us how complex the solution could be relative to the problem, which can be made precise in terms of computability.

Here is an example of reverse mathematics in ring theory. The usual way to prove that every commutative ring has a prime ideal is to prove that it has a maximal ideal (Krull's theorem), and then prove every maximal ideal is prime. However, Friedman, Simpson and Smith showed that the existence of maximal ideals is equivalent to the system $\ACA$, whereas the existence of prime ideals is equivalent to the strictly weaker system $\WKL$ \cite{friedman_countable_1983}. This shows the usual proof strategy is not optimal---there is a ``better'' way to prove the existence of prime ideals, which doesn't require the stronger assumption that maximal ideals exist. In terms of computability, this shows that maximal ideals can be more ``noncomputable'' than prime ideals---more precisely, given a computable ring, its maximal ideals could all be as complex as the \textit{halting problem}, while we can always compute a prime ideal from a \textit{PA degree}.

In this thesis, we study the reverse-mathematical content of various theorems of ring theory. We will begin by reviewing basic ideas from ring theory (\S\ref{sec:rings}), and from logic, computability and reverse mathematics (\S\ref{sec:logic}). We then proceed to study the following key ideas from commutative algebra:
\begin{itemize}
    \item \textit{Primary ideals} (\S\ref{sec:prim-rad}), which are the ideals $I$ such that whenever $ab \in I$, then $a \in I$ or $b^n \in I$ for some $n$.
    
    \item The \textit{radical} $\sqrt{I}$ of an ideal $I$ (\S\ref{sec:prim-rad}), which is the set of $r$ such that some $r^n$ is in $I$.
    
    \item \textit{Noetherian rings} (\S\ref{sec:noeth}), which have many equivalent definitions; one of the more popular is that every ideal is finitely generated.
    
    \item Several classes of \textit{integral domains} (\S\ref{sec:int-doms}), including PIDs, UFDs, \Bez\ and GCD domains, their properties, and the relations between them.
\end{itemize}

In relation to the philosophy of reverse mathematics, we find that all the results we examine are provable in $\ACA$.\footnote{This seems to be true for algebra in general, with the notable exception of some structure theorems and (ordinal) invariant results \cite{simpson_subsystems_2009}.} We will show that many important results, such as the equivalence of different notions of Noetherian, actually \textit{require} $\ACA$. Thus, we conclude that $\ACA$ is the right axiom system in which to develop (most) of classical commutative algebra. This thesis includes many new, original results---some of the more important ones are:
\begin{itemize}
    \item Theorem \ref{thm:wkl-maximal-primary}:
    $\WKL$ is equivalent to ``$\sqrt{I}$ maximal $\implies$ $I$ primary''.
    
    \item Theorem \ref{thm:chain-s1-ideals}: $\RCA+\ISi{2}$ proves the equivalence of weak and strict chain conditions on $\Si{1}$-ideals.
    
    \item Theorems \ref{thm:seq->d1n} and \ref{thm:wkl-d1acc->seq}, showing that $\WKL$ or $\ACA$ are equivalent to the agreement of several definitions of Noetherian.
    
    \item Theorem \ref{thm:aca-prime-princ}: $\ACA$ proves that if $R$ is an integral domain in which every prime $\Si{1}$-ideal is principal, then $R$ is a $\Si{1}$-PID.
    
    \item Theorem \ref{thm:pid->dhd} and Corollary \ref{cor:aca-pid->dhd}: $\ACA$ is equivalent to ``every PID admits a Dedekind--Hasse norm''.
    
    \item Theorem \ref{thm:pid-irr-pi2}, where we construct a PID whose set of primes is $\Ppi{2}$ complete.
\end{itemize}

	
	

\section{Notational conventions}

The following notational conventions will apply to this thesis:

\begin{itemize}
    \item We will use $\N$ to denote (the underlying set of) the model of arithmetic we are working inside, and $\omega$ to denote the standard model $\{ 0, 1, 2, \ldots \}$. Generally, the distinction will not be important.
    
    \item For a ring $R$, we will use $R[\bar{x}]$ as an abbreviation for $R[x_0, x_1, \ldots]$, the polynomial ring over $R$ in infinitely many indeterminates.

    \item In mathematics, pairs, tuples and sequences are commonly denoted using parentheses, e.g.\ $(x,y)$, $(a_0, a_1, \ldots)$. However, in ring theory, it is also common to use parentheses to denote the \textit{ideal} $(A)$ generated by a collection of elements $A \subseteq R$. To avoid confusion, we will try to consistently use angle brackets $\ang{x,y}$, $\ang{a_0, a_1, \ldots}$ to denote a pair, tuple or sequence, and reserve parentheses for ideals.
    
    \item $\varepsilon$ or $\ang{}$ will denote the empty sequence or tuple.
    
    \item The $\mathsf{SANS}\ \mathsf{SERIF}$ font will generally be reserved for subsystems and axioms of second-order arithmetic.
    
    \item We will use $\varphi_0, \varphi_1, \ldots, \varphi_e, \ldots$ to denote a standard listing of the partial computable functions, and $W_0, W_1, \ldots, W_e, \ldots$ to denote a listing of the c.e.\ sets.
    
    \item We may use $A^\complement = \{ x: x \notin A \}$ to denote the (absolute) complement of a set $A$, particularly for sets of natural numbers.
	
	\item For mathematical statements $\varphi$ and $\psi$, we use $\varphi \proves \psi$ (``$\varphi$ proves $\psi$'') to mean there is a proof of $\psi$ from $\varphi$. This notation extends to formal systems, e.g.\ $\Ss \proves \varphi$ means there is a proof of $\varphi$ in the formal system $\Ss$.
	
	\item For a statement $\varphi$ and a structure $\M$, we use $\M \models \varphi$ (``$\M$ models $\varphi$'') to mean the statement $\varphi$ is true in $\M$. Similarly, $\M \models \Ss$ means that all axioms of the formal system $\Ss$ are true in $\M$.
\end{itemize}
\chapter{Ring theory}
\label{sec:rings}

Here, we quickly review the basic notions of ring theory that we will need, as covered in any basic algebra textbook \cite{atiyah_introduction_1994}. For us, ``ring'' will mean ``commutative ring with unity'', unless explicitly stated otherwise.

\begin{definition}
    A \textit{ring} is a set $R$, equipped with constants $0_R, 1_R \in R$, and binary operations $+, \cdot$ on $R$ (called \textit{addition} and \textit{multiplication}, respectively), such that:
    \begin{enumerate}
        \item $R$ is an abelian group under addition, with additive identity $0_R$.
        
        \item Multiplication is associative: $(a \cdot b) \cdot c = a \cdot (b \cdot c)$ for all $a, b, c \in R$.
        
        \item Multiplication is commutative: $a \cdot b = b \cdot a$ for all $a, b \in R$.
        
        \item $1_R$ is a two-sided multiplicative identity.
        
        \item Multiplication distributes both ways over addition: $a \cdot (b + c) = (a \cdot b) + (a \cdot c)$ and $(a + b) \cdot c = (a \cdot c) + (b \cdot c)$ for all $a, b, c \in R$.
    \end{enumerate}
\end{definition}

As usual, we will often omit the dot for multiplication, and instead denote it by juxtaposition, i.e.\ $ab$ instead of $a \cdot b$.

Given a ring $R$, we can construct a larger ring $R[x]$ of ``polynomials over $R$'' in the variable $x$. We do this by ``freely'' adding the variable $x$, i.e.\ asserting no relationship between $x$ and elements of $R$. A formal construction follows.

\begin{definition}
    For a ring $R$, the \textit{polynomial ring} $R[x]$ is defined as follows:
    \begin{itemize}
        \item The underlying set of $R[x]$ is the collection $$\big\{ \ang{a_0, \ldots, a_n} \in R^{<\omega} : a_n \neq 0_R \big\} \cup \{ \ang{} \}$$
        We denote $\ang{a_0, \ldots, a_n}$ by $a_n x^n + a_{n-1} x^{n-1} + \cdots + a_1 x + a_0$.
        
        \item $0_{R[x]} = 0_R = \ang{}$ and $1_{R[x]} = 1_R = \ang{1_R}$.
        
        \item Addition and multiplication in $R[x]$ are defined as follows:
        \begin{align*}
            \left( \sum_{i=0}^n a_i x^i \right) + \left( \sum_{i=0}^m b_i x^i \right) &= \sum_{i=0}^{\max\{n,m\}} (a_i+b_i) x^i\\[2mm]
            \left( \sum_{i=0}^n a_i x^i \right) \cdot \left( \sum_{i=0}^m b_i x^i \right) &= \sum_{i=0}^{m+n} \left[ \sum_{j=0}^i a_j b_{i-j} \right] x^i
        \end{align*}
    \end{itemize}
\end{definition}

Given a polynomial ring $R[x]$, we could repeat the construction to get $\big( R[x] \big)[y]$, which we write simply as $R[x,y]$ for brevity. Iterating this construction, we get an increasing sequence $R[x_0], R[x_0,x_1], R[x_0,x_1,x_2], \ldots$. We will use $R[\bar{x}] = R[x_0,x_1,\ldots]$ to refer to the limit of this sequence.

\begin{definition}
    In a ring $R$, an \textit{ideal} $I \subseteq R$ is a subset of $R$ such that for all $a,b \in I$ and $r,s \in R$, we have $ar + bs \in I$.
\end{definition}

Ideals are important since they give us a way to create new rings:

\begin{definition}
    Given a ring $R$ and ideal $I \subseteq R$, the \textit{quotient ring} $R/I$ is defined as follows:
    \begin{itemize}
        \item The underlying set is the quotient of $R$ by the equivalence relation $r \sim s \iff r-s \in I$. We denote the equivalence class of $r$ by $r+I$.
        
        \item $0_{R/I} = 0_R + I$ and $1_{R/I} = 1_R + I$.
        
        \item $(a+I)+(b+I) = (a+b)+I$ and $(a+I)(b+I) = (ab)+I$. These operations are well-defined.
    \end{itemize}
\end{definition}

We recall some important ways of creating ideals.

\begin{definition}
    Let $R$ be a ring and $A \subseteq R$ an arbitrary subset. The \textit{ideal generated by $A$} is the set $$(A) \defeq \{ a_1 r_1 + \cdots + a_n r_n : a_1, \ldots, a_n \in A, r_1, \ldots,\, r_n \in R \}$$
    If $A = \{ a_1, \ldots, a_n \}$ is a finite set, we write $(a_1, \ldots, a_n)$ and say this ideal is \textit{finitely generated}. If $A = \{ a \}$, we say $(a)$ is \textit{principal}.
\end{definition}

\begin{definition}\label{defn:idquo}
    Given two ideals $I, J \subseteq R$, the \textit{ideal quotient of $I$ by $J$} is the set $$\idquo{I}{J} \defeq \{ r \in R : (\forall j \in J)(rj \in I) \}$$
\end{definition}

The most common case of Definition \ref{defn:idquo} is when $J = (a)$ is a principal ideal. In this case, we will abuse notation and write $\idquo{I}{a}$ instead of $\idquo{I}{(a)}$. $\idquo{I}{a}$ also admits a simpler definition here, as $\idquo{I}{a} = \{ r \in R : ra \in I \}$.

We now recall an important subclass of the commutative rings.

\begin{definition}
    An element $a \in R$ is a \textit{zero-divisor} if there is $b \neq 0_R$ such that $ab = 0$.
\end{definition}

\begin{definition}
    A ring $R$ is an \textit{integral domain} if it has no nonzero zero-divisors, i.e.\ whenever $ab = 0$, then $a=0$ or $b=0$.
\end{definition}

Integral domains satisfy cancellation of multiplication: if $a \neq 0$ and $ab = ac$, then $b = c$. Indeed, this is an alternative characterisation of integral domains.

\begin{definition}\label{defn:field-frac}
    Given an integral domain $R$, its \textit{field of fractions} $\Frac(R)$ is the ring defined as follows:
    \begin{itemize}
        \item The underlying set of $\Frac(R)$ is the quotient of $\{ (r,s) \in R^2 : s \neq 0_R \}$ by the equivalence relation $(r,s) \sim (r',s') \iff rs' = r's$. We denote the equivalence class of $(r,s)$ by $r/s$.
        
        \item $0_{\Frac(R)} \defeq 0_R / 1_R$ and $1_{\Frac(R)} \defeq 1_R / 1_R$.
        
        \item $(r/s) + (r'/s') = (rs' + r's)/(ss')$ and $(r/s)(r'/s') = (rr')/(ss')$. These operations are well-defined.
    \end{itemize}
\end{definition}

We can verify that $\Frac(R)$ is indeed a ring, and in fact, it is a field---every nonzero element has an inverse. We can naturally view $R$ as a subring of $\Frac(R)$ via the embedding $r \mapsto r/1$.

There is an important generalisation of Definition \ref{defn:field-frac}, which covers both the case when $R$ is not an integral domain, and when we don't want every element of $R$ to be a denominator.

\begin{definition}\label{defn:loc}
    Suppose $M \subseteq R$ is multiplicatively closed, contains $1_R$, and contains no zero-divisors. The \textit{localisation of $R$ at $M$}, $\loc{R}{M}$, is the ring whose underlying set is the quotient of $\{ (r,m) \in R^2 : m \neq 0_R \}$ by the equivalence relation $(r,m) \sim (r',m') \iff rm' = r'm$, and where the operations are defined as in Definition \ref{defn:field-frac}.
\end{definition}

Intuitively, $\loc{R}{M}$ is obtained from $R$ by allowing division by the elements of $M$. As before, $R$ is naturally a subring of $\loc{R}{M}$ via the embedding $r \mapsto r/1$.

\newcommand{\sat}[1]{\overline{#1}}
\begin{definition}\
    \begin{enumerate}
        \item A multiplicatively closed set $M$ is \textit{saturated} if whenever $ab \in M$, then both $a,b \in M$.
        
        \item $\sat{M} \defeq \{ r \in R : (\exists s \in R)(rs \in M) \}$ is the smallest saturated set containing $M$, and is called the \textit{saturation} of $M$.
    \end{enumerate}
\end{definition}

The concept of saturation is important because:

\begin{theorem}
    Given multiplicatively closed sets $M, N \subseteq R$, we have that $\loc{R}{M} \cong \loc{R}{N}$ canonically iff $\sat{M} = \sat{N}$. In particular, it is always true that $\loc{R}{M} \cong \loc{R}{\sat{M}}$.
\end{theorem}

Essentially, when we localise at $M$, we are really localising at its saturation $\sat{M}$. To see why, take $r \in \sat{M}$, i.e.\ there is $s \in R$ with $rs \in M$. Then, $s/rs$ is an element of $\loc{R}{M}$. But $r \cdot (s/rs) = 1$, so we have introduced an inverse for $r$, even if $r \notin M$. If we check that $M$ is saturated, then we know that we've \textit{only} added inverses for elements of $M$, and nothing else.

An important case of localisation is the so-called ``localisation at a prime ideal'', which is really the localisation at the complement of a prime ideal.

\begin{definition}
    An ideal $I \subseteq R$ is \textit{prime} if whenever $ab \in I$, then $a \in I$ or $b \in I$.
\end{definition}

\begin{example}
    Suppose $P \subseteq R$ is a prime ideal. Then, $M \defeq R \setminus P$ is multiplicatively closed, so we can take the localisation $\loc{R}{M}$. This is sometimes called the \textit{localisation of $R$ at $P$} and denoted $R_P$.
\end{example}

Another important example of localisation is the total quotient ring, which generalises the field of fractions construction to non-integral domains.

\begin{definition}\label{defn:tot-quo-ring}
    For any ring $R$, the set $M$ of all non-zero-divisors of $R$ is multiplicatively closed and contains $1_R$. The \textit{total quotient ring} $\Frac(R)$ is the localisation $\loc{R}{M}$.
\end{definition}

When $R$ is an integral domain, $0_R$ is the only zero-divisor, so Definition \ref{defn:tot-quo-ring} reduces to Definition \ref{defn:field-frac}.
\chapter{Logical prerequisites}
\label{sec:logic}

Here, we review the necessary background material from model theory \cite{marker_model_2002}, first- and second-order arithmetic \cite{hajek_metamathematics_2017,simpson_subsystems_2009}, reverse mathematics \cite{simpson_subsystems_2009}, and computability \cite{soare_recursively_1987, soare_turing_2016}.

\section{Second-order arithmetic}

For us, \textit{arithmetic} will refer to the model-theoretic study of the theory of the natural numbers, $\mathrm{Th}(\N)$. Our reverse-mathematical studies will be done in the traditional setting of \textit{second-order arithmetic}. In first-order arithmetic, we are only allowed to quantify over \textit{elements} of $\N$, while in second-order arithmetic, we may also quantify over \textit{subsets} of $\N$. This greatly increases the expressive power of our logic (for example, we can define well-foundedness, or completeness of $\R$).

We work in a two-sorted model theory, whose sorts are \textit{numbers}, denoted with lowercase letters $n,m,k,\ldots$, and \textit{sets}, denoted in uppercase $A,B,C,\ldots$. Our language is $\Ls \defeq \{ 0,1,+,\cdot,<,\in \}$, where the symbols have the expected types, e.g. $0$ is a number, $+$ takes two numbers and returns another, $\in$ is a binary relation between a number and a set, etc.
%
%

Numerical $\Ls$-terms are defined inductively: $0,1$, and variable symbols $x_i$ are numerical $\Ls$-terms, and if $s,t$ are numerical $\Ls$-terms, then $(s+t)$, $(s \cdot t)$ are too. Numerical $\Ls$-terms represent (possibly nonstandard) natural numbers, and we will use $k$ to abbreviate the numerical term $$k\ \defeq\ \underbrace{1+1+\cdots+1}_{k\text{ times}}$$ The only set $\Ls$-terms are variables $X_i$. We also define $\Ls$-formulae inductively:

\begin{definition}
	The collection of \textit{$\Ls$-formulae} is defined as follows:
	\begin{enumerate}
		\item If $s$, $t$ are numerical terms, and $X$ is a set variable symbol, then $(s=t)$, $(s < t)$ and $(s \in X)$ are formulae.
		
		\item If $\varphi$, $\psi$ are formulae, then $(\lnot \varphi)$, $(\varphi \land \psi)$, $(\varphi \lor \psi)$, $(\varphi \to \psi)$ and $(\varphi \leftrightarrow \psi)$ are formulae.
		
		\item If $\varphi$ is a formula, then $(\forall x)\, \varphi$ and $(\exists x)\, \varphi$ are formulae.
		
		\item If $\varphi$ is a formula, then $(\forall X)\, \varphi$ and $(\exists X)\, \varphi$ are formulae.
	\end{enumerate}
\end{definition}

A \textit{theory} $\Ss$ is simply a set of formulae. We may also use the terms \textit{subsystem} or \textit{formal system}, especially when we are considering the members of $\Ss$ as axioms.

There are a wide variety of $\Ls$-formulae, which we will now classify based on their ``complexity''. 
Our chosen measure of complexity will be based on how many times the quantifiers alternate. This defines a structure known as the \textit{arithmetical hierarchy}.

The lowest level of complexity consists of formulae containing only bounded quantifiers: those of the form $(\forall x)(x < k\, \to\, \psi)$ or $(\exists x)(x < k\, \to\, \psi)$ for some numerical term. We will often abbreviate these to $(\forall x < k)\, \psi$ and $(\exists x < k)\, \psi$ respectively. From there, universal formulae are given $\Pi$ classifications, and existential formulae given $\Sigma$ classifications.

\begin{definition}[(arithmetical hierarchy)]
	Let $\varphi$ be an $\Ls$-formula. We say $\varphi$ is $\D[1]{0}$ or \textit{arithmetical} if it contains no set quantifiers. If $\varphi$ is arithmetical, we assign further classifications to it as follows:
	\begin{enumerate}
		\item $\varphi$ is called $\Si[0]{0}$ and $\Ppi[0]{0}$ if it only contains bounded quantifiers.
		
		\item $\varphi$ is called $\Si[0]{n+1}$ if it is of the form $\varphi = \exists x_1 \cdots \exists x_n\ \psi$, where $\psi$ is $\Ppi[0]{n}$.
		
		\item $\varphi$ is called $\Ppi[0]{n+1}$ if it is of the form $\varphi = \forall x_1 \cdots \forall x_n\ \psi$, where $\psi$ is $\Si[0]{n}$.
		
		\item $\varphi$ is called $\D[0]{n}$ if it is $\Si[0]{n}$, and logically equivalent to a $\Ppi[0]{n}$ formula.
	\end{enumerate}
\end{definition}

\begin{figure}
	\centering
	\begin{tikzpicture}
	\matrix(ah)[matrix of math nodes,row sep=6mm,column sep=6.5mm]{
		\Si{0} & & \Si{1} & & \Si{2} & & \Si{3} & \cdots \\
		\D{0} & \D{1} & & \D{2} & & \D{3} & & \cdots \\
		\Ppi{0} & & \Ppi{1} & & \Ppi{2} & & \Ppi{3} & \cdots \\
	};
	\begin{scope}[every node/.style={sloped,gray}]
	\path (ah-1-1) -- node{$=$} (ah-2-1) -- node{$=$} (ah-3-1);
	\path (ah-2-1) -- node{$\subsetneq$} (ah-2-2) -- node{$\subsetneq$} (ah-1-3) -- node{$\subsetneq$} (ah-2-4) -- node{$\subsetneq$} (ah-1-5) -- node{$\subsetneq$} (ah-2-6) -- node{$\subsetneq$} (ah-1-7);
	\path (ah-2-2) -- node{$\subsetneq$} (ah-3-3) -- node{$\subsetneq$} (ah-2-4) -- node{$\subsetneq$} (ah-3-5) -- node{$\subsetneq$} (ah-2-6) -- node{$\subsetneq$} (ah-3-7);
	\end{scope}
	\end{tikzpicture}
	\caption{The arithmetical hierarchy.}
\end{figure}
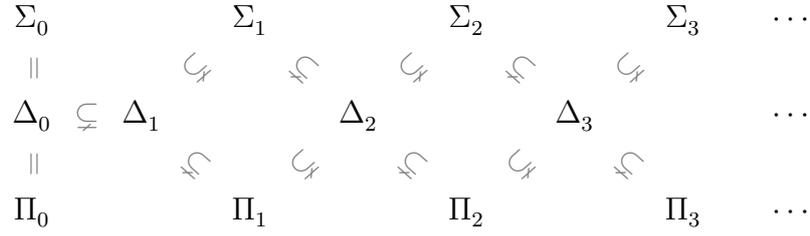

Often
, we will drop the superscript $0$, and just speak of $\Si{n}$ or $\Ppi{n}$ formulae. Now we look at some important principles of arithmetic, which will be pertinent in our study of reverse mathematics.

\begin{definition}\label{defn:ind-ln-bd}
    Let $\varphi(x)$, $\psi(x,y)$ be formulae of arithmetic.
    \begin{enumerate}
        \item The \textit{induction principle for $\varphi$}, $\Ind{\varphi}$, is the statement $$\varphi(0) \land (\forall n)[\varphi(n) \to \varphi(n+1)]\, \to (\forall n)\, \varphi(n)$$
        
        \item The \textit{strong induction principle for $\varphi$}, $\Ind'{\varphi}$, is the statement $$(\forall n) \big( [(\forall m<n)\, \varphi(m)] \to \varphi(n) \big)\, \to (\forall n)\, \varphi(n)$$
        
        \item The \textit{least number principle for $\varphi$}, $\Ln{\varphi}$, is the statement $$ (\exists n)\, \varphi(n) \to (\exists n') \big( \varphi(n') \land (\forall m < n')\, \lnot \varphi(m) \big) $$
        
        \item The \textit{bounding principle for $\psi$}, $\Bd{\psi}$, is the statement $$(\forall k) \big[ (\forall n<k)(\exists m)\, \psi(n,m)\, \to\, (\exists \ell)(\forall n<k)(\exists m<\ell)\, \psi(n,m) \big]$$
    \end{enumerate}
\end{definition}

All of the above principles trivially hold in the standard model of arithmetic $\omega$. However, recall that we have \textit{nonstandard} models of arithmetic, which satisfy the same basic axioms as $\omega$. These nonstandard models are linear orders of type $\omega + \Z \cdot K$ for some linear order $K$ \cite{henkin_completeness_1950}. In one of these models, the above principles could fail. For example, $\Ln{\varphi}$ fails if $\varphi$ is true exactly on the $\Z \cdot K$ part, while $\Bd{\psi}$ could fail if the bound $k$ is nonstandard.

For a \textit{class} of formulae $\Gamma$, we define $\Ind{\Gamma}$ to be the theory consisting of all the statements $\Ind{\varphi}$ for all (appropriate) $\varphi \in \Gamma$. The theories $\Ind'{\Gamma}$, $\Ln{\Gamma}$, $\Bd{\Gamma}$ are defined analogously. We will usually take $\Gamma$ to be a classification in the arithmetical hierarchy: for example, $\Ind{\Si{2}}$ is induction for all $\Si{2}$ formulae.

These arithmetical principles are closely related to each other:

\begin{theorem}\label{thm:ind-ln}
    For every formula of arithmetic $\varphi(x)$, $\Ind'{\varphi} \equiv \Ln{(\lnot\varphi)}$.
\end{theorem}

\begin{theorem}[\cite{paris_sigma_n-collection_1978,hajek_metamathematics_2017}]\label{thm:paris-kirby}
    Over $\PA^- + \ISi{0}$, for all $n \in \N$:
    \begin{enumerate}
        \item $\ISi{n}$, $\Ind\Ppi{n}$, $\Ind'\Si{n}$, $\Ind'\Ppi{n}$, $\LSi{n}$, $\LPi{n}$ are all equivalent.
        
        \item $\Ln\D{n+1}$, $\BSi{n+1}$, $\Bd\Ppi{n}$, $\Bd\D{n+1}$ are all equivalent, and imply $\Ind\D{n+1}$.
        
        \item $\ISi{n+1} \implies \BSi{n+1} \implies \ISi{n}$, and these implications are strict.
    \end{enumerate}
\end{theorem}

\begin{theorem}[\cite{slaman_sigma_n-bounding_2004}]\label{thm:slaman}
    Over $\PA^- + \ISi{0} + \mathsf{exp}$, $\Ind\D{n+1}$ is equivalent to $\Ln\D{n+1}$ (and hence to $\BSi{n+1}$, $\Bd\Ppi{n}$, $\Bd\D{n+1}$).
\end{theorem}

\begin{figure}[t!]
    \centering
    \includegraphics[width=0.95\textwidth]{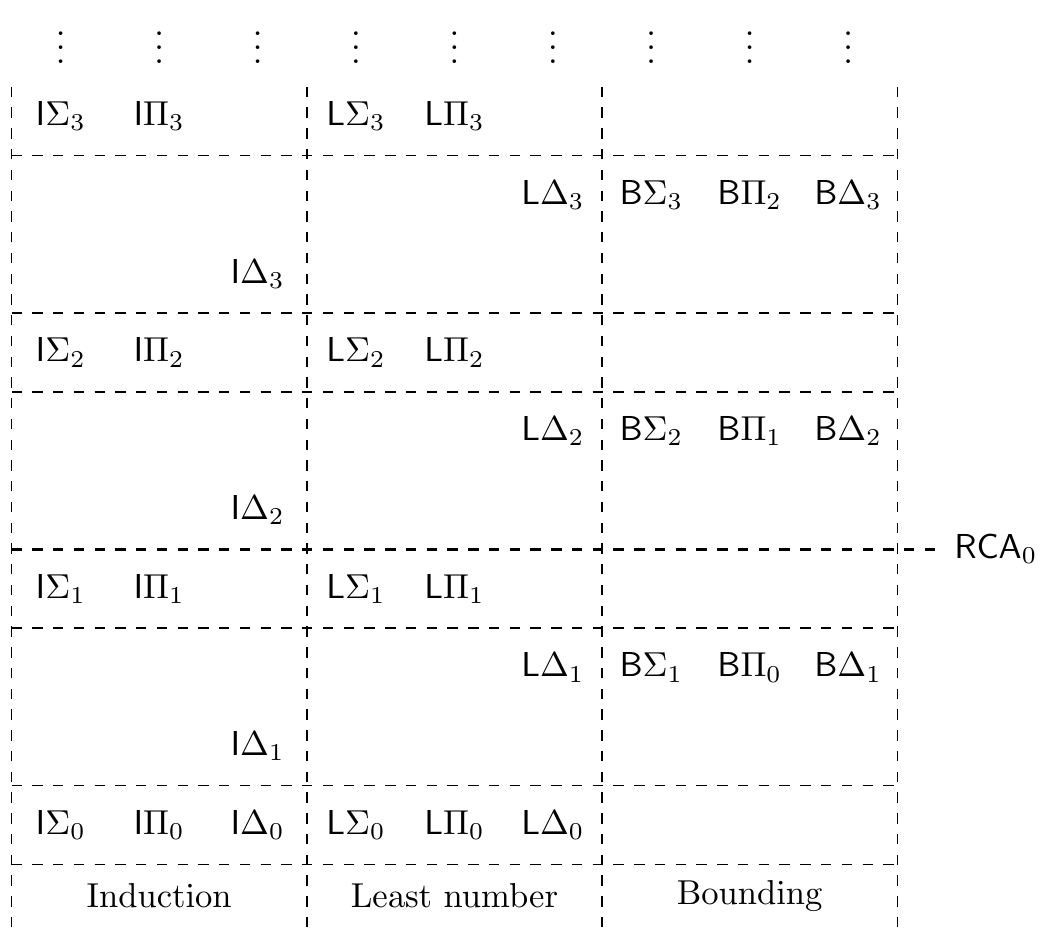}
    \caption{The induction, bounding, and least number principles. Vertical positioning denotes equivalence over $\PA^- + \ISi{0}$, and the horizontal dotted lines denote equivalence over $\PA^- + \ISi{0} + \mathsf{exp}$.}
    \label{fig:ind-bd-ln}
\end{figure}

We will not define $\PA^-$ and $\mathsf{exp}$, but they are extremely weak base theories, which will be subsumed by our chosen base theory $\RCA$. Hence, we can simply assume the equivalences in Theorems \ref{thm:paris-kirby} and \ref{thm:slaman}, and our arithmetical principles are split into two families of equivalence classes, as shown in Figure \ref{fig:ind-bd-ln}.

We now prove a lemma about induction, which will be useful later.

\begin{definition}
    For a formal system $\Ss$, the \textit{inductive formulae $\I(\Ss)$ for $\Ss$} is the collection of formulae $\varphi$ such that $\Ss \models \Ind'\varphi$.
\end{definition}

\begin{lemma}\label{lem:ind-bool-comb}
    If $\I(\Ss)$ is closed under negation, then it is closed under arbitrary Boolean combinations.
\end{lemma}

\begin{proof}
    Let $\Gamma \defeq \I(\Ss)$: then $\Ss \models \Ind'\Gamma$ (and indeed, $\Gamma$ is maximal with this property). Since $\Gamma$ is closed under negation, $\Ss \models \Ln\Gamma$ too by Theorem \ref{thm:ind-ln}. Since the connectives $\lnot$ and $\land$ are complete for propositional logic \cite[98]{smith_introduction_2003}, and we already know $\Gamma$ is closed under $\lnot$, it suffices to prove closure under $\land$.
    
    So, pick $\varphi(x), \psi(x) \in \Gamma$: we want to show $\Ss \models \Ind'(\varphi \land \psi)$. For each $n$, assume $$(\forall m<n)[\varphi(m) \land \psi(m)] \to \varphi(n) \land \psi(n)$$
    
    By contradiction, suppose there is $a$ with $\lnot \big( \varphi(a) \land \psi(a) \big)$, and without loss of generality, suppose $\lnot \varphi(a)$. By $\Ln\Gamma$, we can assume $a$ is minimal. By the inductive assumption, there is $b < a$ with $\lnot \big( \varphi(b) \land \psi(b) \big)$, hence $\lnot \psi(b)$ by minimality of $a$. Again by $\Ln\Gamma$, assume $b$ is minimal. Then, by the inductive assumption, there is $c < b$ with $\lnot \big( \varphi(c) \land \psi(c) \big)$. But this contradicts minimality of $a$ and $b$.
\end{proof}

Lemma \ref{lem:ind-bool-comb} allows us to induct on arbitrary Boolean combinations of $\Si{1}$ and $\Ppi{1}$ formulae in $\RCA$, which will prove useful.

\section{Reverse mathematics}

For a system $\Ss$ of second-order arithmetic, a \textit{model} $\M$ of $\Ss$ consists of a set $\N$ and a collection of subsets of $\N$, with appropriate interpretations for the $\Ls$-symbols, so that all the axioms of $\Ss$ hold. All the systems $\Ss$ we consider will include the basic axioms of first-order arithmetic, hence $\N$ will be restricted to range over nonstandard models $\omega + \Z \cdot K$, where the $\Ls$-symbols are given the usual meanings. If $K = \varnothing$, i.e.\ $\N = \omega$, we call $\M$ an \textit{$\omega$-model}, and then $\M$ is determined by a collection of subsets of $\omega$.

We are almost ready to define the subsystems $\RCA$, $\WKL$, $\ACA$ of second-order arithmetic, which will form the basis of our work. Some of these subsystems will include the following second-order comprehension principle:

\begin{definition}
    Let $\varphi(x)$ be a formula of arithmetic. The \textit{comprehension principle for $\varphi$} is the sentence $\Cmp\varphi \defeq (\exists X)(\forall n)(n \in X \leftrightarrow \varphi(n))$.
\end{definition}

Essentially, $\Cmp\varphi$ asserts that the set $X = \{ n : \varphi(n) \}$ exists. As before, given a class $\Gamma$ of formulae, we write $\Cmp\Gamma$ for the theory consisting of $\Cmp\varphi$ for all (appropriate) $\varphi \in \Gamma$.

\begin{definition}
	$\RCA$ is the subsystem consisting of the basic axioms of first-order arithmetic, $\ISi[0]{1}$, and $\Cmp\D[0]{1}$.
\end{definition}

$\RCA$ is the system in which reverse mathematics is usually done. Although $\RCA$ doesn't give us $\Si{1}$ comprehension, we do get \textit{bounded} $\Si{1}$ comprehension:

\begin{definition}
    Let $\varphi(x)$ be a formula of arithmetic. The \textit{bounded comprehension principle for $\varphi$} is the sentence $\BCmp\varphi \defeq (\forall k)(\exists X)(\forall n)[ n \in X \leftrightarrow (\varphi(n) \land n < k) ]$.
\end{definition}

$\BCmp\varphi$ asserts that for every $k$, the set $X_k = \{ n < k : \varphi(n) \}$ exists.

\begin{lemma}[{\cite[Thm II.3.9]{simpson_subsystems_2009}}]\label{lem:bdd-s1-comp}
    $\RCA$ proves $\BCmp\Si{1}$, i.e.\ $\BCmp\varphi$ for every $\Si{1}$ formula $\varphi$.
\end{lemma}

\begin{definition}
    $\ACA$ is the subsystem consisting of the basic axioms of first-order arithmetic, $\Ind\D[1]{0}$, and $\Cmp\D[1]{0}$.
\end{definition}

$\WKL$ has a slightly different definition. Recall that (finitary) Cantor space $\CSf$ is the set of all finite binary strings. A \textit{tree} is a set $T \subseteq \CSf$ that is closed under taking initial segments. A \textit{path through $T$} is an infinite binary sequence $\alpha$ such that all initial segments are in $T$.

\begin{definition}
	\textit{Weak K\H onig's lemma} is the statement that every infinite tree $T \subseteq \CSf$ contains a path. $\WKL$ is the subsystem consisting of $\RCA$ plus weak K\H onig's lemma.
\end{definition}

Every infinite tree $T \subseteq \CSf$ has a path which is arithmetical relative to $T$ \cite[Example I.8.8]{simpson_subsystems_2009}. Hence, $\ACA$ implies $\WKL$, which in turn implies $\RCA$. In fact, all these implications are strict. Along with the stronger systems $\ATR$ and $\PiCA$, these make up the ``big five'' subsystems of second-order arithmetic. $\RCA$ has a standard $\omega$-model $\mathrm{REC}$, consisting of the subsets of $\omega$ definable by a $\D{1}$ formula. Similarly, $\ACA$ has a standard $\omega$-model $\mathrm{ARITH}$, consisting of the subsets of $\omega$ definable by an arithmetical formula. We also have the full $\omega$-model $\Pow(\omega)$, which is a model of all these axiom systems.

The idea of reverse mathematics is: for a known theorem $\varphi$ of mathematics, find the weakest formal system $\Ss$ such that $\Ss \proves \varphi$. Given such a proof, we show $\Ss$ is optimal by demonstrating a \textit{reversal} of $\varphi$: a proof $\varphi \proves \Ss$. In practice, no single theorem can axiomatise all of mathematics, and so we have to supplement $\varphi$ with a base theory $\B$ (i.e.\ a reversal is actually a proof $\B + \varphi \proves \Ss$). Throughout this report, we will take $\B = \RCA$, unless stated otherwise.

If we prove $\Ss \proves \varphi$, and then reverse this (over $\B$), we say that $\varphi$ is \textit{equivalent} to $\Ss$ (over $\B$). A remarkable empirical result of reverse mathematics is that almost all theorems of mathematics are equivalent (over $\RCA$) to one of the big five systems, though there are exceptions. In the case of abstract algebra, most results can be proven in $\ACA$, except (ordinal) invariant results, which are usually equivalent to $\ATR$, and structure theorems, which often fall at the level of $\PiCA$. To demonstrate, we now review some reverse-mathematical studies of ring theory in the literature.

\begin{theorem}
	The following are provable in $\RCA$:
	\begin{enumerate}
		\item A field has no nontrivial proper ideals \cite{downey_ideals_2007}.
		
		\item If $I \subseteq R$ is an ideal and $\q{R}{I}$ is a field, then $I$ is maximal \cite{downey_ideals_2007}.
		
		\item An ideal $I \subseteq R$ is prime if and only if $\q{R}{I}$ is an integral domain \cite{conidis_chain_2010}.
		
		\item Maximal ideals are prime \cite{conidis_chain_2010}.
		
		
		
		
		\item Euclidean domains are PIDs \cite{sato_reverse_2016}.
		
		
		
		
		
	\end{enumerate}
\end{theorem}

\begin{theorem}[($\RCA$)]
	$\WKL$ is equivalent to each of the following statements:
	\begin{enumerate}
		\item Every commutative ring has a prime/radical ideal \cite{friedman_countable_1983,friedman_addendum_1985}.
		
		
		
		\item If $I \subseteq R$ is a maximal ideal, then $\q{R}{I}$ is a field \cite{downey_ideals_2007}.
		
		\item In Artinian rings, prime ideals are maximal \cite{conidis_chain_2010}.
		
		\item Every Artinian integral domain is a field \cite{conidis_chain_2010}.
		
		
		
		
		\item $r \in R$ is nilpotent if and only if $r$ belongs to every prime ideal \cite{sato_reverse_2016}.
		
		\item (Local) Artinian rings are Noetherian \cite{conidis_chain_2010, conidis_computability_2019}.
		
		
		\item Every Artinian ring is a finite direct product of local Artinian rings \cite{conidis_computability_2019}.
		
		
		
	\end{enumerate}
\end{theorem}



\begin{theorem}[($\RCA$)]
	$\ACA$ is equivalent to each of the following statements:
	\begin{enumerate}
		\item Every commutative ring/domain has a maximal ideal \cite{friedman_countable_1983}.
		
		
		\item Every commutative ring has a minimal prime ideal \cite{hatzikiriakou_commutative_1989}.
		
		\item If $R$ has no nontrivial, proper, principal ideals, it is a field \cite{downey_ideals_2007}.
		
		
		\item An integral domain is a UFD iff it has the a.c.c.p.\ and all irreducibles are prime \cite{greenberg_proper_2017}.
		
		
		\item $r \in R$ belongs to every maximal ideal of $R$ if and only if for all $s \in R$, $1-rs$ is a unit \cite{sato_reverse_2016}.
	\end{enumerate}
\end{theorem}

\section{Computability}

Here, we review the basic notions of computability theory that we will need. The presentation will be very brief, so the reader unfamiliar with computability theory is urged to consult a textbook on the subject \cite{soare_recursively_1987,soare_turing_2016}.

We will assume the Church--Turing thesis, and work with an informal notion of computability. Therefore, we say a set $A \subseteq \omega$ is \textit{computable} if there is an algorithm which, given $n$, always terminates and tells us if $n \in A$. Similarly, a function $f\colon \omega \to \omega$ is \textit{(total) computable} if there is an algorithm which, given $n$, always halts and outputs $f(n)$.

By this definition, most sets and functions arising in mathematics are computable. We can construct non-computable sets: the archetypal example is the halting problem $\zerojump$, which is the set of pairs $\ang{e,n}$ for which the $e$th computable function halts on input $n$. Nonetheless, $\zerojump$ is still \textit{computably enumerable}, or c.e.: there is an algorithm which lists its elements. There are several equivalent definitions of c.e.:

\begin{proposition}
    The following are equivalent for a set $A \subseteq \omega$:
    \begin{enumerate}
        \item $A$ is c.e., i.e.\ there is an algorithm which lists the elements of $A$.
        
        \item $A$ is the domain of a partial computable function.
        
        \item $A$ is the range of a partial computable function.
        
        \item $A$ is empty, or the range of a total computable function.
        
        \item $A$ is finite, or the range of a total computable injection.
        
        \item $A$ has a \textit{computable enumeration}: a uniformly computable sequence $A_0 \subseteq A_1 \subseteq A_2 \subseteq \cdots$ such that $A = \bigcup_{s \in \omega} A_s$. We can further require that $\abs{A_s} = s$.
    \end{enumerate}
\end{proposition}

There is a close correspondence between computability and the arithmetical hierarchy: a set $A \subseteq \omega$ is computable if and only if it can be defined by a $\D{1}$ formula, and c.e.\ if and only if it can be defined by a $\Si{1}$ formula.

We recall some of the basic results of computability theory. Every algorithm can be coded by a natural number: for example, by writing it in a fixed (Turing-complete) programming language. This gives a listing $\varphi_e$ of the partial computable functions. We can furthermore get a \textit{uniformly computable} listing, i.e.\ the function $\ang{e,n} \mapsto \varphi_e(n)$ is computable. Similarly, we can get a uniformly c.e.\ listing $W_e$ of the c.e.\ sets.

We can also consider \textit{relativised computations}: those with access to an \textit{oracle} $A \subseteq \omega$, so that the algorithm can ask at any point if any natural number $n$ is in $A$. We say $A$ is \textit{Turing reducible} to $B$, and write $A \leq_\mathrm{T} B$, if $A$ can be computed with $B$ as an oracle. The relation $\leq_\mathrm{T}$ is a preorder, so we obtain an equivalence relation $\equiv_\mathrm{T}$ in the standard way: $$A \equiv_\mathrm{T} B \iff A \leq_\mathrm{T} B\ \text{ and }\ B \leq_\mathrm{T} A$$ The $\equiv_\mathrm{T}$-equivalence-classes are called \textit{Turing degrees}, and are partially ordered by $\leq_\mathrm{T}$. The degree of all computable sets is called $\mathbf{0}$.

For any oracle $A \subseteq \omega$, we can similarly list all the partial $A$-computable functions $\varphi_e^A$. Thus, for any $A$, we define the \textit{Turing jump} $A'$ as the halting problem relativised to $A$, i.e. $$A \defeq \{ \ang{e,n} : \varphi_e^A(n) \text{ halts} \}$$
Then, $A <_\mathrm{T} A'$, i.e.\ $A'$ is always strictly above $A$ in the Turing degrees. The Turing jump is also a well-defined function on Turing degrees.

We also have a stronger notion of computable reduction: we say $A$ is \textit{m-reducible} to $B$ ($A \leq_m B$) if there is a computable function $f\colon \omega \to \omega$ such that $n \in A \iff f(n) \in B$. Informally, the idea is that we can consult the oracle $B$ only \textit{once}. For some level $\Gamma$ in the arithmetical hierarchy, we say that a $\Gamma$ set $A$ is \textit{$\Gamma$ complete} if for every $\Gamma$ set $B$, $B \leq_m A$.

\begin{proposition}\
    \begin{enumerate}
        \item $\zerojump$ is $\Si{1}$ complete.
        
        \item $\Inf = \{ e: W_e \text{ is infinite} \}$ is $\Ppi{2}$ complete.
    \end{enumerate}
\end{proposition}

Another important class of Turing degrees are the PA degrees, which are those that can compute a complete consistent extension of Peano arithmetic. By G\"odel's famous incompleteness theorem, $\mathbf{0}$ is not PA. $\mathbf{0}'$ is PA, but there are also PA degrees strictly below $\mathbf{0}'$. There are also many equivalent characterisations of PA degrees:

\begin{proposition}\label{prop:pa-deg}
    For a Turing degree $\mathbf{d}$, the following are equivalent:
    \begin{enumerate}
        \item $\mathbf{d}$ is PA, i.e.\ it computes a complete consistent extension of $\PA$.
        
        \item $\mathbf{d}$ computes a path through any computable (or $\Ppi{1}$) tree $T \subseteq \CSf$.
        
        \item For any disjoint c.e.\ sets $A,B$, $\mathbf{d}$ computes a \textit{separating set} for $A$ and $B$, i.e.\ a set $C$ such that $A \subseteq C$, $B \cap C = \varnothing$.
    \end{enumerate}
\end{proposition}

We can also speak about PA degrees relative to some oracle $A$, or PA degrees \textit{over $A$}. Proposition \ref{prop:pa-deg} relativises: $\mathbf{d}$ is PA over $A$ iff it computes an extension of ``$\PA$ plus a predicate for $A$'', iff it computes a path through any $A$-computable tree, iff it can separate any two disjoint $A$-c.e.\ sets.

Having reviewed classical computability, we now point out the close correspondence between computability and the subsystems $\RCA$, $\WKL$, $\ACA$ of second-order arithmetic. This connection arises because our chosen base system $\RCA$ somehow corresponds to the ``computable world''. Indeed, $\D{1}$ comprehension allows us to define \textit{exactly} the subsets of $\omega$ which are computable (with parameters/oracles in the model).

Many of the theorems $T$ studied in reverse mathematics have the form $$T = (\forall X) \big[ \varphi(X) \to (\exists Y)\, \psi(X,Y) \big]$$ where $\varphi(X)$ and $\psi(X,Y)$ are properties of the sets $X$ and $Y$. We could view $T$ as a problem or challenge: given a set $X$ such that $\varphi(X)$, find a set $Y$ such that $\psi(X,Y)$. Now the connection is thus: $\RCA \proves T$ when we can always choose $Y$ to be $X$-computable. $\WKL \proves T$ when we can choose $Y$ computable from a PA degree over $X$, and $\ACA \proves T$ if we can choose $Y$ arithmetical in $X$.

This correspondence is also useful in reversals of $T$. Over $\RCA$, $T \proves \WKL$ if for every $A$, we can construct $A$-computable $X$ so that every suitable $Y$ is PA over $A$. Similarly, $T \proves \ACA$ if we can construct $X$ so that every $Y$ computes $A'$.

While this correspondence is completely precise only in $\omega$-models, a nearly identical proof can generally be made to work in nonstandard models. Furthermore, we will generally just take $X$ to be computable (i.e.\ $A = \varnothing$), as the fully relativised version will again follow nearly identically. So in practice, we will prove something like ``there is a computable $X$ with $\varphi(X)$, such that every $Y$ with $\varphi(X,Y)$ has PA degree [computes $\zerojump$]'', and we will take this as evidence that $T \proves \WKL$ [resp.\ $T \proves \ACA$].

Now, we will construct various computability-theoretic objects, which will be useful later when proving reversals. We saw that a degree $\mathbf{d}$ is PA iff it can separate any two disjoint c.e.\ sets $A,B$. Our first lemma is that there is a ``universal pair'' $A,B$ so that \textit{any} separator has PA degree:

\begin{lemma}\label{lem:pa-sep}
    There are disjoint c.e.\ sets $A$, $B$ such that whenever $C$ has $A \subseteq C$, $B \cap C = \varnothing$, then $C$ has PA degree.
\end{lemma}

\begin{proof}
    Fix a coding $F_n$ of all sentences in the language of $\PA$. Let $A = \{ n : \PA \proves F_n \}$ and $B = \{ n : \PA \proves \lnot F_n \}$. Then, $A,B$ are c.e.\ since $\PA$ is computably axiomatisable, so we can search for a proof of $F_n$ from $\PA$. $A \cap B = \varnothing$ since $\PA$ is consistent.
    
    We claim $A,B$ are as required. Fix a separator $C$ of $A$ and $B$. There is no reason $C$ should be consistent, but we can use it to compute a complete consistent extension $D$, as follows. We build $D$ in stages $D_s$, with $D_0 = \varnothing$. At stage $s$, let $G_s \defeq (\bigwedge_{F \in D_s} F) \to F_s$. If $G_s \in C$, let $D_{s+1} = D_s \cup \{ F_s \}$; else, let $D_{s+1} = D_s \cup \{ \lnot F_s \}$. By construction, $D$ is a complete consistent extension of $\PA$, hence $C$ is PA.
\end{proof}

Lemma \ref{lem:pa-sep} will be useful in showing a theorem implies $\WKL$. To show a theorem implies $\ACA$, we will employ a few different techniques. The first is simply to code $\zerojump$ itself into a computable ring. Often it is not possible to do this directly, so we have a few other methods.

\begin{definition}
    Fix a computable enumeration $\varnothing_s'$ of $\zerojump$. The \textit{modulus of $\zerojump$} is the function $\mu_\zerojump$ mapping $n$ to the least $s$ such that $\substr{\varnothing_s'}{n} = \substr{\zerojump}{n}$.
\end{definition}

\begin{theorem}
    Suppose $f\colon \omega \to \omega$ \textit{dominates} $\mu_\zerojump$: for every $n$, $f(n) \geq \mu_\zerojump(n)$. Then, $f \geq_\mathrm{T} \zerojump$.
\end{theorem}

\begin{proof}
    We compute if $n \in \zerojump$ by asking if $n \in \varnothing_{f(n+1)}'$. By definition of $\mu_\zerojump$, this computation is always correct.
\end{proof}

Another useful trick to code $\varnothing'$ into a construction is using a c.e.\ set that is so ``dense'' that its complement dominates $\mu_\varnothing'$.

\begin{lemma}\label{lem:ce-dense}
    There is a c.e.\ set $A$ such that
        \begin{enumerate}
        \item $A^\complement$ is infinite.
        
        \item Any infinite subset $B \subseteq A^\complement$ computes $\varnothing'$.
        
    \end{enumerate}
\end{lemma}

\begin{proof}
    Fix an enumeration of $\zerojump$, and enumerate $A$ starting with $A_0 = \varnothing$. Now, at stage $s$, suppose $A_s^\complement = \{ c_0 < c_1 < \cdots \}$. If $n$ enters $\zerojump$ at stage $s$, we put $c_n, c_{n+1}, \ldots, c_s$ into $A$. We now verify that $A$ has the required properties.
%
    \begin{enumerate}
        \item $A^\complement$ is infinite: by induction, all the $A_s$ are finite, so each $A_s^\complement$ is infinite. Now given $k \in \N$, $A_{\mu_\zerojump(k+1)}^\complement$ contains some $j > k$ which will never leave $A^\complement$.
        
        \item Any infinite subset $B \subseteq A^\complement$ computes $\zerojump$ in a very simple way: if $B = \{ b_0 < b_1 < \cdots \}$, then the function $n \mapsto b_n$ dominates $\mu_{\varnothing'}$. The proof is as follows. If $A^\complement = \{ c_0 < c_1 < \cdots \}$, then each $b_i = c_j$ for some $j$ with $i \leq j \leq b_i$.
        
        Now if any $k \leq i$ entered $\zerojump$ at a stage $s \geq b_i$, then $c_k,\ldots,c_i,\ldots,c_j=b_i,\ldots,c_{b_i}$ would have been put into $A$ (by construction). As $b_i \notin A$, this can't have happened, so $\substr{\varnothing_{b_i}}{i} = \substr{\varnothing}{i}$, i.e.\ $b_i \geq \mu_{\varnothing'}(i)$.\qedhere
%
        
    \end{enumerate}

\end{proof}

To finish this section, we will review some useful ideas from a sub-branch of computability called \textit{computable structure theory} \cite{montalban_computable_2021}. Essentially, this is the study of model theory from the point of view of computability, and the primary objects of study are computable structures. For what follows, let $\T$ be a computable theory in a finite language $\Ll$.

\begin{definition}
    A \textit{computable (presentation of a) $\T$-structure} $\M$ consists of:
    \begin{enumerate}
        \item A computable subset $M \subseteq \N$.
        
        \item For every constant $c \in \Ll$, an element $c^\M \in M$.
        
        \item For every $k$-ary function $f \in \Ll$, a computable function $f^\M\colon M^k \to M$.
        
        \item For every $k$-ary relation $R \in \Ll$, a computable function $R^\M\colon M^k \to \{ 0, 1 \}$.
    \end{enumerate}
    such that, with these interpretations of the $\Ll$-symbols, $\M$ satisfies all formulae in $\T$.
    
    A \textit{c.e.\ $\T$-structure} consists of the same data, but we allow $M$ to be c.e.\ instead.
\end{definition}

We will be primarily concerned with the language $\Ll_\text{ring} = \{ 0, 1, +, \cdot \}$, and the $\Ll_\text{ring}$-theory $\T_\text{ring}$ consisting of the usual commutative ring axioms. A \textit{computable ring} is just a computable $\T_\text{ring}$-structure, i.e.\ a computable set $R \subseteq \N$ with elements $0_R, 1_R \in R$ and computable binary operations $+, \cdot$ on $R$ which form a ring.

When doing reversals, we will often need to construct a computable ring having certain properties. However, the following result shows that it is sufficient to construct a c.e.\ ring:
\begin{theorem}\label{thm:ce-struct}
    Any c.e.\ $\T$-structure $\M$ is (computably) isomorphic to a computable $\T$-structure $\M'$.
\end{theorem}

\begin{proof}
    If $\M$ is finite, this is trivial, so assume $\M$ is infinite, and fix an injective computable enumeration $\varphi\colon \N \to \M$. The inverse $\varphi^{-1}\colon \M \to \N$ is partial computable, since given $m \in \M$, we can search for the $n$ such that $\varphi(n) = m$.
    
    Now, define $\M'$ as the structure with support $\N$, and $\Ll$-symbols interpreted
    \begin{align*}
        c^{\M'} &\defeq \varphi^{-1}(c^\M) \\ f^{\M'}(n_1,\cdots,n_k) &\defeq \varphi^{-1} \big[ f^\M \big( \varphi(n_1), \ldots, \varphi(n_k) \big) \big] \\
        R^{\M'}(n_1,\cdots,n_k) &\defeq R^\M \big( \varphi(n_1), \ldots, \varphi(n_k) \big)
    \end{align*}
    Then, $\M'$ is a computable $\T$-structure, and by construction, $\varphi$ is a computable $\T$-isomorphism $\M' \to \M$.
\end{proof}

When using Theorem \ref{thm:ce-struct} in practice, we will often abuse notation, and identify $\M$ and $\M'$ via the computable isomorphism $\varphi$. This means we will treat c.e.\ rings as if they are computable. We should note that the isomorphism will also preserve many computable subobjects of $\M$. For example, if $I \subseteq \M$ is a computable ideal, then so is $\varphi^{-1}(I) \subseteq \M'$.

A particularly useful case of Theorem \ref{thm:ce-struct} is localising a computable ring by a c.e.\ subset:

\begin{corollary}\label{cor:ce-loc}
    Suppose $R$ is a computable ring, and $M \subseteq R$ is a multiplicatively closed c.e.\ subset containing $1_R$ but no zero-divisors. Then $\loc{R}{M}$ is (computably) isomorphic to a computable ring $S$.
\end{corollary}

\begin{proof}
    Define $\loc{R}{M}$ exactly as in Definition \ref{defn:loc}. This is a computable quotient of a c.e.\ structure, hence we get a c.e.\ structure. Now apply Theorem \ref{thm:ce-struct}.
\end{proof}
\chapter{Radicals of ideals and primary ideals}
\label{sec:prim-rad}

We now begin our study of ring theory in second-order arithmetic. We can use the standard definitions of rings, polynomial rings, ideals, etc.\ verbatim in second-order arithmetic. The definition of quotient ring might cause some concern, since we naively define $R/I$ as a set of sets. However, we can amend this by defining the elements of $R/I$ to be minimal representatives of their equivalence class \cite[Defn III.5.2]{simpson_subsystems_2009}.

\begin{definition}[($\RCA$)]\label{defn:quo-ring-rca}
    Let $I \subseteq R$ be an ideal. The \textit{quotient ring} $R/I$ is the set $$R/I\ \defeq\ \{ r \in R : (\forall s <_\N r)(s \in R\, \to\, r-s \notin I) \}$$
    which exists by $\D{1}$ comprehension. $\RCA$ can define a function $q\colon R \to R/I$, called the \textit{quotient map}, so that for every $r \in R$, $q(r)$ is the unique element of $R/I$ such that $r - q(r) \in I$. The ring operations on $R/I$ are then the operations induced by $q$ from $R$. 
\end{definition}

A similar trick can be used to construct the field of fractions or localisations in $\RCA$. As usual, we use $r+I$ (or simply $r$) to denote $q(r)$. In $\RCA$, it follows from Definition \ref{defn:quo-ring-rca} that $r+I = s+I$ iff $r-s \in I$. We will also make wide use of the following theorem:

\begin{theorem}[(ideal correspondence theorem; $\RCA$)]\label{thm:id-corr}
    For a ring $R$ and ideal $I \subseteq R$, the quotient map $q\colon R \to R/I$ is an isomorphism between the ideals of $R$ containing $I$, and the ideals of $R/I$.
\end{theorem}

Furthermore, the quotient map preserves many properties of ideals, such as maximality, primality, primary-ness, being the radical of another ideal, etc. As a result, this frequently gives an equivalence between statements of the form $\forall R\ \forall I\! \subseteq\! R\ \varphi(R,I)$ and $\forall R\ \varphi(R,\{0\})$, with the equivalence provable in $\RCA$.

By fixing the parameters $a$ and $b$, the following lemmas can be proved by $\D{0}$ induction.

\begin{lemma}[($\RCA$)]
    For all $a,b$ in a commutative ring $R$, $(ab)^n = a^n b^n$.
\end{lemma}

\begin{lemma}[(binomial theorem; $\RCA$)]
    For all $a,b$ in a commutative ring $R$, $(a+b)^n = \sum_{j=0}^n \binom{n}{j} a^j b^{n-j}$.
\end{lemma}

We now examine some properties of ideals and the relationships between them, and show that most of these relationships are provable in $\RCA$.

\begin{definition}
    A ring $R$ is \textit{reduced} if it has no nontrivial nilpotent elements, i.e.\ whenever $x^n=0$ for some $n>0$, then $x=0$.
\end{definition}

\begin{definition}
    Let $I \subseteq R$ be an ideal.
    \begin{enumerate}
        \item $I$ is \textit{prime} if whenever $ab \in I$, then $a \in I$ or $b \in I$.
        
        \item $I$ is \textit{primary} if whenever $ab \in I$, then $a \in I$ or $b^n \in I$ for some $n$.
        
        \item The \textit{radical of $I$} is defined $\sqrt{I} = \{ a \in R: \exists\, n>0 \text{ such that } a^n \in I \}$.
        
        \item $I$ is \textit{semiprime} or \textit{radical} if whenever $a^n \in I$ for some $n>0$, then $a \in I$ (equivalently, $I = \sqrt{I}$).
        
        \item The \textit{adjoint of $I$} is defined $I^{\not\perp} = \{ a \in R: \exists\, b \notin I \text{ such that } ab \in I \}$.
        
        \item $I$ is \textit{primal} if $I^{\not\perp}$ forms an ideal.
    \end{enumerate}
\end{definition}

\begin{theorem}\label{thm:radicals-rca}
    The following are provable in $\RCA$:
    \begin{enumerate}
        \item $\sqrt{I}$ is an ideal.\label{thm:radical-is-ideal}
        
        \item $I \subseteq R$ is semiprime if and only if $R/I$ is reduced.
        
        \item If $I$ is primary, then $\sqrt{I}$ is prime\footnote{We also say that $I$ is \textit{quasi-primary}.}.
        
        \item $I$ is prime iff it is primary and semiprime.
        
        \item $I$ is primary iff every zero divisor in $R/I$ is nilpotent.\label{thm:primary-nilp-zd}
        
        \item If $I \neq R$ is primary, then $I$ is primal (and furthermore, $I^{\not\perp} = \sqrt{I}$).
        
        \item If $I$ is primal, then $I^{\not\perp}$ is prime.
    \end{enumerate}
\end{theorem}

\begin{proof}\
    \begin{enumerate}
        \item We prove only closure under addition. Suppose $a,b \in \sqrt{I}$, i.e.\ $a^n, b^m \in I$. Let $k \defeq n+m-1$. Using the binomial formula, we can write $(a+b)^k$ as
        \begin{center}
        \arraycolsep=0pt
        $\begin{array}{lrlrl}
            \multicolumn{2}{c}{a^k\ +\ \cdots\ +\ \binom{k}{n}\, a^n\, b^{m-1}} & \multicolumn{2}{c}{\ \ +\ \ \binom{k}{m}\, a^{n-1}\, b^{m}\ +\ \cdots\ +\ b^k} \\[2mm]
            a^n\big( & \big) & \ \ +\ \ b^m \big( & \big) &\ \ \in\ I \\
        \end{array}$
        \end{center}
        hence $a+b \in \sqrt{I}$.
        
        \item By the ideal correspondence theorem, this is equivalent to saying $\{0\}$ is semiprime iff $R$ is reduced, which is trivial.
%
        
        \item Suppose $ab \in \sqrt{I}$: then $(ab)^n = a^n b^n \in I$, so either $a^n \in I$, whence $a \in \sqrt{I}$, or $(b^n)^k = b^{nk} \in I$, whence $b \in \sqrt{I}$.
        
        \item \begin{itemize}
            \ifff Any prime ideal is trivially primary (take $n=1$), and semiprime by $\D{0}$ induction on $\varphi(n): a^n \in I \to a \in I$ for a fixed parameter $a \in I$.
            
            \iffb If $ab \in I$, then by primary-ness, either $a \in I$ or $b^n \in I$, whence $b \in I$ by semiprimality.
        \end{itemize}
        
        \item By the ideal correspondence theorem, this is equivalent to saying $\{ 0 \}$ is primary iff every zero divisor in $R$ is nilpotent, which is immediate from the definition of primary.
        
%
        
        \item By \ref{thm:radical-is-ideal}, it is enough to prove the assertion in brackets. Pick $a \in I^{\not\perp}$, i.e.\ there is $b \notin I$ with $ab \in I$. Since $I$ primary, either $b \in I$ (which is a contradiction), or $a \in \sqrt{I}$ as required. Now if $a \in \sqrt{I}$, by $\LSi{1}$, there is a minimal $n>0$ such that $a^n \in I$. Then $a^{n-1}$ witnesses that $a \in I^{\not\perp}$.
        
        \item Suppose $ab \in I^{\not\perp}$, i.e.\ there is $c \notin I$ such that $abc \in I$. If $bc \in I$, then $c$ witnesses that $b \in I^{\not\perp}$. If $bc \notin I$, then $bc$ witnesses that $a \in I^{\not\perp}$.\qedhere
    \end{enumerate}
\end{proof}

\noindent Now, we present some basic results about operations and relations between \textit{two} ideals and their radicals, which are provable in $\RCA$.

\begin{definition}
    Ideals $I,J \subseteq R$ are \textit{comaximal} if every $r \in R$ can be written as $r = i+j$ for $i \in I$, $j \in J$.
\end{definition}

It suffices (under $\RCA$) to show that $1=i+j$.

\begin{theorem}
    The following are provable in $\RCA$:
    \begin{enumerate}
        \item If $\sqrt{I}$ and $\sqrt{J}$ are comaximal, then $I$ and $J$ are comaximal.
        
        \item $\sqrt{I \cap J} = \sqrt{I} \cap \sqrt{J}$.
    \end{enumerate}
\end{theorem}

\begin{proof}\
    \begin{enumerate}
        \item Suppose $1 = a+b$ where $a^n \in I$, $b^m \in J$. Let $k \defeq n+m-1$. Then, using the binomial formula, we can write $1 = 1^k = (a+b)^k$ as
        \begin{center}
        \arraycolsep=0pt
        $\begin{array}{lrlr}
            \multicolumn{2}{c}{a^k\ +\ \cdots\ +\ \binom{k}{n}\, a^n\, b^{m-1}} & \multicolumn{2}{c}{\ \ +\ \ \binom{k}{m}\, a^{n-1}\, b^{m}\ +\ \cdots\ +\ b^k} \\[2mm]
            a^n\big( & \big) & \ \ +\ \ b^m \big( & \big) \\
        \end{array}$
        \end{center}
        Thus, $1 = i+j$ for suitable $a^n \mid i$, $b^m \mid j$.
        
        \item \begin{enumerate}
            \item[($\subseteq$)] Trivial.
            \item[($\supseteq$)] If $a \in \sqrt{I} \cap \sqrt{J}$, then $a^n \in I$, $a^m \in J$, so $a^{\max\{n,m\}} \in I \cap J$.\qedhere
        \end{enumerate}
    \end{enumerate}
\end{proof}

\noindent To conclude this section, we analyse the following characterisation of the radical: $$\sqrt{I}\ = \bigcap_{\substack{P \supseteq I \\
P \text{ prime}}} P$$

However, we require $\ACA$ to show the LHS exists, and $\PiCA$ to show the RHS exists. So, we will analyse the part of this theorem that doesn't require comprehension for either side. The two containments can be written respectively as:
\begin{itemize}
    \item $I$ and $\sqrt{I}$ are contained in exactly the same prime ideals.
    
    \item If $x \notin \sqrt{I}$, then there is a prime ideal $P \supseteq I$ such that $x \notin P$.
\end{itemize}

\noindent By the ideal correspondence theorem, these are respectively equivalent to:
\begin{itemize}
    \item If $r$ is nilpotent, it belongs to every prime ideal of $R$.
    
    \item If $r$ belongs to every prime ideal of $R$, then it is nilpotent.
\end{itemize}

\cite[62]{sato_reverse_2016} showed that the first statement is provable in $\RCA$, and the second is equivalent to $\WKL$.

\section{The \texorpdfstring{$\RAD$}{RAD} principle}
\label{sec:rad}

So far, we have seen that almost all basic facts about radicals and primary ideals are provable in $\RCA$. The notable exception is the following:

\begin{theorem}[($\RCA$)]\label{thm:wkl-maximal-primary}
    $\WKL$ is equivalent to ``for all ideals $I \subseteq R$, if $\sqrt{I}$ is maximal, then $I$ is primary''.
\end{theorem}

\begin{proof}\
    \begin{enumerate}
        \ifff We prove the contrapositive of the consequent. Suppose $I$ is not primary: then there are $a,b \in R$ such that $ab \in I$, $a \notin I$, $b \notin \sqrt{I}$. Using $\WKL$, we will construct an ideal $J \supsetneq \sqrt{I}$.
        
        Define $T \subseteq \CSf$ as the set of all $F \subseteq \{ 0,1,\ldots,n-1 \}$ such that the following hold:\footnote{Here, we identify binary strings of length $n$ with subsets $F \subseteq \{ 0,1,\ldots,n-1 \}$.}
        \begin{enumerate}[label=(\roman*)]
            \item If $0_R < n$, then $0_R \in F$.
            
            \item For all $c,d \in F$, if $c+d < n$, then $c+d \in F$.
            
            \item For all $c \in F$, $d < n$, if $cd < n$, then $cd \in F$.
            
            \item For all $c<n$, if $c^n \in I$, then $c \in F$.
            
            \item If $b < n$, then $b \in F$.
            
            \item If $1_R < n$, then $1_R \notin F$.
        \end{enumerate}
        $T$ is computable since for a finite string, we can simply check all the above conditions exhaustively. A path $J$ through $T$ is a proper ideal containing $b$ and $\sqrt{I}$: since $b \notin \sqrt{I}$, it follows that $J \supsetneq \sqrt{I}$. Also, $T$ is downwards closed, hence a tree.
        
        $T$ is infinite since every level is nonempty. Given $n$, 
        let
        \begin{equation*}
            A_n \defeq \begin{cases}
                \{ c < n: c^n \in I \} \cup \{ b \}\quad & b < n \\
                \{ c < n: c^n \in I \} & b \geq n \\
            \end{cases}
        \end{equation*}
        Enumerating $A_n = \{ c_1, \ldots, c_k \}$, let $$F_n\ \defeq\ \{ c < n: (\exists\, d_1,\ldots,d_k \in R)(c = c_1 d_1 + \cdots + c_k d_k) \}\ =\ \substr{(A_n)}{n}$$
        which exists by bounded $\Si{1}$ comprehension (Lemma \ref{lem:bdd-s1-comp}).
        
        We claim $F_n \in T$ at level $n$. Conditions (i)--(v) are evidently satisfied. For condition (vi), note that it can only fail if $1_R < n$ and $1_R \in F$, i.e.\ $1_R = c_1 d_1 + \cdots + c_k d_k$ for some $d_1,\ldots,d_k \in R$. Expanding $$1_R = (1_R)^{nk} = (c_1 d_1 + \cdots + c_k d_k)^{nk}$$ note that each term is divisible by either $b$ or some $c^n \in I$, i.e.\ $1_R \in (I,b)$. Then $a = a 1_R \in (I,ab) = I$, giving a contradiction.\qed

        
        
        

        \iffb We construct a computable ring $R$ and a computable, proper ideal $I \subsetneq R$ which is not primary, such that any proper ideal $J \supsetneq \sqrt{I}$ has PA degree. Thus, the assumption that $\sqrt{I}$ is not maximal gives a larger ideal $J \supsetneq \sqrt{I}$ of PA degree, which computes $\WKL$.
        
        Fix disjoint c.e.\ sets $A, B$ as in Lemma \ref{lem:pa-sep}. We will build $R$ and $I$ such that any proper ideal $J \supsetneq \sqrt{I}$ computes a separator of $A$ and $B$.
    
        Let $R = \Z[x,y,z_0, z_1,\ldots]$. To begin, let $I = ( xy )$ (the ideal generated by $xy$). This is to ensure $I$ is not primary (with witness $xy$). Now if $j$ enters $A$ at stage $s$, add $z_j^s$ to the list of generators for $I$. If $j$ enters $B$ at stage $s$, add $(z_j - 1)^s$ to the list of generators for $I$.

        Then, $I$ is computable - to work out if $p \in I$, we only need to run the construction up to stage $s = \degr(p)$. At this stage, $I$ is finitely generated, so we can compute if $p \in I$. We can assume nothing is enumerated when $s=0$, meaning $I$ (and $\sqrt{I}$) are proper ideals.
        
        Now, let $J \supseteq \sqrt{I}$ be a proper ideal of $R$. Let $X_J = \{ n \in \N: z_n \in J \}$. Then, $X_J$ separates $A$ and $B$: $A \subseteq X_J$, and $B \cap X_J = \varnothing$, because if there were $n \in B \cap X_J$, then $z_n, (z_n - 1) \in J$ so $1 \in J$, contradicting that $J$ is proper. Since $X_J$ is $J$-computable, it follows that $J$ has PA degree.\qedhere
    \end{enumerate}
\end{proof}

In the ($\Leftarrow$) construction of Theorem \ref{thm:wkl-maximal-primary}, note that $\sqrt{I}$ itself is of PA degree. Thus, the computational power lies in comprehension for an ideal containing $\sqrt{I}$, possibly not strictly. By asserting that $\sqrt{I}$ itself must be computable, the focus shifts more to the (non-)primary-ness of $I$. In other words, we want to establish the reverse-mathematical strength of the statement
\begin{equation}\label{rad-principle}
    \text{``if $\sqrt{I}$ \textbf{exists} and is maximal, then $I$ is primary''}\tag{$\RAD$}
\end{equation}
or equivalently,
\begin{equation*}
    \text{``if $I$ is not primary and $\sqrt{I}$ \textbf{exists}, then $\sqrt{I}$ is not maximal''}
\end{equation*}
By the ideal correspondence theorem, $\RAD$ is also equivalent to:
\begin{gather*}
    \text{``if $\sqrt{0}$ \textbf{exists} and is maximal, then $\{ 0 \}$ is primary''} \\
    \text{``if $\{ 0 \}$ is not primary and $\sqrt{0}$ \textbf{exists}, then $\sqrt{0}$ is not maximal''}
\end{gather*}

\noindent From Theorem \ref{thm:wkl-maximal-primary}, it follows that:

\begin{proposition}
    $\WKL$ proves $\RAD$.
\end{proposition}

Conversely, we might try to argue that $\RAD$ implies $\WKL$. One strategy is to use a construction similar to \cite[Thm 3.2]{downey_ideals_2007}. As before, we fix disjoint, c.e.\ sets $A, B \subseteq \N$ such that any separator is of PA degree. Starting with some ``ring of coefficients'' $U$, we set $R_0 \defeq U[x_0,x_1,\ldots]$. Then, we enumerate $R$ as a c.e.\ subring of the total quotient ring of $R_0$, adding $x_n / f$ for $n \in A$, $f \in U[x_0,\ldots,x_{n-1}]$ a non-zero-divisor, and $(x_n-1) / f$ for $n \in B$, $f \in U[x_0,\ldots,x_{n-1}]$ a non-zero-divisor.

We want to choose $U$ so that $R$ satisfies the following:

\begin{enumerate}[label=(\alph*)]
    \item $I = \{0\}$ is not primary (equivalently, there are non-nilpotent zero divisors).
    
    \item The nilradical $\sqrt{0}$ is computable.
    
    \item Any proper ideal $J \supsetneq \sqrt{0}$ has PA degree.
\end{enumerate}
        
\noindent To ensure $R$ satisfies condition (c), we will force it to satisfy:
\begin{enumerate}[label=(\alph*)]\setcounter{enumi}{3}
    \item Any proper ideal $J \supsetneq \sqrt{0}$ contains a non-zero-divisor.
\end{enumerate}
Then, given an ideal $J \supsetneq \sqrt{0}$, we have a non-zero-divisor $f \in J$, so $x_n = f \cdot (x_n/f) \in J$ for sufficiently large $n \in A$, and similarly $x_n \notin J$ for sufficiently large $n \in B$. Thus, $J$ computes a separator for $A$ and $B$ up to finitely many differences. Sensible though it may seem, this strategy is doomed to fail because of the following:

\begin{proposition}
    $\{ 0 \} \subseteq R$ is primary iff any proper ideal $J \supsetneq \sqrt{0}$ contains a non-zero-divisor.
\end{proposition}

\begin{proofcite}{guest123456789_answer_2021}\
    \begin{enumerate}
        \ifff By Theorem \ref{thm:radicals-rca}.\ref{thm:primary-nilp-zd}, every zero-divisor is nilpotent, so since $J$ contains a non-nilpotent, it contains a non-zero-divisor.
        
        \iffb Assuming any proper ideal $J \supsetneq \sqrt{0}$ contains a non-zero-divisor, we will show any non-nilpotent is a non-zero-divisor. Let $x$ be non-nilpotent: then the ideal $\big( x,\sqrt{0} \big) \supsetneq \sqrt{0}$ contains a non-zero-divisor by assumption. In other words, there are $a, b \in R$ so that $ax + b$ is a non-zero-divisor, and $b^n = 0$.
        
        By induction, we show that for every $d>0$, $(ax+b)^d$ is not a zero-divisor. The base case $d=1$ is true by assumption. Now, suppose $(ax+b)^d$ is not a zero-divisor, but $(ax+b)^{d+1}$ is. Let $r \neq 0$ be such that $(ax+b)^{d+1} r = 0$. Then, $(ax+b)^d (ax+b) r = 0$. By assumption, $ax+b$ is not a zero-divisor, so $(ax+b) r \neq 0$: this shows that $(ax+b)^d$ is a zero-divisor, a contradiction.

        In particular $(ax + b)^n$ is a non-zero-divisor. Expanding $(ax + b)^n$ using the binomial formula, and using the fact that $b^n = 0$, we have $(ax + b)^n = cx$ for some $c \neq 0$. Since $cx$ is not a zero-divisor, it follows that $x$ is not a zero-divisor, as required.\qedhere
    \end{enumerate}
\end{proofcite}

\comment{
\noindent With that attempt unsuccessful, we now try to establish a lower bound for $\RAD$. First, we show that it does not follow from $\RCA$.

\begin{proposition}
    $\RCA$ does not imply $\RAD$.
\end{proposition}

\begin{proof}
    We show that $\RAD$ fails in the standard $\omega$-model $\REC$ of $\RCA$, consisting of the computable subsets of $\N$. To do so, we build a computable ring $R$ for which:
    \begin{enumerate}
        \item $\{ 0 \}$ is not primary.
        
        \item $\sqrt{0}$ is computable.
        
        \item There is no computable proper ideal $J$ strictly containing $\sqrt{0}$.
    \end{enumerate}
    
    Let $R_0 \defeq \Z[a,b,x_i : i \in \N]/(ab)$. Then, $\{ 0 \}$ is not primary, since $ab=0$ but $a \neq 0$, $b^n \neq 0$. Note that 
    the nilradical of $R_0$ is $\{ 0 \}$, and the $x_i$ are not zero-divisors. We will construct $R$ from $R_0$, preserving the property that $\sqrt{0} = \{ 0 \}$ is not primary, while ensuring that $R$ satisfies property (iii).
    
    We diagonalise against a listing $J_0, J_1, \ldots$ of all partial computable functions $R \to 2$. Wait for $J_e$ to give an opinion on $x_e$:
    \begin{itemize}
        \item If $J_e(x_e) = 1$, then localise at $x_e$.
        
        \qedhere
    \end{itemize}
    
    
    
    
\end{proof}

\comment{
If we wanted to show that, in fact, $\RCA$ proves $\RAD$, we would want a method of constructing a $\sqrt{0}$-computable ideal $J \supsetneq \sqrt{0}$ in any ring where $\{ 0 \}$ is not primary. One idea is the following. Since $\{ 0 \}$ is not primary, there are non-nilpotent zero divisors---let $x \notin \sqrt{0}$ be one. Let $k \neq 0$ be such that $xk = 0$. Then, $\idquo{\sqrt{0}}{k} = \{ r \in R: rk \in \sqrt{0} \}$ is an ideal strictly containing $\sqrt{0}$.

The problem is that $\idquo{\sqrt{0}}{k}$ might not be a proper ideal: in fact, $\idquo{\sqrt{0}}{k} = R$ iff $k \in \sqrt{0}$. The next example shows that we can't always take $k$ to be non-nilpotent, and thus this argument fails in general.

\begin{example}
    Let $R = \Z[x,y]/(x^2, xy)$. An arbitrary element of $R$ can be written in the form $ax + p(y) + b$, where $p$ is a polynomial with no constant term. Using this ``normal form'', we can verify that:
    \begin{itemize}
        \item $\sqrt{0} = (x)$.
        
        \item $rs = 0$ iff $x \mid r$ or $x \mid s$.
    \end{itemize}
    Thus, for any non-nilpotent zero-divisor (e.g. $y$), all its annihilators $k$ are nilpotent. Thus the ideal $\{ r: kr \in \sqrt{0} \}$ is in fact $R$ itself.
\end{example}}

Using <>, we can show that $\WWKL$ does not imply $\RAD$.

\begin{definition}
    Given $X \in \CS$, define $\M_X \subseteq \CS$ by $$\M_X \defeq \left\{ Y \in \CS: (\exists n)\, Y \leq_\mathrm{T} X^{[0]} \oplus X^{[1]} \oplus \cdots \oplus X^{[n]} \right\}$$
\end{definition}

\begin{lemma}
    If $X$ is difference random, then $\M_X \models \WWKL$.
\end{lemma}

\begin{proof}
    $\WWKL$ is equivalent to the principle $\MLR$ which says for every $X$, there is $A$ which is Martin-L\"of random relative to $X$ \cite{nies_randomness_2020}. So pick
\end{proof}

\comment{
Difference randoms \cite{diffRand}. Following \cite{nies_shafer_2020}, we treat difference randomness as a reverse mathematical principle, i.e.\ consider the statement
\begin{equation}
    \text{for every $A \in \CS$, there exists a difference random relative to $A$}\tag{$\DIF$}
\end{equation}
in second-order arithmetic

\begin{theorem}
    $\DIF$ does not imply $\RAD$
\end{theorem}

Recall that $X$ is difference random relative to $A$ iff $X$ is Martin-L\"of random relative to $A$, and $A \oplus X \ngeq_\mathrm{T} A'$ \cite{diffRand}. Hence, $\DIF$ implies
\begin{equation}
    \text{for every $A \in \CS$, there exists a Martin-L\"of random relative to $A$}\tag{$\MLR$}
\end{equation}
which is known to be equivalent to weak weak K\H onig's lemma $\WWKL$ \cite{nies_shafer_2020}. Thus, it follows that

\begin{corollary}
    $\mathsf{WWKL}$ does not imply $\RAD$
\end{corollary}
}}

\noindent Currently, we know that $\RAD$ lies between $\RCA$ and $\WKL$, but its exact reverse-mathematical strength remains open. It's possible that $\RAD$ falls into the ``zoo'' of reverse-mathematical principles lying strictly between $\RCA$ and $\WKL$, which have been keenly studied in recent years \cite{sanders_refining_2018,noauthor_reverse_nodate}.

\begin{question}
    What is the exact reverse-mathematical strength of $\RAD$?
\end{question}

\comment{
This version of the theorem is ostensibly weaker, but we show it is still equivalent to $\WKL$.

\begin{theorem}\label{thm:wkl-maximal-primary-exists}
    Let $R$ be a commutative ring with unity, and $I \subseteq R$ an ideal. The following are equivalent:
    \begin{enumerate}
        \item $\WKL$.
        
        \item If $\sqrt{I}$ is maximal, then $I$ is primary.
        
        \item If $\sqrt{I}$ \textbf{exists} and is maximal, then $I$ is primary.
    \end{enumerate}
\end{theorem}

\begin{proof}\
    \begin{description}
        \tfae{i}{ii} This was the ($\Rightarrow$) direction in the proof of Theorem \ref{thm:wkl-maximal-primary} above.
        
        \tfae{ii}{iii} Trivial.
        \pagebreak
        
        \tfae{iii}{i} We construct a computable ring $R$ such that:
        \begin{enumerate}[label=(\alph*)]
            \item $I = \{0\}$ is not primary (equivalently, there are non-nilpotent zero divisors).
            
            \item The nilradical $\sqrt{0}$ is computable.
            
            \item Any proper ideal $J \supsetneq \sqrt{0}$ has PA degree.
        \end{enumerate}
        
        $R$ will be constructed as in \cite[Thm 3.2]{downey_ideals_2007}, but starting with a different ring of coefficients $B$. First, we define $B$, then set $R_0 \defeq B[x_0,x_1,\ldots]$, and finally enumerate $R$ as a c.e.\ subring of the total quotient ring\footnote{$B$ will not be an integral domain.} of $R_0$.
        
        To ensure $R$ satisfies condition (c), we will force it to satisfy:
        \begin{enumerate}[label=(\alph*)]\setcounter{enumi}{3}
            \item Any proper ideal $J \supsetneq \sqrt{0}$ contains a non-zero-divisor.
        \end{enumerate}
        
        $R$ satisfies condition (a) iff $B$ does, and likewise with condition (b). If $B$ satisfies (d), then so does $R$.
        
    \end{description}
\end{proof}
}
\chapter{Noetherian rings}
\label{sec:noeth}

As we proceed, we will consider ideals $\I \subseteq R$ which may not be \textit{computable} relative to $R$, but only \textit{c.e.}\ relative to $R$.

\begin{definition}
    A \textit{$\Si{1}$-ideal} $\I \subseteq R$ is a sequence $\ang{a_0, a_1, a_2, \ldots}$ of elements of $R$, such that for all $i,j \in \N$ and $r,s \in R$, there exists $k = k(i,j,r,s) \in \N$ such that $a_k = a_i r + a_j s$.
\end{definition}

To be precise, we will define basic notions such as subset and equality for $\Si{1}$-ideals. Note that these are defined \textit{non-uniformly}---we don't require there to exist a function witnessing the inclusion/equality.

\begin{definition}\label{defn:s1-id-subs}
    Let $\I = \ang{a_0, a_1, \ldots}$, $\J = \ang{b_0, b_1, \ldots}$ be $\Si{1}$-ideals in $R$, and $r \in R$.
    \begin{enumerate}
        \item We say $r \in \I$ if there exists $m$ such that $r = a_m$.
    
        \item We say $\I \subseteq \J$ if for all $n$, $a_n \in \J$.
        
        \item We say $\I \subsetneq \J$ if $\I \subseteq \J$ and $\J \nsubseteq \I$.
        
        \item We say $\I = \J$ if $\I \subseteq \J$ and $\J \subseteq \I$.\label{defn:s1-id-eq}
    \end{enumerate}
\end{definition}

\begin{definition}\
    \begin{enumerate}
        \item For any $n \in \N$ and $a_0, \ldots, a_{n-1} \in R$, the set $$(a_0,\ldots,a_{n-1})\ \defeq\ \big\{ a_0 r_0 + \cdots + a_{n-1} r_{n-1} : r_0, \ldots, r_{n-1} \in R \big\}$$ is a $\Si{1}$-ideal, called \textit{the $\Si{1}$-ideal generated by $a_0, \ldots, a_{n-1}$}.
        
        \item A $\Si{1}$-ideal $\I$ is \textit{finitely generated} if there exist $n \in \N$ and $a_0, \ldots, a_{n-1} \in R$ such that $\I = (a_0, \ldots, a_{n-1})$, in the sense of Definition \ref{defn:s1-id-subs}.\ref{defn:s1-id-eq}.
    \end{enumerate}
    
\end{definition}

\noindent Typical examples of $\Si{1}$-ideals are principal ideals $(a) \subseteq R$, which in general are not computable. $\Si{1}$-ideals are thus the right notion of ideal for PIDs.

As \cite{sato_reverse_2016} remarks, defining Noetherian rings is difficult, because there are many different, classically equivalent notions which are not equivalent over $\RCA$. For a countable commutative ring $R$, \cite{sato_reverse_2016} considers the following eight definitions of Noetherian:\footnote{\cite{sato_reverse_2016} actually considered their negations, as well as two further conditions which don't characterise Noetherian, but instead a (weaker) condition called a.c.c.p..}

\begin{enumerate}
    \item Every $\D{1}$-ideal $I \subseteq R$ is finitely generated.
    
    \item $R$ has no strictly increasing chain of $\D{1}$-ideals $I_0 \subsetneq I_1 \subsetneq \cdots$.
    
    \item Every chain of $\D{1}$-ideals $I_0 \subseteq I_1 \subseteq \cdots$ in $R$ eventually stabilises.
    
    \item $R$ contains no sequence $a_0, a_1, \ldots$ such that for all $i$, $a_i \notin (a_0,\ldots,a_{i-1})$.
    
    \item For all $a_0, a_1, \ldots$ in $R$, there is $i$ such that $\{ a_j: j \in \N \} \subseteq (a_0,\ldots,a_{i-1})$.\footnote{Simpson \cite{simpson_ordinal_1988} called property (v) \textit{``Hilbertian''}, and noted the equivalence (v)$\Leftrightarrow$(viii).}
    
    \item $R$ has no strictly increasing chain of $\Si{1}$-ideals $\I_0 \subsetneq \I_1 \subsetneq \cdots$.
    
    \item Every chain of $\Si{1}$-ideals $\I_0 \subseteq \I_1 \subseteq \cdots$ in $R$ eventually stabilises.
    
    \item Every $\Si{1}$-ideal $\I \subseteq R$ is finitely generated.
\end{enumerate}

We conduct a full reverse-mathematical analysis of these conditions. Firstly, we reduce them to just five conditions, by showing that, over $\RCA$, (ii) and (iii) are equivalent, and (v), (vii) and (viii) are equivalent.

\begin{description}
    \tfae{ii}{iii} By contradiction, suppose there is a chain of $\D{1}$-ideals $I_0 \subseteq I_1 \subseteq \cdots$ in $R$ which never stabilises. Fix an enumeration $a_0, a_1, \ldots$ of $R$. Define $n_0 = 0$, $j_0$ least such that $a_{j_0} \in I_{n_0}$ and $\ang{n_{k+1},j_{k+1}}$ the least pair\footnote{Such a pair always exists, by assumption. We can find it simply by brute force search.} such that $a_{j_{k+1}} \in I_{n_{k+1}} \setminus I_{n_k}$. Then, $I_{n_0} \subsetneq I_{n_1} \subsetneq \cdots$ is a strictly increasing, uniformly $\D{1}$ chain of ideals. \label{d1-ideals-unbdd-search}
    
    \tfae{iii}{ii} Trivial.
    
    \tfae{v}{viii} Let $a_0, a_1, \ldots$ be an enumeration of $\I$. By assumption, there is $i$ such that $\I \subseteq (a_0,\ldots,a_{i-1})$, whence $\I = (a_0,\ldots,a_{i-1})$ is finitely generated.
    
    \tfae{viii}{vii} Given a chain of $\Si{1}$-ideals $\I_0 \subseteq \I_1 \subseteq \cdots$ in $R$, $\I \defeq \bigcup_{n=0}^\infty \I_n$ is a $\Si{1}$-ideal. Suppose $\I = (b_0,\ldots,b_{m-1})$. Then, in particular, each $b_i \in \I$, so $(\forall i < m) (\exists k_i) (b_i \in I_{k_i})$. The predicate ``$b_i \in I_{k_i}$'' is $\Si{1}$, hence by $\BSi{1}$ (see Definition \ref{defn:ind-ln-bd}), we conclude that $(\exists n)(\forall i < m) (\exists k_i < n) (b_i \in I_{k_i})$. Hence all $b_i$ are in $I_n$, whence the chain stabilises at $n$.
    
    \tfae{vii}{v} $(a_0) \subseteq (a_0,a_1) \subseteq (a_0,a_1,a_2) \subseteq \cdots$ is a nested chain of $\Si{1}$-ideals, so there is some $i$ such that for all $j$, $a_j \in (a_0,\ldots,a_j) = (a_0,\ldots,a_{i-1})$.
\end{description}

(vii)$\Rightarrow$(vi) follows trivially in $\RCA$. The converse also seems computably true; however, the induction necessary to prove it appears to go beyond $\ISi{1}$.

\begin{theorem}[($\RCA+\ISi{2}$)]\label{thm:chain-s1-ideals}
    (vi)$\Rightarrow$(vii).
\end{theorem}

\begin{proof}
    By contradiction, suppose there is a uniformly c.e.\ chain of $\Si{1}$-ideals $\I_0 \subseteq \I_1 \subseteq \cdots$ in $R$ which never stabilises. We will build a strictly increasing, uniformly c.e.\ subsequence $\J_0 \subsetneq \J_1 \subsetneq \cdots$ via a ``moving marker'' priority argument \cite[\S4.3.2]{soare_turing_2016}.
    
    Each $\J_n$ will have a ``marker'' $m_n$ pointing to the ideal $\I_{m_n}$, which it will copy: any element enumerated into $\I_{m_n}$ is also enumerated into $\J_n$. Additionally, each $\J_n$, $n>0$ keeps track of a ``witness'' $x_n \in \J_n$ which it believes is not in $\J_{n-1}$. If we see $x_n$ enter $\J_{n-1}$, then we increment $m_n$, and enumerate $\J_n$ till we find a new witness. We will assume that every element enumerated into $\J_n$ also enters $\J_m$ for $m>n$.
    
    To begin, let $m_0 = 0$, and enumerate the first element of $\I_0$ into $\J_0$. At stage $s>0$, suppose we have finitely enumerated $\J_0, \ldots, \J_{s-1}$, and defined markers $m_0, \ldots, m_{s-1}$ and witnesses $x_1, \ldots, x_{s-1}$. Define $\J_s$ as the current value of $\J_{s-1}$, and $m_s = m_{s-1}+1$. For each $i \leq s$, enumerate the next element of $\I_{m_i}$ into $\J_i, \ldots, \J_s$, and let $x_s$ be some element enumerated\footnote{Choosing $x_s \in \J_{s-1}$ ensures that the process in the next paragraph will be carried out for $i = s$.} into $\J_{s-1}$.
    
    Now, for $i = 1,\ldots,s$ in increasing order, we check if $x_i$ has been enumerated into $\J_{i-1}$ during this stage. If not, nothing need be done. If so, increment $m_i$ by 1, and start enumerating the new $\I_{m_i}$ into $\J_i, \ldots, \J_s$. Since $\J_{i-1}$ is finite so far, we will eventually see an element not in $\J_{i-1}$: this is our new witness $x_i$. We stop enumerating $\I_{m_i}$ when $x_i$ is found. This concludes stage $s$.
    
    Now, we prove by induction on $n$ that every marker $m_n$ eventually stabilises. Formally, we induct on the $\Si{2}$ formula $$\varphi(n) \defeq (\exists s) (\forall t > s) (m_{n,t} = m_{n,s})$$ where $m_{i,s}$ is the value of the $i$th marker at stage $s$.
    
    Firstly, note that $m_0 = 0$ never changes. Now, suppose that $m_0, \ldots, m_{n-1}$ have stabilised already. There are infinitely many elements in $\bigcup \I_k \setminus \I_{m_{n-1}}$: let $k^*$ be least such that $\I_{k^*}$ contains one. Then, $m_n$ will stabilise at some point $\geq k^*$, since after this point, we will eventually enumerate an element of $\bigcup \I_k \setminus \I_{m_{n-1}}$, and thus take a witness $x_n \notin \I_{m_{n-1}}$, which will never be discarded. Therefore, the $x_n$ witness that the chain $\J_0 \subseteq \J_1 \subseteq \cdots$ is strictly increasing.
\end{proof}

We leave open the question of whether $\ISi{2}$ is actually necessary to prove Theorem \ref{thm:chain-s1-ideals}, and instead assume a base theory of $\RCA+\ISi{2}$ for the rest of this chapter. Over $\RCA+\ISi{2}$, we then get four distinct notions of Noetherian, which we name \textit{$\D{1}$-Noetherian}, \textit{$\D{1}$-a.c.c.}, \textit{sequentially Noetherian} and \textit{$\Si{1}$-Noetherian} respectively.

\begin{definition}\label{defn:d1-noeth}
    A ring $R$ is \textit{$\D{1}$-Noetherian} if every $\D{1}$-ideal $I \subseteq R$ is finitely generated.
\end{definition}

\begin{definition}[($\RCA$)]\label{defn:d1-acc}
    A ring $R$ has the \textit{$\D{1}$-a.c.c.}\ if either of the following equivalent conditions holds:
\begin{enumerate}\setcounter{enumi}{1}
    \item $R$ has no strictly increasing chain of $\D{1}$-ideals $I_0 \subsetneq I_1 \subsetneq \cdots$.
    
    \item Every chain of $\D{1}$-ideals $I_0 \subseteq I_1 \subseteq \cdots$ in $R$ eventually stabilises.
\end{enumerate}
\end{definition}

\begin{definition}\label{defn:seq-noeth}
    A ring $R$ is \textit{sequentially Noetherian} if it contains no sequence $a_0, a_1, \ldots$ such that for all $i$, $a_i \notin (a_0,\ldots,a_{i-1})$.
\end{definition}

\begin{definition}[($\RCA+\ISi{2}$)]\label{defn:s1-noeth}
    A ring $R$ is \textit{$\Si{1}$-Noetherian} if any of the following equivalent conditions holds:
\begin{enumerate}\setcounter{enumi}{4}
    \item For all $a_0, a_1, \ldots$ in $R$, there is $i$ such that $\{ a_j: j \in \N \} \subseteq (a_0,\ldots,a_{i-1})$.
    
    \item $R$ has no strictly increasing chain of $\Si{1}$-ideals $\I_0 \subsetneq \I_1 \subsetneq \cdots$.
    
    \item Every chain of $\Si{1}$-ideals $\I_0 \subseteq \I_1 \subseteq \cdots$ in $R$ eventually stabilises.
    
    \item Every $\Si{1}$-ideal $\I \subseteq R$ is finitely generated.
\end{enumerate}
\end{definition}


\begin{figure}[t]
\centering

\def\crad{5}
\def\ballrad{1.1}
\pgfmathsetmacro{\radshift}{0.2+\ballrad}
\def\latshift{0.12}
\pgfmathsetmacro{\outrad}{\crad+\latshift}
\pgfmathsetmacro{\inrad}{\crad-\latshift}

\begin{tikzpicture}[
    ball/.style = {circle, draw, thick, minimum size=2*\ballrad cm,inner sep=0mm,align=center},
    inner/.style = {-latex,thick}
]
    \def\startp(#1)(#2){ ($ (#1) + (#2:\radshift) + (#2+90:\latshift) $) }
    \def\endp(#1)(#2){ ($ (#1) + (#2:\radshift) + (#2-90:\latshift) $) }
    \pgfmathsetmacro{\angshift}{ 2*asin(\radshift/\crad/2) }
    \def\centerarc[#1](#2)(#3:#4:#5)
    { \draw[#1] ($(#2)+(#3:#5)$) arc (#3:#4:#5) }

    \node[ball] at (0,0) (d1n) {\hyperref[defn:d1-noeth]{$\D{1}$-}\\\hyperref[defn:d1-noeth]{Noeth.}\vspace{1mm}};
    \node[ball] at (330:\crad) (d1a) {\hyperref[defn:d1-acc]{$\D{1}$-}\\\hyperref[defn:d1-acc]{a.c.c.}};
    \node[ball] at (210:\crad) (seq) {\hyperref[defn:seq-noeth]{seq.}\\\hyperref[defn:seq-noeth]{Noeth.}\vspace{0.5mm}};
    \node[ball] at (90:\crad) (s1n) {\hyperref[defn:s1-noeth]{$\Si{1}$-}\\\hyperref[defn:s1-noeth]{Noeth.}\vspace{1mm}};

    \centerarc[inner,blue](0,0)(210-\angshift:90+\angshift:\outrad) node[midway,above left] {\hyperref[thm:seq->d1n]{$\ACA$}}; 
    \centerarc[inner](0,0)(90+\angshift:210-\angshift:\inrad) node[midway,below right] {\hyperref[prop:rca-noeth]{$\RCA$}}; 
    
    \centerarc[inner](0,0)(90-\angshift:-30+\angshift:\outrad) node[midway,above right] {\hyperref[prop:rca-noeth]{$\RCA$}}; 
    \centerarc[inner,blue](0,0)(-30+\angshift:90-\angshift:\inrad) node[midway,below left] {\hyperref[thm:seq->d1n]{$\ACA$}}; 
    
    \centerarc[inner,dred](0,0)(330-\angshift:210+\angshift:\outrad) node[midway,below] {\hyperref[thm:wkl-d1acc->seq]{$\WKL$}}; 
    \centerarc[inner](0,0)(210+\angshift:330-\angshift:\inrad) node[midway,above] {\hyperref[prop:rca-noeth]{$\RCA$}}; 

    \draw[inner] \startp(s1n)(270) -- \endp(d1n)(90) node[midway,right] {\hyperref[prop:rca-noeth]{$\RCA$}};
    \draw[inner,dgreen] \startp(d1n)(90) -- \endp(s1n)(270) node[midway,left] {\hyperref[thm:d1n->d1acc]{$\ACA$}};
    
    \draw[inner,blue] \startp(seq)(30) -- \endp(d1n)(210) node[pos=0.65,above left] {\hyperref[thm:seq->d1n]{$\ACA$}};
    \draw[inner,dgreen] \startp(d1n)(210) -- \endp(seq)(30) node[pos=0.65,below right] {\hyperref[thm:d1n->d1acc]{$\ACA$}};
    
    \draw[inner,blue] \startp(d1a)(150) -- \endp(d1n)(330) node[pos=0.35,below left] {\hyperref[thm:seq->d1n]{$\ACA$}};
    \draw[inner,dgreen] \startp(d1n)(330) -- \endp(d1a)(150) node[pos=0.35,above right] {\hyperref[thm:d1n->d1acc]{$\ACA$}};
\end{tikzpicture}

\caption{The different definitions of Noetherian, and logical implications between them over $\RCA+\ISi{2}$. {\color{dgreen} Green} arrows are conjectured.}
\label{fig:noetherian}
\end{figure}
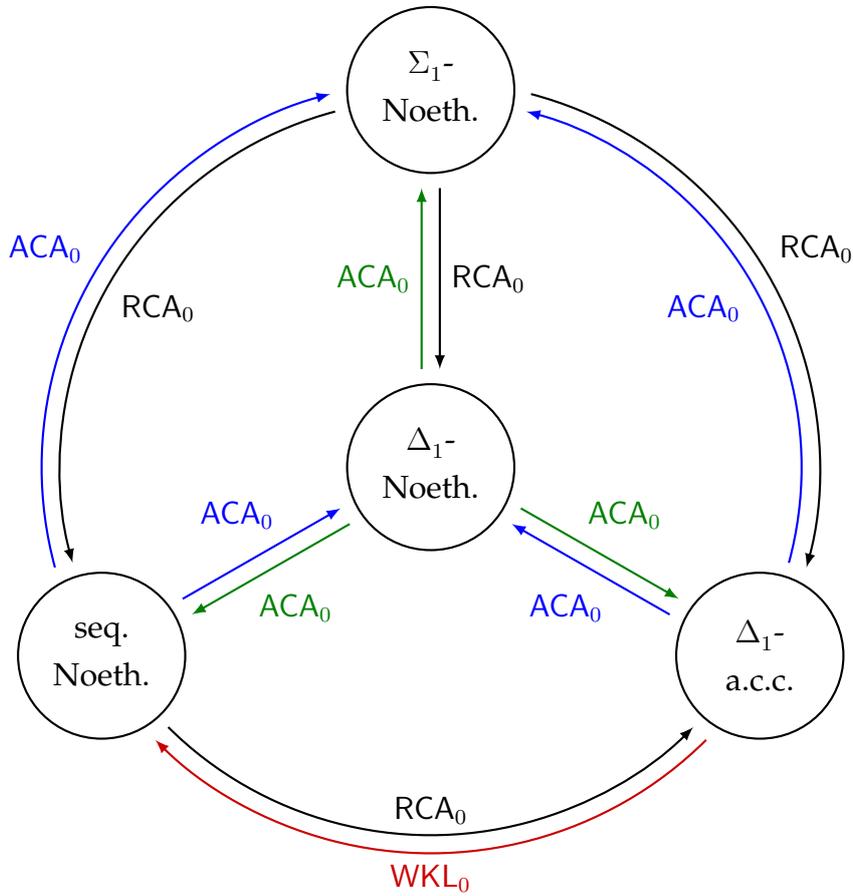

We first establish the trivial (true in $\RCA$) relations between the different notions of Noetherian.

\begin{proposition}[($\RCA$)]\label{prop:rca-noeth}\
    \begin{enumerate}
        \item Every $\Si{1}$-Noetherian ring is $\D{1}$-Noetherian.
        
        \item Every $\Si{1}$-Noetherian ring is sequentially Noetherian.
        
        \item Every sequentially Noetherian ring has the $\D{1}$-a.c.c..
        
        \item Every $\Si{1}$-Noetherian ring has the $\D{1}$-a.c.c..
    \end{enumerate}
\end{proposition}

\begin{proof}\
    \begin{enumerate}
        \item Trivial from the last definition of $\Si{1}$-Noetherian.
        
        \item Trivial from the first definition of $\Si{1}$-Noetherian.
        
        \item Repeat the construction in the proof of (ii)$\Rightarrow$(iii) on page \pageref{d1-ideals-unbdd-search}. Note that $(a_{j_k})_{k \in \N}$ is a ``bad'' sequence: $a_{j_\ell} \in I_{n_k}$ for $\ell \leq k$. Hence, we have $(a_{j_0},\ldots,a_{j_k}) \subseteq I_{n_k}$, but $a_{j_{k+1}} \notin I_{n_k}$.
        
        \item Follows from the previous two items.\qedhere
    \end{enumerate}
\end{proof}

$\ACA$ is enough to prove the equivalence of \textit{all} the given definitions of Noetherian, so it is an upper bound for all arrows in Figure \ref{fig:noetherian}. We now determine the strength of the other implications.

\newtheoremstyle{blue}{\topsep}{\topsep}{\normalfont}{0pt}{\color{blue}\bfseries}{.}{5pt plus 1pt minus 1pt}{}
\theoremstyle{blue}
\newtheorem{bluethm}[theorem]{Theorem}
\begin{bluethm}[($\RCA$)]\label{thm:seq->d1n}
    The following are equivalent:
    \begin{enumerate}
        \item $\ACA$.
        
        \item Every ring with the $\D{1}$-a.c.c.\ is $\Si{1}$-Noetherian.
        
        \item Every ring with the $\D{1}$-a.c.c.\ is $\D{1}$-Noetherian.
        
        \item Every sequentially Noetherian ring is $\Si{1}$-Noetherian.
        
        \item Every sequentially Noetherian ring is $\D{1}$-Noetherian.
    \end{enumerate}
\end{bluethm}

\begin{proof}
    (i)$\Rightarrow$(ii) was observed above, and (ii)$\Rightarrow$(iii), (iii)$\Rightarrow$(v), (ii)$\Rightarrow$(iv), and (iv)$\Rightarrow$(v) are trivial in light of Proposition \ref{prop:rca-noeth}. Therefore, it just remains to prove (v)$\Rightarrow$(i). We will use the contrapositive of (v), and construct a computable ring $R$ with a computable, non-finitely-generated ideal $I \subseteq R$, such that every ``independent sequence'' $a_0, a_1, \ldots$ computes $\varnothing'$.
    
    We use the ring $R$ from \cite[Thm 6.1]{conidis_chain_2010}. In short, start with the ring $R_0 \defeq \Q[\bar{x}] / (x_i x_j : i, j \in \N)$, which consists solely of linear polynomials with the multiplication
    $$\left( q + \sum a_i x_i \right) \left( r + \sum b_i x_i \right)\ =\ qr + \sum (ra_i + qb_i) x_i$$
    Then, we enumerate $A$ as in Lemma \ref{lem:ce-dense}, and when we see $n$ enter $A$, quotient $R_0$ by $x_n = k x_{n+1}$ for an appropriate choice of $k \in \Q$. We can choose the $k$ in such a way to ensure the final ring $R$ is computable---see \cite{conidis_chain_2010} for details.
    
    The ideal $I = (x_0, x_1, \ldots) \subseteq R$ is computable, as it consists of all polynomials with zero constant term, and not finitely generated since $A$ is co-infinite. Now, given an independent sequence $a_0, a_1, \ldots \in R$, we can obtain an independent sequence $a'_0, a'_1, \ldots \in I$ as follows. Set $a'_i = a_i$ till we find $a_n$ with nonzero constant term $q$. Then, for all $i \geq n$, we set $a'_i \defeq q a_{i+1} - r_{i+1} a_n$, where $r_j$ is the constant term of $a_j$.
    
    Now, $I$ is an $\omega$-dimensional vector space over $\Q$ with basis $\{ x_n : n \notin A \}$, and the $a'_i$ are a linearly independent sequence in $I$. For each $n$, let $f(n)$ be the largest variable appearing in $a'_0, \ldots, a'_{n+1}$. By independence, $f(n)$ must be greater than the $n$th element of $A^\complement$. Thus, $f$ dominates $\mu_\zerojump$, and so $a'_0, a'_1, \ldots$ computes $\zerojump$.
\end{proof}

\newtheoremstyle{dred}{\topsep}{\topsep}{\normalfont}{0pt}{\color{dred}\bfseries}{.}{5pt plus 1pt minus 1pt}{}
\theoremstyle{dred}
\newtheorem{dredthm}[theorem]{Theorem}
\begin{dredthm}\label{thm:wkl-d1acc->seq}
    $\WKL$ is equivalent to ``every ring with the $\D{1}$-a.c.c.\ is sequentially Noetherian''.
\end{dredthm}

\begin{proof}
    We will actually work with the contrapositive of the given statement, i.e.\ ``if $R$ has a sequence $a_0, a_1, \ldots$ such that for all $i$, $a_i \notin (a_0,\ldots,a_{i-1})$, then $R$ has a strictly increasing chain of $\D{1}$-ideals''.
    
    \begin{enumerate}
        \ifff Let $\C \subseteq \CS$ consist of the sequences of sets $I = (I_0,I_1,I_2,\ldots)$ such that
        \begin{enumerate}
            \item For all $k \in \N$, $I_k$ is an ideal.
            
            \item For all $k \in \N$, $I_k \subseteq I_{k+1}$.
            
            \item For all $k \in \N$, $a_k \in I_k$ but $a_{k+1} \notin I_k$.\label{thm:wkl-d1acc->seq:strict}
        \end{enumerate}
        By writing the above conditions in first-order logic, we can verify that $\C$ is a $\Ppi[0]{1}$ class. Classically, $\C$ is nonempty, since it contains the sequence $(a_0) \subsetneq (a_0,a_1) \subsetneq (a_0,a_1,a_2) \subsetneq \cdots$. Therefore, $\WKL$ gives a member of $\C$, which is a strictly increasing chain by condition \ref{thm:wkl-d1acc->seq:strict}.\footnote{This can be done more rigorously \`a la Theorem \ref{thm:wkl-maximal-primary}, by building a computable tree $T \subseteq \CSf$ such that $\C = [T]$, and verifying every level of $T$ is nonempty.}

        \iffb We build a computable ring $R$, with a computable bad sequence $a_0, a_1, \ldots$, such that every strictly increasing chain of computable ideals $I_0 \subsetneq I_1 \subsetneq \cdots$ is of PA degree. As in Theorem \ref{thm:wkl-maximal-primary}, fix disjoint c.e.\ sets $A, B$ such that any separator has PA degree. Without loss of generality, we may assume that the complement of $A \cup B$ contains a computable increasing sequence\footnote{If not, redefine $A' = \{ 2a: a \in A \}$ and $B' = \{ 2b: b \in B \}$: then $1, 3, 5, \ldots$ is such a sequence.} $n_0 < n_1 < \cdots$.
        
        The construction is identical to \cite[Thm 3.2]{downey_ideals_2007}. In short, we first set $R_0 \defeq \Z[x_0,x_1,\ldots]$, and add in elements of its field of fractions $\Frac(R_0)$. Begin enumerating $A$ and $B$, and:
        \begin{itemize}
            \item If $n$ enters $A$, add to $R$ all elements of the form $\dfrac{x_n}{p(x_0,\ldots,x_{n-1})}$.
            
            \item If $n$ enters $B$, add to $R$ all elements of the form $\dfrac{x_n-1}{p(x_0,\ldots,x_{n-1})}$.
        \end{itemize}
        Let $R \subseteq \Frac\! \big( \Z[\bar{x}] \big)$ be the subring generated by all the above additions. Then, $R$ is c.e., and hence computably isomorphic to a computable ring \cite[8]{downey_ideals_2007}.
        
        We claim the sequence $x_{n_0}, x_{n_1}, x_{n_2}, \ldots$ is bad. By $\Ppi{1}$ induction on $i$, we will show that for all $k \geq i$, $x_{n_k} \notin (x_{n_0},\ldots,x_{n_{i-1}})$. The base case $i=0$ follows since $x_{n_k} \notin \{ 0 \}$. The case $i=1$ follows since $n_k \notin A$, so $x_{n_k} / x_{n_0} \notin R$ and thus $x_{n_k} \notin (x_{n_0})$.
        
        Now, assume the inductive hypothesis for $i-1$. Fixing $k \geq i$, we have, in particular, that $x_{n_k} \notin (x_{n_0},\ldots,x_{n_{i-2}})$. By contradiction, we will show that $x_{n_k} \notin (x_{n_0},\ldots,x_{n_{i-1}})$. That is, we suppose that $x_{n_k} \in (x_{n_0},\ldots,x_{n_{i-1}})$, so we can write $x_{n_k} = \sum_{j<i} r_j x_{n_j}$ for some $r_j \in R$. We will rewrite $x_{n_k} = \sum_{j<i-1} r'_j x_{n_j}$ for some $r'_j \in R$, showing $x_{n_k} \in (x_{n_0},\ldots,x_{n_{i-2}})$ and contradicting our inductive assumption.
        
        As an arbitrary element of $R$, $r_{i-1}$ must have the form
        $$r_{i-1}\ =\ f\ +\ \sum_{\ell \in A'} g_\ell x_\ell / p_\ell\ +\ \sum_{\ell \in B'} h_\ell (x_\ell-1) / q_\ell$$
        where $f \in \Z[\bar{x}]$, $A' \subseteq A$ and $B' \subseteq B$ are finite, and $p_\ell, q_\ell$ are elements of $\Z[x_0,\ldots,x_{\ell-1}]$. We will show that each summand in $r_{i-1}$ can be ``moved'' into a different $r_j$, hence we can write $x_{n_k}$ as a linear combination of $x_{n_0},\ldots,x_{n_{i-2}}$. Let $c_\ell / d_\ell$ be an arbitrary summand in $r_{i-1}$, for $c_\ell, d_\ell \in \Z[\bar{x}]$.
        \begin{itemize}
            \item If $d_\ell = d'_\ell\, x_{n_{i-1}}$ for some $d_\ell \in \Z[\bar{x}]$, then we can take $c_\ell / d_\ell$ out of $r_{i-1}$, and put $c_\ell / d'_\ell x_{n_j}$ into $r_j$ for some $j < i-1$.
            
            \item Otherwise, the term $(c_\ell / d_\ell) x_{n_{i-1}}$ contains a factor of $x_{n_{i-1}}$, which must be cancelled out by another $r_j$, $j < i-1$, since $x_{n_k}$ contains no factor of $x_{n_{i-1}}$. Hence, $c_\ell = c'_\ell x_{n_j}$, and we can take $c_\ell / d_\ell$ out of $r_{i-1}$, and put $c'_\ell x_{n_{i-1}} / d_\ell$ into $r_j$.
        \end{itemize}
        
        This rewrite shows that $x_{n_k} \in (x_{n_0},\ldots,x_{n_{i-2}})$, contradicting our inductive assumption. Hence, $x_{n_0}, x_{n_1}, x_{n_2}, \ldots$ is a computable bad sequence, as required.

        
%
%

        
        As shown in \cite[Thm 3.2]{downey_ideals_2007}, every nontrivial proper ideal of $R$ has PA degree. If $I_0 \subsetneq I_1 \subsetneq \cdots$ is a strictly increasing chain, then $I_1$ is a nontrivial proper ideal of $R$: hence the chain has PA degree.\qedhere
        
%
%
%
%
%
%
%
    \end{enumerate}
\end{proof}

\newtheoremstyle{dgreen}{\topsep}{\topsep}{\normalfont}{0pt}{\color{dgreen}\bfseries}{.}{5pt plus 1pt minus 1pt}{}
\theoremstyle{dgreen}
\newtheorem{dgreenthm}[theorem]{Conjecture}
\begin{dgreenthm}[($\RCA$)]\label{thm:d1n->d1acc}
    The following are equivalent:
    \begin{enumerate}
        \item $\ACA$.
        
        \item Every $\D{1}$-Noetherian ring is $\Si{1}$-Noetherian.
        
        \item Every $\D{1}$-Noetherian ring is sequentially Noetherian.
        
        \item Every $\D{1}$-Noetherian ring has the $\D{1}$-a.c.c..
    \end{enumerate}
\end{dgreenthm}

    We previously observed (i)$\Rightarrow$(ii), and by Proposition \ref{prop:rca-noeth}, (ii)$\Rightarrow$(iii) and (iii)$\Rightarrow$(iv) are trivial. Therefore, it just remains to prove (iv)$\Rightarrow$(i). Essentially, this would require us to %
%
    construct a computable ring $R$ with a uniformly computable, nonstabilising chain of ideals $I_0 \subseteq I_1 \subseteq I_2 \subseteq \cdots$, such that every non-finitely-generated ideal $J \subseteq R$ computes $\varnothing'$. 
    Following Conidis \cite{conidis_chain_2010,conidis_infinite_2014}, one idea would be to take $A$ as in Lemma \ref{lem:ce-dense}, and construct a ring so that every non-finitely-generated ideal computes an infinite subset of $A^\complement$. However, we have not yet found a construction that works.

Noetherian rings have been studied before in reverse math, but the definitions have not been standard until now. For example, Conidis \cite{conidis_chain_2010,conidis_computability_2019,conidis_computability_2021} took $\D{1}$-a.c.c.\ as his definition of ``Noetherian'', while Simpson \cite{simpson_ordinal_1988} used $\Si{1}$-Noetherian. We hope the results of this section allow a finer analysis of Noetherian rings in reverse mathematics. In particular, one could analyse the reverse-mathematical strength of previously studied theorems, but with a different notion of Noetherian. For example, Conidis \cite{conidis_chain_2010, conidis_computability_2019} proved that ``every Artinian ring has the $\D{1}$-a.c.c'' is equivalent to $\WKL$. From this result, we can deduce:

\begin{corollary}
    $\WKL$ is equivalent to
    \begin{equation}
        \text{``every Artinian ring is sequentially Noetherian''}\tag{$*$}
    \end{equation}
\end{corollary}

\begin{proof}
    In one direction, given an Artinian ring $R$, $\WKL$ proves that $R$ has the $\D{1}$-a.c.c.\ \cite{conidis_chain_2010, conidis_computability_2019}, and then that $R$ is sequentially Noetherian by Theorem \ref{thm:wkl-d1acc->seq}.
    
    In the other direction, $\RCA$ proves that sequentially Noetherian implies $\D{1}$-a.c.c. (Proposition \ref{prop:rca-noeth}), hence over $\RCA$, ($*$) implies ``every Artinian ring has the $\D{1}$-a.c.c'', which implies $\WKL$.
\end{proof}

Simpson \cite{simpson_ordinal_1988} showed that, over $\RCA$, the theorem ``for every field $K$ and $n \in \N$, $K[x_1,\cdots,x_n]$ is $\Si{1}$-Noetherian'' is equivalent to $\mathrm{WO}(\omega^\omega) \defeq$ ``$\omega^\omega$ is well-ordered''. Proposition \ref{prop:rca-noeth} shows that $\RCA + \mathrm{WO}(\omega^\omega)$ also proves this statement for the other notions of Noetherian, but we don't know if these reverse. We would also like to see a study of the more general version of Hilbert's basis theorem: ``if $R$ is $X$-Noetherian, then $R[x]$ is $Y$-Noetherian'', where $X$, $Y$ are chosen from our four notions of Noetherian.

\section{The a.c.c.p.}
\label{sec:accp}

Classically, there is a weakening of the Noetherian chain condition, which only requires every ascending chain of \textit{principal} ideals to stabilise. This is called the \textit{ascending chain condition on principal ideals}, or a.c.c.p.\ for short. The a.c.c.p.\ is of interest because it is often sufficient to prove many of the consequences of Noetherian-ness.

In $\RCA$, there are two sensible notions of a.c.c.p., as noted by \cite{sato_reverse_2016}. Again, assuming $\RCA+\ISi{2}$ as our base theory, we can show that they are equivalent.

\begin{theorem}[($\RCA+\ISi{2}$)]\label{thm:accp}
    The following are equivalent for a ring $R$:
    \begin{enumerate}
        \item There is no sequence $a_0, a_1, \ldots$ in $R$ such that $(a_0) \subsetneq (a_1) \subsetneq \cdots$.
        
        \item Every sequence $(a_0) \subseteq (a_1) \subseteq \cdots$ eventually stabilises.
    \end{enumerate}
\end{theorem}

\begin{proof}
    (ii)$\Rightarrow$(i) is trivial. For (i)$\Rightarrow$(ii), given a nonstabilising sequence $(a_0) \subseteq (a_1) \subseteq \cdots$, we can use the same priority argument as in the proof of Theorem \ref{thm:chain-s1-ideals} to construct a strictly increasing sequence.
\end{proof}

\begin{definition}[($\RCA+\ISi{2}$)]
    $R$ has the \textit{a.c.c.p.}\ if the conditions of Theorem \ref{thm:accp} are satisfied.
\end{definition}

$\RCA+\ISi{2}$ clearly proves that a $\Si{1}$-Noetherian ring has the a.c.c.p.. $\ACA$ appears to be necessary to show the other notions of Noetherian imply a.c.c.p.. However, we do not prove this---we leave it as an open question.
\chapter{Integral domains}
\label{sec:int-doms}

\begin{definition}[($\RCA$)]
    An \textit{integral domain} is a ring $R$ with no nonzero zero-divisors.
\end{definition}

Integral domains are those rings satisfying one of the most basic laws of arithmetic---the cancellation of multiplication. There are a wide variety of different subclasses of integral domains, each generalising properties of our favourite rings: $\Z$, fields, polynomial rings, etc. These form a complex web of implications---an extensive diagram is shown on page \pageref{fig:int-dom}.

\begin{figure}
    \centering
    \includegraphics[width=\textwidth]{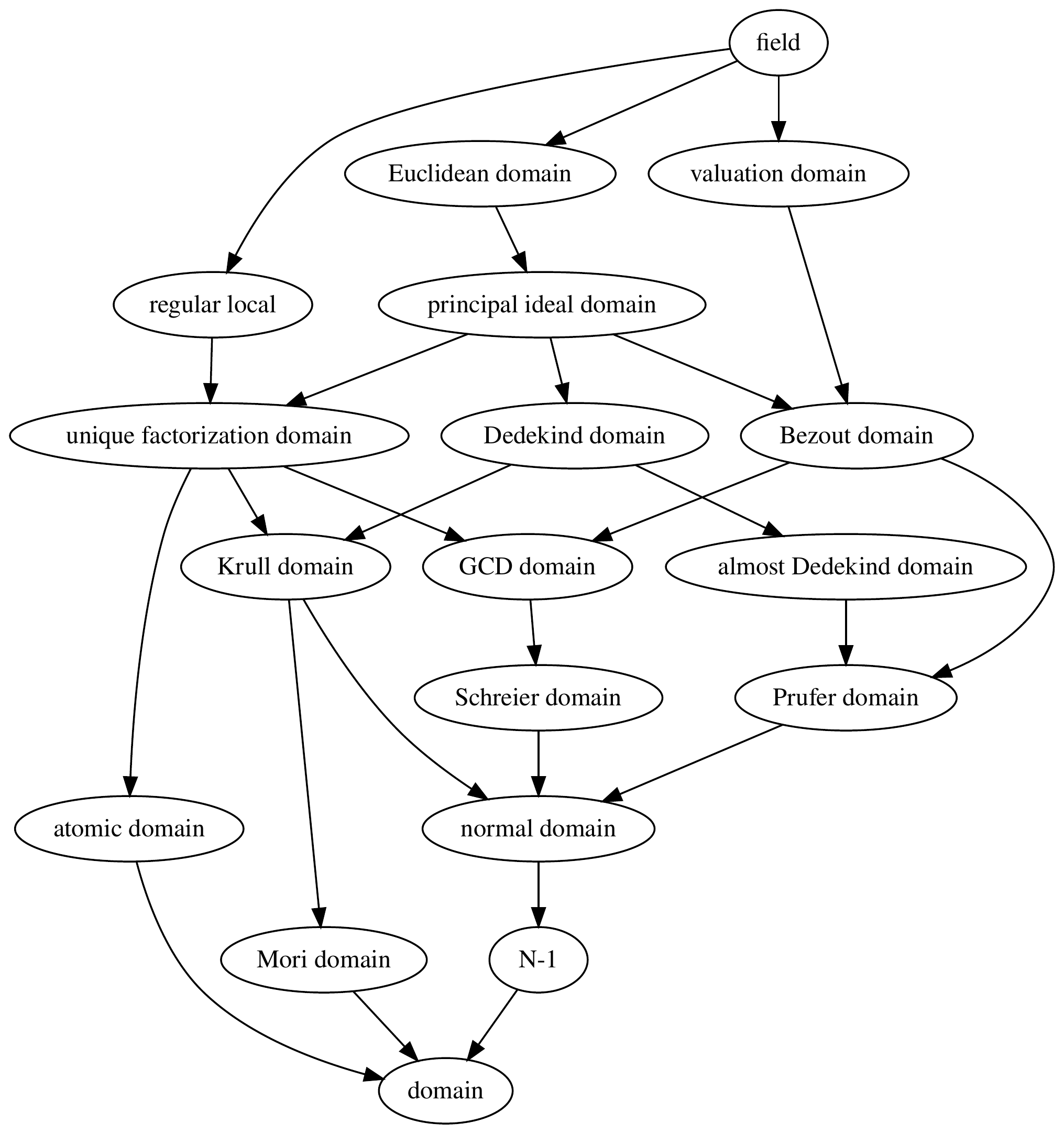}
    \vspace{5mm}
    \caption{Some of the different types of integral domain, and the logical implications between them \cite{schwiebert_field_nodate}. Modified and reproduced with permission.}
    \label{fig:int-dom}
\end{figure}

The aim of this chapter is to initiate a systematic study of the different integral domain properties in reverse mathematics. Some special classes of integral domain have already seen study---Euclidean domains \cite{downey_euclidean_2011, sato_reverse_2016}, PIDs \cite{sato_reverse_2016}, and UFDs \cite{bura_reverse_2013,greenberg_proper_2017}. We extend this study by:
\begin{itemize}
    \item considering additional classes of integral domains (\Bez\ and GCD domains).
    
    \item analysing reverse-mathematically the containments between these different classes of integral domains.
\end{itemize}

\noindent The new types of integral domains we consider often have several (classically) equivalent definitions. Hence, part of the analysis is to determine how hard it is to prove these equivalences, and if they are not equivalent in $\RCA$, to determine which is the right notion in $\RCA$.

First, we will review the existing reverse-mathematical work on integral domains. The first work in this area appears to be the study of Euclidean domains, initiated in \cite{downey_euclidean_2011} and continued in \cite[\S6.4]{sato_reverse_2016}.

\begin{definition}[\cite{downey_euclidean_2011,sato_reverse_2016}]
    An integral domain $R$ is a \textit{Euclidean domain} if there is a function $f\colon R\setminus\{0_R\} \to \N$ such that for all $a,d \in R$ with $d \neq 0_R$, there are $q, r \in R$ such that $a = dq + r$ and $f(r) < f(d)$ if $r \neq 0_R$.
\end{definition}
%

\noindent\cite{downey_euclidean_2011} were concerned with the reverse mathematical strength of
\begin{equation}
    \text{every Euclidean domain has a \textit{minimal} Euclidean function}\tag{$\mathsf{MEF}$}
\end{equation}
They determined that $\mathsf{MEF}$ proves $\ACA$, and conjectured that $\mathsf{MEF}$ is equivalent to $\ACA^+$, the system consisting of $\ACA$ plus the assertion that the $\omega$th Turing jump of any set exists. Meanwhile, \cite{sato_reverse_2016} proved in $\RCA$ that Euclidean domains satisfy a version of \Bez's lemma.

\begin{proposition}[($\RCA$)]
    $R/I$ is an integral domain iff $I \subseteq R$ is a prime ideal.
\end{proposition}

\begin{proof}
    Trivial by the ideal correspondence theorem (Theorem \ref{thm:id-corr}).
\end{proof}

\begin{definition}
    Let $R$ be a ring, and fix $r, s \in R$.
    \begin{enumerate}
        \item $r$ is a \textit{unit} if there exists $s \in R$ such that $rs = 1_R$.
        
        \item $r$ is \textit{irreducible} if $r \neq 0_R$, $r \nmid 1_R$, and for all $a,b \in R$, whenever $r = ab$, then at least one of $a$, $b$ is a unit.
        
        \item $r$ is \textit{prime} if $r \neq 0_R$, $r \nmid 1_R$, and if for all $a,b \in R$, whenever $r \div ab$, then $r \div a$ or $r \div b$.
        
        \item $r, s$ are \textit{associates} (written $r \sim s$) if $r = su$ for a unit $u \in R$.
    \end{enumerate}
\end{definition}

These notions can also be characterised by conditions on the corresponding principal ideals, and these equivalences are trivially provable in $\RCA$.

\begin{proposition}[($\RCA$)]\
    \begin{enumerate}
        \item $r$ is a unit iff $(r) = R$.
    
        \item $r$ is irreducible iff $(r)$ is maximal among proper principal ideals.
        
        \item $r$ is prime iff $(r)$ is a nonzero proper prime ideal.
        
        \item In an integral domain, $r \sim s$ iff $(r) = (s)$.
    \end{enumerate}
\end{proposition}

$\RCA$ can prove basic properties of the associate relation:

\begin{proposition}[($\RCA$)]\
    \begin{enumerate}
        \item $\sim$ is an equivalence relation.
        
        \item If $a \sim b$ and $c \sim d$, then $ac \sim bd$.
        
        \item In an integral domain, if $a \sim b$ and $ac \sim bd$, then $c \sim d$.
    \end{enumerate}
\end{proposition}

\begin{proof}\
    \begin{enumerate}
        \item Reflexivity follows since $1_R$ is a unit. For symmetry, if $a \sim b$, then $a = bu$ where $uv = 1_r$, so $b = buv = av$, whence $b \sim a$. Transitivity follows since the product of units is a unit.
        
        \item We have $a = bu$, $c = dv \implies ac = bd(uv) \implies ac \sim bd$.
        
        \item We have $a = bu$, $ac = bdv$. Then $buc = bdv \implies uc = dv \implies c = d(u^{-1}v)$, so $c \sim d$.\qedhere
    \end{enumerate}
\end{proof}

$\RCA$ can also prove many basic facts about primes and irreducibles. Henceforth, given elements $a,b \in R$, we will use ``$(a/b) \in R$'' as a shorthand for ``there exists $c \in R$ such that $a = bc$'', and use ``$(a/b)$'' as a name for $c$.

\begin{proposition}[($\RCA$)]\label{prop:basic-facts-prime}\
    \begin{enumerate}
        \item If $p$ is prime and $p \div a_1 \cdots a_n$, then there is $j \leq n$ such that $p \div a_j$.\label{prop:prime-div-mult}
        
        \item If $r$ is irreducible, $a$ not a unit, and $a \div r$, then $a \sim r$.\label{prop:irr-div}
        
        \item In an integral domain $R$, every prime element is irreducible. \label{prop:prime->irr}
    \end{enumerate}
\end{proposition}

\begin{proof}\
    \begin{enumerate}
        \item Fixing $a_1, a_2, \ldots$, we proceed by induction on $$\varphi(n)\ \defeq\ \big( p \mid a_1 \cdots a_n \to (\exists j \leq n) (p \div a_j) \big)$$ This formula is $\Si{1} \to \Si{1}$, so by Theorem \ref{lem:ind-bool-comb}, $\RCA$ can carry out this induction.
        
        \item Since $a \div r$, we have $(r/a) \in R$ and $r = a(r/a)$. Then, $(r/a)$ must be a unit since $r$ is irreducible.
        
        \item Suppose $p$ is prime and $p = ab$. Then, $p \div ab$ in particular, so $p \div a$ or $p \div b$ by assumption. WLOG, suppose $p \div a$, so $(a/p) \in R$. Then, $a = p(a/p) = ab(a/p)$. Hence, $b(a/p) = 1$ and $b$ is a unit.\qedhere
    \end{enumerate}
\end{proof}

The converse of Proposition \ref{prop:prime->irr} is not true in general. Integral domains for which the converse holds are called \textit{AP domains}.\footnote{AP is an abbreviation for ``Atoms are Prime'', ``atom'' being an older term for ``irreducible element''.}

\begin{definition}
    An integral domain $R$ is an \textit{AP domain} if every irreducible element in $R$ is prime.
\end{definition}

\section{\Bez\ and GCD domains}
\label{sec:bezout-gcd}

\cite[Thm 6.34]{sato_reverse_2016} proved in $\RCA$ that a version of \Bez's lemma holds in Euclidean domains. This inspired us to take up a reverse-mathematical analysis of \textit{\Bez\ domains}---those in which \Bez's identity holds. However, we will first analyse the weaker notion of \textit{GCD domains}---those in which every pair of elements has a gcd. GCD domains can also be characterised by existence of lcms, or in terms of ideals.

\begin{definition}\label{defn:gcd-lcm}
    Let $R$ be a ring, and fix $b \in R$ and a finite subset $A \subseteq R$.
    \begin{enumerate}
        \item $b$ is the \textit{greatest common divisor} of $A$, written $b = \gcd(A)$, if: \label{defn:gcd}
        \begin{enumerate}
            \item For all $a \in A$, $b \div a$, and
            
            \item For every $c$ satisfying property (i)(a), $c \div b$.
        \end{enumerate}
        
        \item $b$ is the \textit{least common multiple} of $A$, written $b = \lcm(A)$, if:
        \begin{enumerate}
            \item For all $a \in A$, $a \div b$, and
            
            \item For every $c$ satisfying property (ii)(a), $b \div c$.
        \end{enumerate}
    \end{enumerate}
\end{definition}

We first need a simple lemma about gcds.

\begin{lemma}[($\RCA$)]\label{lem:n-gcd}
    For all $a,b,c \in R$, if $\gcd(a,b)$ and $\gcd(ac,bc)$ both exist, then $\gcd(ac,bc) = \gcd(a,b) \cdot c$.
\end{lemma}

\begin{proof}
    Let $d \defeq \gcd(a,b)$ and $e \defeq \gcd(ac,bc)$: we will show that $e \sim dc$. In one direction, since $d \div a,b$, we have $dc \div ac,bc \implies dc \div e$. This implies $(e/dc) \in R$. Conversely, since $e \div ac,bc$, we have $ed \div adc, bdc \implies (e/dc)d \div a,b$. Therefore, $(e/dc)d \div d$, so $e \div dc$ as required.
\end{proof}

\begin{theorem}[($\RCA$)]\label{thm:gcd-lcm}
    For a ring $R$, the following are equivalent:
    \begin{enumerate}
        \item Any two elements of $R$ have a gcd.
        
        \item Any two elements of $R$ have an lcm.
        
        \item For all $a,b$, there is a unique minimal principal $\Si{1}$-ideal containing $(a,b)$.
        
        \item The intersection of two principal $\Si{1}$-ideals is principal.
        
    \end{enumerate}
\end{theorem}

\begin{proof}\
    \begin{description}
        \tfae{i}{ii} Let $d = \gcd(a,b)$. In particular, $d \div a,b$, so $(a/d),(b/d) \in R$. We claim $a(b/d) = \lcm(a,b)$. Clearly $a \div a(b/d)$, but also $a(b/d)d = ab = (a/d)bd$, so $b \div a(b/d) = (a/d)b$.
        
        Now, suppose $c \in R$ is such that $a,b \div c$. Then $ab \div ac, bc$, so $ab \div \gcd(ac,bc) = dc$ by Lemma \ref{lem:n-gcd}. %
        Hence $a(b/d) \div c$ as required.

        \tfae{ii}{i}Let $\ell = \lcm(a,b)$. In particular, $a,b \div \ell$, so $(\ell/a), (\ell/b) \in R$. Also, since $a \div ab$ and $b \div ab$, we have $\ell \div ab$ by definition, so $(ab/\ell) \in R$.
        
        We claim $(ab/\ell) = \gcd(a,b)$. We have $(ab/\ell)(\ell/b)b = ab$, which implies $(ab/\ell)(\ell/b) = a$, so $(ab/\ell) \div a$, and similarly $(ab/\ell) \div b$. Now, suppose $e \div a,b \implies (a/e),(b/e) \in R$. Clearly $a \div a(b/e)$, but also $a(b/e)e = ab = (a/e)be \implies a(b/e) = (a/e)b$, hence $b \div a(b/e) = (a/e)b$. By definition of $\ell$, $\ell \div a(b/e) \implies e\ell \div ab = (ab/\ell)\ell \implies e \div (ab/\ell)$ as required.

%
    \end{description}

    (i)$\Leftrightarrow$(iii) and (ii)$\Leftrightarrow$(iv) are straightforward since $a \div b \iff (b) \subseteq (a)$. We have
    \begin{equation*}
        (a,b) \subseteq \big( \gcd(a,b) \big) \qquad \text{and} \qquad (a) \cap (b) = \big( \lcm(a,b) \big) \qedhere
    \end{equation*}
\end{proof}

\begin{definition}
    An integral domain $R$ is a \textit{GCD domain} if $R$ satisfies any of the conditions in Theorem \ref{thm:gcd-lcm}.
\end{definition}

We can extend each of the conditions in Theorem \ref{thm:gcd-lcm} from two elements to an arbitrary finite number, but in each case, we seem to require $\Si{3}$ induction.

\begin{proposition}[($\RCA+\ISi{3}$)]\label{prop:gcd-induct}
    An integral domain $R$ is a GCD domain iff any of the following holds:
    \begin{enumerate}
        
        \setcounter{enumi}{4}\item Any finite subset of $R$ has a gcd.
        
        \item Any finite subset of $R$ has an lcm.
        
        \item Every finitely generated $\Si{1}$-ideal in $R$ is contained in a unique minimal principal $\Si{1}$-ideal.
        
        \item Finite intersections of principal $\Si{1}$-ideals are principal.
    \end{enumerate}
\end{proposition}

\begin{proof}
    The equivalences (i)$\Leftrightarrow$(v), (ii)$\Leftrightarrow$(vi), (iii)$\Leftrightarrow$(vii), (iv)$\Leftrightarrow$(viii) all follow by induction. In each case, the formula we induct over is $\Si{3}$.
\end{proof}

GCD domains also satisfy one of the most fundamental properties of PIDs and UFDs: every irreducible element is prime.

\begin{proposition}[($\RCA$)]
    GCD domains are AP domains.
\end{proposition}

\begin{proof}
    Suppose $R$ is a GCD domain, $r \in R$ is irreducible, and $r \div ab$. Let $d = \gcd(ar,ab)$. In particular, we have $d \div ar \implies (ar/d) \in R$. The case $ab = 0$ is trivial, so suppose $ab \neq 0$. Then $a,b \neq 0$ since we are in an integral domain, and hence $d \neq 0$.
    
    Now, since $a \div ar$, $a \div ab$ and $d$ is the gcd, we must have $a \div d$, so $(d/a) \in R$. Then, $dr = (d/a)ar = (d/a)(ar/d)d \implies r = (d/a)(ar/d)$. Since $r$ is irreducible, either $(ar/d)$ or $(d/a)$ is a unit. If $(ar/d)$ is a unit, then $ar \sim d$, so $ar \div ab \implies r \div b$. If $(d/a)$ is a unit, then $a \sim d$, and $r \div ar, ab \implies r \div d \implies r \div a$.
\end{proof}

\Bez\ domains are a special type of GCD domains, in which the gcd of $A$ is required to be a \textit{linear combination} of elements of $A$. There are several definitions of \Bez, which turn out to be equivalent in $\RCA$. In particular, this additional condition lowers the amount of induction required to prove the equivalence between the case $\abs{A} = 2$ and the case for arbitrary finite $A$. Thus, the equivalence between these two cases goes through in $\RCA$.

\begin{theorem}\label{thm:bezout-eq}
    For a commutative ring $R$ with unity, the following are equivalent over $\RCA$:
    \begin{enumerate}
        \item Every pair $a,b \in R$ has a gcd $d$, and there are $x,y \in R$ s.t. $ax + by = d$.
        
        \item Every finite set $A \subseteq R$ has a gcd $d$, and $d \in (A)$.
        
        \item For all $a,b \in R$, there exists $d$ such that $(a) + (b) = (d)$.
        
        \item Every finitely generated $\Si{1}$-ideal is principal.
    \end{enumerate}
\end{theorem}

\begin{proof}\
    \begin{description}
        \tfae{i}{iii} Suppose $d \mid a$, $d \mid b$ (i.e.\ $(a/d),(b/d) \in R$), and $d = ax+by$, for some $x,y \in R$. We claim that $(a) + (b) = (d)$. Picking some $ak + b\ell \in (a) + (b)$, we have $$ak + b\ell\ =\ d(a/d)k + d(b/d)\ell\ =\ d \big( (a/d)k + (b/d)\ell \big)\ \in\ (d)$$
        Conversely, picking $dk \in (d)$, we have $$dk\ =\ (ax+by)k\ =\ a(xk) + b(yk)\ \in\ (a) + (b)$$

        \tfae{iii}{iv} Fix a finitely generated $\Si{1}$-ideal $\I = (a_1, \ldots, a_n)$. We induct on the formula $$\varphi(k) = (\exists s) \big[ (a_1,\ldots,a_k)=(s_k) \big]$$
        which is $\Si{1}$, since it can be written as
        \begin{align*}
            \varphi(k) = (\exists s)&(\exists b_1,\ldots,b_k,c_1,\ldots,c_k) \\ &\big( s = b_1 a_1 + \cdots + b_k a_k\, \land (\forall j \leq k)(a_j = c_j s) \big)
        \end{align*}
        
        $\varphi(0)$ is witnessed by $s = 0_R$, and $\varphi(1)$ by $s = a_1$. By induction, assume $(a_1,\ldots,a_k) = (s')$. Then, by assumption, there is $s$ such that $(s') + (a_{k+1}) = (s)$, whence $(a_1,\ldots,a_k,a_{k+1}) = (s)$.


        \tfae{iv}{ii} The ideal $(A)$ is finitely generated, so by assumption, fix $d$ such that $(A) = (d)$. We claim $d = \gcd(A)$. For every $a \in A$, we have $a \in (A) = (d)$, hence $d \mid a$, so $d$ is indeed a common divisor of $A$.
        
        Since $d \in (d) = (A)$, this immediately implies there exist $b_1,\ldots,b_n \in R$ such that $d = a_1 b_1 + \cdots + a_n b_n$, where $A = \{ a_1, \ldots, a_n \}$. It then follows that $d = \gcd(A)$. Fix any other common divisor $e \mid A$: then $(a_i/e) \in R$ for each $i \leq n$. Hence $$d\ =\ \sum_{i=1}^n a_i b_i\ =\ \sum_{i=1}^n e (a_i/e) b_i\ =\ e \cdot \sum_{i=1}^n (a_i/e) b_i$$ so $e \mid d$ as required.

        \tfae{ii}{i} Trivial.\qedhere
    \end{description}
\end{proof}

\begin{definition}
    An integral domain $R$ is \textit{\Bez} if any of the equivalent conditions in Theorem \ref{thm:bezout-eq} holds.
\end{definition}

From definition (ii) of \Bez\ and (v) of GCD domain, it follows trivially (in $\RCA$) that any \Bez\ domain is a GCD domain.

\section{UFDs and PIDs}
\label{sec:ufds}

In this section, we discuss two important types of integral domains: \textit{unique factorisation domains} (UFDs) and \textit{principal ideal domains} (PIDs). These were some of the earliest types of integral domains considered, and they both have rich and well-developed theories with a lot of overlap. Consequently, we have had to order the results in this section carefully to make it clear that there is no circularity.

\subsection{UFDs}

In what follows, let $[n] \defeq \{ 1, \ldots, n \}$.

\begin{definition}\label{defn:ufd}
    An integral domain $R$ is a \textit{unique factorisation domain} (UFD) if every nonzero element $r \in R$ can be written $r = u q_1 \cdots q_n$ for a unit $u$ and irreducibles $q_1, \ldots, q_n$, and this factorisation is unique, i.e.\ for any other irreducible factorisation $r = u' q'_1 \cdots q'_m$, we have $n=m$, and there is a bijective map $h\colon [n] \to [n]$ such that for all $i \leq n$, $q_i$ and $q'_{h(i)}$ are associates.
\end{definition}

UFDs can also be characterised in terms of the existence of \textit{prime} factorisations, and then uniqueness automatically follows. $\RCA$ can prove this characterisation, but first we need some lemmas about UFDs.

\begin{lemma}[($\RCA$)]\label{lem:ufd->ap}
    UFDs are AP domains.
\end{lemma}

\begin{proof}
    Let $q \in R$ be irreducible: we will show $q$ is prime. If $q \div ab$, then $(ab/q) \in R$. Substituting in the unique factorisation of $(ab/q)$ into $q(ab/q) = ab$, it follows that $q$ is in the unique factorisation of $ab$. Multiplying the factorisations of $a$ and $b$ respectively gives another factorisation of $ab$. Hence, up to units, $q$ is in the unique factorisation of $a$ or $b$, so $q \div a$ or $q \div b$ as required.
\end{proof}

\begin{lemma}[($\RCA$)]\label{lem:factor-unq}
    In an integral domain $R$, suppose that $p_1, \ldots, p_n$ are prime, $q_1, \ldots, q_m$ are irreducible, and $p_1 \cdots p_n \sim q_1 \cdots q_m$. Then, $n=m$, and there is a bijective map $h\colon [n] \to [n]$ such that for all $i \leq n$, $p_i$ and $q_{h(i)}$ are associates.
\end{lemma}

\begin{proof}
    Fixing prime $p_1, \ldots, p_n$ and irreducible $q_1, \ldots, q_m$ with $p_1 \cdots p_n \sim q_1 \cdots q_m$, we induct up to $n$ on the $\Si{1}$ formula $$\varphi(k) = (\exists \text{ injective } h\colon [k] \to [m])(\forall i \leq n)(p_i \sim q_{h(i)})$$
    See \cite[50]{bura_reverse_2013} for details.
\end{proof}

\begin{theorem}[($\RCA$)]\label{thm:ufd-defns}
    For an integral domain $R$, the following are equivalent.
    \begin{enumerate}
        \item $R$ is a UFD, in the sense of Definition \ref{defn:ufd}.
        
        \item Every nonzero element $r \in R$ factors into primes and a unit.
    \end{enumerate}
\end{theorem}

\begin{proof}\
    \begin{description}
        \tfae{i}{ii} Follows from Lemma \ref{lem:ufd->ap}.
        
        \tfae{ii}{i} By Proposition \ref{prop:basic-facts-prime}.\ref{prop:prime->irr}, a prime factorisation is an irreducible factorisation. Uniqueness follows from Lemma \ref{lem:factor-unq}.\qedhere
    \end{description}
\end{proof}

By Lemma \ref{lem:ufd->ap}, we will use ``prime'' and ``irreducible'' interchangeably when discussing UFDs.

\begin{proposition}
    UFDs are GCD domains.
\end{proposition}

\begin{proof}
%
    Suppose $R$ is a UFD, and fix finite $A \subseteq R$. By assumption, each $a \in A$ has a unique factorisation $a = u p_1 \cdots p_n$. By bounded $\Si{1}$ comprehension, for each $a \in A$, $\RCA$ can recursively define finite sets $A_0 \supseteq \cdots \supseteq A_n$ by $A_0 = \{ 1, \ldots, n \}$ and $$A_{k+1}\ =\ A_k\ \setminus\ \{ j \in A_k : (\exists b \in R)(p_j = b \cdot p_{\min A_k}) \}$$
    Letting $q_i = p_{\min A_{i-1}}$ and $b_i = \abs{A_i \setminus A_{i-1}}$ as long as $A_{i-1} \neq \varnothing$, we can instead factorise $a = u q_1^{b_1} \cdots q_m^{b_m}$ for non-associate primes $q_1, \ldots, q_m$ and $b_i > 0$. Then $\gcd(A) = q_1^{c_1} \cdots q_m^{c_m}$ where $c_i$ is the minimum of the $b_i$ for all $a \in A$.
\end{proof}

In fact, the above proof shows that UFDs satisfy the stronger definition of GCD domain (Proposition \ref{prop:gcd-induct}), without the assumption of $\ISi{3}$. Thus, we can freely assume the existence of arbitrary gcds in a UFD.

\comment{\color{red}
Definition \ref{defn:ufd} doesn't guarantee the existence of an algorithm to find the factorisation of an element. Indeed, this cannot be done in general: \cite[Thm 6.2]{sato_reverse_2016} proved that $\ACA$ is needed to compute the set of irreducibles in arbitrary rings, and \cite{dzhafarov_complexity_2018} sharpened this result by constructing a UFD whose set of irreducibles is $\Pi_2$-complete.

\begin{definition}
    Given a UFD $R$, a factorisation function on $R$ is a function $f\colon R \to R^{<\omega}$ mapping every
\end{definition}

\begin{definition}
    Given a UFD $R$, a UFD norm on $R$ is a function $f\colon R \to \N$ mapping every
\end{definition}
}

\subsection{PIDs}

\cite{sato_reverse_2016} furthermore considered \textit{principal ideal domains} (PIDs)---integral domains in which every ideal is principal. As mentioned in the previous section, the most natural notion of ideal for PIDs is $\Si{1}$-ideal, since the principal ideal $(a)$ is $\Si{1}$ in general, and may not be computable for every $a$.

\begin{definition}[\cite{sato_reverse_2016}]
    An integral domain $R$ is a ($\Si{1}$-)\textit{PID} if every $\Si{1}$-ideal $\I \subseteq R$ is principal.
\end{definition}

Right from the definition, we see that $\RCA$ proves every PID is \Bez\ (using definition (iv) of \Bez), and that every PID is $\Si{1}$-Noetherian. Recall that $\Si{1}$-Noetherian was the strongest notion of Noetherian, and so $\RCA$ also proves that PIDs satisfy all the other definitions of Noetherian, and that all PIDs have the a.c.c.p..


Classically, every PID is a UFD. The usual proof can be broken down into four steps:
\begin{enumerate}
    \item Show every PID is a Noetherian AP domain.
    
    \item Show every element $r$ has an irreducible factor.
    
    \item Then, show $r$ can be written as a product of irreducibles.
    
    \item Show that any prime factorisation of $r$ is unique up to order and units.
\end{enumerate}

$\RCA$ can carry out step (i): this is Theorem \ref{thm:sato-pid}.\ref{thm:pid-ap}, along with the observation that any PID is Noetherian. $\RCA$ can also carry out step (iv): this is Lemma \ref{lem:factor-unq}. The usual proofs of steps (ii) and (iii) require us to recognise when an element is irreducible or a unit, hence they require $\ACA$ (as these conditions are $\Ppi{2}$ and $\Si{1}$ respectively). \cite{bura_reverse_2013,greenberg_proper_2017} showed that $\ACA$ is necessary for step (iii):

\begin{theorem}[\cite{bura_reverse_2013,greenberg_proper_2017}]
\label{thm:accp->irr-fact}
    $\ACA$ is equivalent to ``in a ring with a.c.c.p., every element has an irreducible factorisation''.
\end{theorem}

\comment{
A simple variation on their argument shows that $\ACA$ is also necessary for step (ii):
\begin{theorem}[($\RCA$)]
\label{thm:accp->irr-fact}
    The following are equivalent:
    \begin{enumerate}
        \item $\ACA$.
        
        \item In a ring with a.c.c.p., every element has an irreducible factorisation.
        
        \item In a ring with a.c.c.p., every non-unit has an irreducible \textit{factor}.
    \end{enumerate}
\end{theorem}

\begin{proof}
    (i)$\Rightarrow$(ii) was observed, and (ii)$\Rightarrow$(iii) is trivial, so we prove (iii)$\Rightarrow$(i). Essentially, we replicate the of Theorem 1.2 in \cite[202]{greenberg_proper_2017}, except $T_s$ is a perfect binary tree of height $s$ instead of a ``fishbone'' of height $s$. We spell out the proof here in elementary terms. Our strategy is to build a computable ring where some element doesn't have an irreducible factor, but every infinite chain of principal ideals computes $\zerojump$.
    
    Informally, we begin with $\Q[x]$, and at each stage, for every variable $x$ added in the previous stage, we introduce a factorisation $x = y_x z_x$. At stage $s$, our variables will resemble a perfect binary tree of height $s$. Now, if we see $n$ enter $\zerojump$ at stage $s$, we ``cut'' the tree above level $n$, by localising at everything above level $n$ which is not in the leftmost branch. We then rebuild the tree with new variables up to height $s$.
    
    Formally, we begin with $\Q[x_{\varepsilon,0}]$. At stage $s$, we add variables $x_{\sigma,s}$ for all $\sigma \in \CSf$ with $\abs{\sigma} \leq s$, and quotient by the ideal $\big( x_{\sigma,s} - x_{\sigma0,s}\, x_{\sigma1,s} : \abs{\sigma} < s \big)$. Now there are two cases:
    \begin{enumerate}
        \item If no $n < s$ entered $\zerojump$ at stage $s$, then quotient by $x_{\sigma,s-1} - x_{\sigma,s}$ for all $\abs{\sigma} < s$.
        
        \item If some $n < s$ entered $\zerojump$ at stage $s$, then we quotient by $x_{\sigma,s-1} - x_{\sigma,s}$ for all $\abs{\sigma} \leq n$, and ``cut'' the tree past the $n$th level, by localising at $x_{\sigma,s-1}$ for all $\sigma$ such that $n < \abs{\sigma} < s$ and $\sigma_i = 1$ for some $n \leq i < s$.
    \end{enumerate}
    
    As the natural maps $R[x] \mapsto R[x,y]/(x-y)$ and $R[x] \mapsto R[x,y,z]/(x-yz)$ are both injective, we can consider every step as an expansion of the ring. Thus, the ring $R$ obtained in the limit is c.e., and by Theorem \ref{thm:ce-struct}, we may assume $R$ is computable.
    
    We claim $x_{\varepsilon,0}$ has no irreducible factor in $R$. The only possible factors of $x_{\varepsilon,0}$ are the $x_{\sigma,s}$, and none of these are irreducible, since:
    \begin{itemize}
        \item If $x_{\sigma,s}$ was localised, then it is a unit, hence not irreducible.
        
        \item If $x_{\sigma,s}$ was not localised, then it is associate to another $x_{\tau,t}$ (possibly itself) which was never ``cut'', and hence has a proper factorisation into non-units.
    \end{itemize}
   
    By assumption, there is an infinite ascending sequence $(a_0) \subsetneq (a_1) \subsetneq (a_2) \subsetneq \cdots$ in $R$. The construction described above forms a so-called ``linear system of trees'' \cite[Defn 3.15]{greenberg_proper_2017}, so Corollary 3.36 of \cite{greenberg_proper_2017} applies, and we obtain a properly decreasing sequence $M_0 \succ M_1 \succ M_2 \succ \cdots$ of multisets of elements of $\{ x_{\sigma,s} \}$ (where $M \succ N$ means $\prod N$ properly divides $\prod M$ in $R$). This sequence of multisets will allow us to compute $\zerojump$.
    
    As in the proof of \cite[Thm 1.2]{greenberg_proper_2017}, for each multiset $M$, for sufficiently large $s$, there is a unique multiset $M[s]$ of elements of $\{ x_{\sigma,s} : \abs{\sigma} = s \}$ so that $\prod M \sim \prod M[s]$. Furthermore, the map $M \mapsto M[s]$ is uniformly computable in $M$ and $s$. To get from $M[s-1]$ to $M[s]$:
    \begin{enumerate}
        \item If no $n < s$ entered $\zerojump$ at stage $s$, then for every $x_{\sigma,s-1} \in M[s-1]$, put $x_{\sigma0,s}$ and $x_{\sigma1,s}$ in $M[s]$.
        
        \item If some $n < s$ entered $\zerojump$ at stage $s$, then for every $x_{\sigma,s-1}$, put $x_{\sigma0,s}$ and $x_{\sigma1,s}$ in $M[s]$ iff $\sigma_i = 0$ for all $n \leq i < s$. (All other $x_{\sigma,s-1}$ are removed).
    \end{enumerate}
    
    \newcommand{\Sig}{\mathrm{Sig}}
    For each $\sigma$ and $s \geq \abs{\sigma}$, let $k_{\sigma,s} \big( M \big)$ count how many copies of $x_{\sigma00\cdots0,s}$ appear in $M[s]$. Define the \textit{$n$-signature of $M$ at stage $s$} to be the $2^n$-tuple $\Sig_{n,s}(M) \defeq \big( k_{\sigma,s}(M) : \abs{\sigma} = n \big)$
    
    
    {\color{orange} dont get $n \to \infty$ and $s \to \infty$ mixed up}
    
    
    \finthis

\end{proof}
}

We have not yet determined the reverse mathematical strength of step (ii), i.e.\ ``in a ring with a.c.c.p., every element has an irreducible \textit{factor}''. We believe that it could be shown equivalent to $\ACA$, using a similar argument to \cite[Thm 1.2]{greenberg_proper_2017}, but using perfect binary trees instead of ``fishbones''. However, the details have not been worked through at the time of writing.

For completeness, here is a proof in $\ACA$ that all PIDs are UFDs.

\begin{corollary}[($\ACA$)]\label{cor:pid->ufd}
    Every PID is a UFD.
\end{corollary}

\begin{proof}
    We saw that $\RCA$ proves every PID $R$ is AP and Noetherian (hence has the a.c.c.p.). By Theorem \ref{thm:accp->irr-fact}, every element of $R$ factors into irreducibles. By AP-ness, this is also a prime factorisation, hence $R$ is a UFD by Theorem \ref{thm:ufd-defns}.
\end{proof}

\subsection{Equivalent definitions of PIDs}

PIDs can be classically characterised in two alternative ways. One is the existence of a \textit{Dedekind--Hasse norm}, a slight generalisation of a Euclidean norm. The other is the (ostensibly weaker) requirement that only every \textit{prime} ideal is principal. We now show that this latter characterisation is provable in $\ACA$.

\begin{lemma}[($\RCA$)]\label{lem:princ-ext}
    Suppose $I \subseteq R$ is an ideal, and $a \notin I$. If $(I,a)$ and $\idquo{I}{a}$ are both principal, then $I$ is principal.
\end{lemma}

\begin{proof}
    Suppose $(I,a) = (b)$ and $\idquo{I}{a} = (c)$. Since $(I,a) = (b)$, we have $(a/b) \in R$, and there are $d \in R$, $i \in I$ such that $b = i + ad$. We will in fact prove that $I = (bc)$. To see that $bc \in I$, note that $c \in \idquo{I}{a}$, so $ac \in I$. Then, $bc = (i + ad)c = ic + (ac)d \in I$.
    
    To see $I \subseteq (bc)$, pick $r \in I$. Then, $r \in (I,a) = (b)$, so $(r/b) \in R$. Since $r \in I$, we get $r(a/b) = (r/b)b(a/b) = (r/b)a \in I$, whence $(r/b) \in \idquo{I}{a} = (c)$. Writing $(r/b) = c\ell$, we have $r = bk = bc\ell \in (bc)$.
\end{proof}

\begin{theorem}[($\ACA$)]\label{thm:aca-prime-princ}
    Suppose $R$ is an integral domain in which every prime $\Si{1}$-ideal is principal. Then, $R$ is a $\Si{1}$-PID.
\end{theorem}

\begin{proof}
    By contrapositive. Suppose $R$ is not a PID: then there is a nonprincipal $\Si{1}$-ideal $\I \subseteq R$. We will construct a nonprincipal prime ideal $\Pp \supseteq \I$.
    
    Fix a standard listing of all pairs $(a,b) \in R^2$. By recursion, we simultaneously build a tree $T \subseteq \CSf$ and associate every finite binary string $\sigma \in T$ with a $\Si{1}$-ideal $\I_\sigma$. As we construct $T$, we will ensure that for every $\sigma \in T$, $\I_\sigma$ is nonprincipal. To begin, we let $T = \{ \varepsilon \}$ and $\I_\varepsilon \defeq \I$. Now, given $\sigma \in T$, there are two cases:
    \begin{enumerate}
        \item If $\I_\sigma$ is prime, set $\Pp \defeq \I_\sigma$, and stop the construction here---we are done.
        
        \item Otherwise, look for the first pair $(a,b) \in R^2$ such that $ab \in \I_\sigma$ but $a,b \notin \I_\sigma$. Then set $\I_{\sigma\ct 0} \defeq (\I_\sigma, a)$ and $\I_{\sigma\ct 1} \defeq \idquo{\I_\sigma}{a}$. Note that $a \in (\I_\sigma, a)$ and $b \in \idquo{\I_\sigma}{a}$. Put $\sigma\ct 0$ into $T$ iff $(\I_\sigma, a)$ is nonprincipal, and put $\sigma\ct 1$ into $T$ iff $\idquo{\I_\sigma}{a}$ is nonprincipal.
    \end{enumerate}
    
    Note that all the $\I_\sigma$ are $\Si{1}$ relative to $\I$. To tell whether a $\Sigma_1^\I$ ideal is prime is $\Pi_2^\I$, and telling if one is principal is $\Sigma_3^\I$. $\ACA$ proves the existence of $\I'''$, which is powerful enough to carry out the construction of $T$ and the $\I_\sigma$. Thus, $\ACA$ proves that $T$ and the $\I_\sigma$ exist.
    
    Now, assume that case (i) never happened. The resulting set $T$ is indeed a tree, and by Lemma \ref{lem:princ-ext}, every $\sigma \in T$ has a successor in $T$. By induction, it follows that $T$ is infinite. By $\WKL$, take a path $\alpha \in [T]$, and define $$\Pp = \bigcup_{n \in \N} \I_{\substr{\alpha}{n}}$$
    $\Pp$ must be a prime ideal, since if there were $a,b \in R$ such that $ab \in \Pp$ but $a,b \notin \Pp$, we would have forced $a \in \Pp$ or $b \in \Pp$ at some stage of the construction. Furthermore, $\Pp$ is nonprincipal, since if $\Pp = (p)$, then $p \in \I_{\substr{\alpha}{n}}$ for some $n$, whence $\I_{\substr{\alpha}{n}} = (p)$. The theorem follows.%
    %
    %
    %
    %
    %
    %
\end{proof}



Now, we consider the characterisation of PIDs in terms of Dedekind--Hasse norms, and show that this is equivalent to $\ACA$ over $\RCA$.

\begin{definition}
    A \textit{Dedekind--Hasse norm} on an integral domain $R$ is a function $f\colon R \to \N$ such that:
    \begin{enumerate}
        \item $f(r) = 0 \iff r = 0_R$.
        
        \item For all nonzero $a, b \in R$, either $b \mid a$ or there exist $x,y \in R$ such that $0 < f(ax + by) < f(b)$.
        
        \item For all nonzero $a, b \in R$, $f(a) \leq f(ab)$.
    \end{enumerate}
    
    A \textit{Dedekind--Hasse domain} (DHD) is an integral domain which admits a Dedekind--Hasse norm.
\end{definition}

Item (iii) is not always included in the definition, since given a function $f\colon R \to \N$ satisfying just (i) and (ii), we can define $f'\colon R \to \N$ satisfying all three by $f'(r) = \min \{ f(ra) : a \neq 0_R \}$. However, this process is not computable, and hence we must assert condition (iii).

Classically, a ring is a DHD if and only if it is a PID. One direction of this equivalence is provable in $\RCA$:

\begin{theorem}
    $\RCA$ proves ``every DHD is a PID''.
\end{theorem}

\begin{proofcite}{henry_455_dedekind--hasse_2013}
    Suppose $(R,f)$ is a DHD, and $\I \subseteq R$ a nonzero $\Si{1}$-ideal. The image $A \defeq f \big( \I \setminus \{ 0 \} \big)$ is a nonempty c.e.\ subset of $\N$. By $\LSi{1}$, $A$ has a least element $n$. Then, we enumerate $\I$ till we find $b \neq 0_R$ with $f(b)=n$.
    
    We claim $\I = (b)$. Pick nonzero $a \in I$. Note that for all $x,y \in R$, $ax + by \in I$, so we can't have $0 < f(ax + by) < f(b)$ by choice of $b$. It follows that $b \mid a$.
\end{proofcite}

Every Euclidean norm is a Dedekind-Hasse norm (choosing $x=1_R$ every time), and hence:

\begin{corollary}[\cite{sato_reverse_2016}]\label{cor:eucl-is-pid}
    $\RCA$ proves ``every Euclidean domain is a PID''.
\end{corollary}

However, $\ACA$ is needed for the converse of this theorem. First, we show that it can be proved in $\ACA$.

\begin{theorem}[($\ACA$)]\label{thm:pid->dhd}
    Every $\Si{1}$-PID is a DHD.
\end{theorem}

\begin{proofcite}{henry_455_dedekind--hasse_2013}
    As we saw in Theorem \ref{cor:pid->ufd}, $\ACA$ proves that every PID is a UFD. In $\ACA$, we can tell which elements are irreducible/prime (as this is $\Pi^0_2$), so given an element $r \in R$, simply search for its factorisation $p_1 \cdots p_k$. We define $f\colon R \to \N$ by mapping $0_R$ to $0$, and $r \neq 0_R$ to $k+1$, where $k$ is the number of irreducibles in the factorisation of $r$.
    
    We claim $f$ is a Dedekind-Hasse norm. (i) is true by definition, and (iii) follows since $f(ab) = f(a) + f(b)$. For (ii), since every PID is \Bez, $d = \gcd(a,b)$ can be written as a linear combination of $a$ and $b$. Then, if $b \sim q_1 \cdots q_k \nmid a \sim p_1 \cdots p_n$, then there is some $q_j$ not associate to any $p_i$. Hence, $f(d) < f(b)$ as required.
\end{proofcite}

For the reversal, we need to construct a computable PID $R$ so that every DHN on $R$ computes $\zerojump$. We want to use our usual method of coding a c.e.\ set $A \subseteq \N$ into the polynomial ring $\Q[\bar{x}]$. Unfortunately, $\Q[\bar{x}]$ is not a PID. However, we can take a localisation to make it a PID, and still retain the ability to code using the $x_i$'s. First, we need a lemma about a certain partial order, which we will use in the construction.

\newcommand{\Nfin}{\N^\mathrm{fin}}
\begin{lemma}[($\RCA$)]\label{lem:nfin}
    Let $\Nfin$ be the collection of sequences in $\N^\N$ which are eventually zero. Define a partial order $\leq$ on $\N^\N$ by $\alpha \leq \beta \iff \forall i\ \alpha_i \leq \beta_i$. Then, $(\Nfin,\leq)$ has a meet for every nonempty $\Si{1}$ subset.
\end{lemma}

\begin{proof}
    Let $A \subseteq \Nfin$ be nonempty. Essentially, the meet $\alpha$ of $A$ is defined by $\alpha_i = \min \{ \beta_i : \beta \in A \}$, but it takes some work to show that this exists in $\RCA$.
    
    Formally, we will construct $\alpha$ by viewing it as a function $\N \to \N$, i.e.\ a set of pairs. Since $A$ is nonempty, fix some $\gamma \in A$. Since $\gamma \in \Nfin$, let $n \in \N$ be such that $(\forall m \geq n)(\gamma_m = 0)$. By $\D{0}$ comprehension, let $K_0 \defeq \{ (m,a) : m < n,\, a < \gamma_m \}$. Via bounded $\Si{1}$ comprehension (Lemma \ref{lem:bdd-s1-comp}), we can define the following finite subsets of $K_0$:
    \begin{align*}
        K_=\ &\defeq\ \big\{ (m,a) \in K_0: (\exists \beta \in A) (a = \beta_m) \big\} \\
        K_>\ &\defeq\ \big\{ (m,a) \in K_0: (\exists \beta \in A) (a > \beta_m) \big\}
    \end{align*}
    
    Again by $\D{0}$ comprehension, let $K \defeq K_= \setminus K_>$. We claim $K$ is a function $n \to \N$: if $(m,a), (m,a') \in K$, then there are $\beta, \beta' \in A$ such that $a = \beta_m$, $a' = \beta'_m$, but for all $\delta \in A$, $a, a' \leq \delta_m$. In particular, $a \leq \beta'_m = a'$ and $a' \leq \beta_m = a$, hence $a = a'$.
    
    Now, fixing $m < n$, the set $A_m = \{ b : (\exists \beta \in A)(b = \beta_m) \}$ is $\Si{1}$, so it has a least element $a_m$ by $\LSi{1}$. Then, $(m,a_m) \in K$. Hence, $K$ is a function $n \to \N$. We finally define $\alpha$ by $$\alpha\ \defeq\ K\ \cup\ \{ (m,0) : m \geq n \}$$
    and this is a function $\N \to \N$, and an element of $\Nfin$. By definition, $\alpha = \bigwedge A$.%
    %
    %
    %
    %
    %
\end{proof}

In fact, $\Nfin$ has a meet for \textit{every} nonempty subset, but $\RCA$ only has enough induction to prove this for $\Si{1}$ subsets.

Before proving the next proposition, we observe that every polynomial $p \in R[\bar{x}]$ can be written as $p = \sum_{\alpha \in F} c_\alpha x^\alpha$ for a unique choice of finite $F \subseteq \Nfin$ and coefficients $c_\alpha \neq 0$ (where $x^\alpha \defeq \prod x_i^{\alpha_i}$ for $\alpha \in \Nfin$). We will call $F$ the \textit{support} of $p$, and denote it $\supp(p)$.

\begin{proposition}[($\RCA$)]\label{prop:pid-polynom}
    There is a computable ring $T \subseteq \Frac \big( \Q[\bar{x}] \big)$ so that:
    \begin{enumerate}
        \item $T$ is a $\Si{1}$-PID.
        \item Every $x_i$ is not a unit in $T$.
        \item For all $i \neq j \in \N$, we have $x_i \nmid x_j$.
    \end{enumerate}
\end{proposition}

\begin{proof}
    Let $M \subseteq \Q[\bar{x}]$ be the set $M \defeq \Q[\bar{x}] \setminus \bigcup_{i=0}^\infty (x_i)$. $M$ is multiplicatively closed: take $p,q \in M$. By assumption, $\supp(p)$ contains at least one element $\alpha$ with $\alpha_i = 0$. Among these, take the ones with maximal $\alpha_0$, then of those, the ones with maximal $\alpha_1$, etc., until we obtain a unique element $\alpha$. Do the same to obtain a ``lexicographically maximal'' element $\beta \in \supp(q)$.
    
    We claim that the pointwise sum $\alpha + \beta \in \supp(pq)$. Otherwise, the cross term $c_\alpha d_\beta x^{\alpha+\beta}$ would have to be cancelled out by another term $c_\gamma d_\delta x^{\gamma+\delta}$, with $\gamma \neq \alpha$, $\delta \neq \beta$, $\gamma+\delta = \alpha+\beta$. But then either $\gamma$ would have to be lexicographically above $\alpha$, or $\delta$ above $\beta$---contradiction. Hence, $pq$ has a term $c_\alpha d_\beta x^{\alpha+\beta}$ not containing $x_i$. As this holds for all $i \in \N$, $pq \in M$.
    
%
%
    
    Let $T$ be the localisation of $\Q[\bar{x}]$ at $M$. We claim $T$ satisfies all the conclusions of the proposition. For (iii), note that the elements of $T$ have the form $p/m$ for $p \in \Q[\bar{x}]$, $m \in M$. Now, consider the product $x_i p/m$. If $x_i p/m = x_j$ for $j \neq i$, then $x_i p = x_j m$. However, this is not possible, since $m$ is not divisible by $x_i$. 
%
This also implies that (ii) holds.
    
    Now, we show $T$ is a $\Si{1}$-PID. Given $p \in \Q[\bar{x}]$, let $\beta \defeq \bigwedge \supp(p)$. Then, $p = x^\beta m$ for some $m \in M$. 
%
%
%
    Thus, up to units, every element of $T$ is a product of $x_i$'s, i.e.\ a monic monomial. %
    Note also that there is a bijective map between these monic monomials and $\Nfin$ (as defined in Lemma \ref{lem:nfin}), where $\alpha \in \Nfin$ corresponds to $x^\alpha$. Furthermore, the relation $\leq$ on $\Nfin$ corresponds exactly to the divisibility order on $T$.
    
    Given a nonzero $\Si{1}$-ideal $I \subseteq T$, let $X_I \defeq \{ \alpha \in \Nfin : x^\alpha \in I \}$. Since $X_I \subseteq \Nfin$ is nonempty and $\Si{1}$, we can define $\alpha \defeq \bigwedge X_I$ by Lemma \ref{lem:nfin}. Now, we claim $I = (x^\alpha)$. The $\subseteq$ direction follows since $x^\alpha$ divides everything in $I$.
    
    Conversely, we show $x^\alpha \in I$. Let $n$ be such that $\alpha_m = 0$ for all $m \geq n$. By definition of $\alpha$, for every $i < n$, there is $\alpha^{(i)} \in X_I$ such that $\alpha^{(i)}_i = \alpha_i$. Now, considering the sum $\sum_{i<n} x^{\alpha^{(i)}}$, we can factor out $x^\alpha$, and we are left with an element $m$ of $M$. Multiplying by $1/m$ gives $x^\alpha \in I$.
\end{proof}

$T$ is very useful in reversals of theorems about PIDs, since we can code into the $x_i$ as we usually would. Here's a simple example.

\begin{theorem}\label{thm:pid-units}
    There is a PID so that the set of units computes $\zerojump$.
\end{theorem}

\begin{proof}
    Localise $T$ at $\{ x_n : n \in \zerojump \}$. Then, $n \in \zerojump \iff x_n$ is a unit.
\end{proof}

Now, we can complete the reversal of ``every PID is a DHD''.

\begin{lemma}[($\RCA$)]\label{lem:dhd-units}
    If $(R,f)$ is a DHD, then $r \in R$ is a unit iff $f(r) = f(1_R)$.
\end{lemma}

\begin{proof}
    First, since $f(a) \leq f(ab)$, we have $f(1_R) \leq f(1_R r) = f(r)$ for any $r \neq 0_R$, hence $1_R$ has minimal norm among nonzero elements of $R$.
    
    \begin{enumerate}
        \ifff If $rs = 1_R$, then $f(r) \leq f(rs) = f(1_R)$. Hence, $f(r) = f(1_R)$ by minimality of $f(1_R)$.
        
        \iffb Suppose $f(r) = f(1_R)$. If $r \nmid 1_R$, then there are $x, y \in R$ such that $0 < f(x+ry) < f(r) = f(1_R)$, contradicting minimality of $f(1_R)$.\qedhere
    \end{enumerate}
\end{proof}

\begin{corollary}\label{cor:aca-pid->dhd}
    ``Every PID is a DHD'' implies $\ACA$.
\end{corollary}

\begin{proof}
    Let $R$ be the PID in Theorem \ref{thm:pid-units}. By assumption, let $f$ be a DHN on $R$. By Lemma \ref{lem:dhd-units}, $f$ computes the set of units of $R$, which computes $\zerojump$.
\end{proof}

In Theorem \ref{thm:pid-units}, we constructed a computable PID $R$ so that every DHN on $R$ computes $\zerojump$ (by computing the units of $R$). However, note that the collection of DHNs on a PID is a $\Pi_2$ set in Baire space. Thus, we expect that this result is not optimal in terms of computability, i.e. one could likely find a computable PID so that every DHN computes some $X >_\mathrm{T} \zerojump$. However, we leave this question open.

Another possible direction is to analyse the strength of ``every PID is a DHD'' for the weakened notion of DHN without condition (iii). Unfortunately, Lemma \ref{lem:dhd-units} fails badly in this case, and it is not clear that such a DHN can determine the units. Hence, one would need a different way of getting computational power from a DHN.

\subsection{Theorems about PIDs}

\cite{sato_reverse_2016} proved several basic results about PIDs in $\RCA$:

\begin{theorem}[{\cite[\S6.4]{sato_reverse_2016}}]\label{thm:sato-pid}
The following are provable in $\RCA$:
    \begin{enumerate}
        \item Every Euclidean domain is a $\Si{1}$-PID.
        
        \item Every PID is an AP domain.\label{thm:pid-ap}
        
        \item For $a$ irreducible/prime in a PID, $(a)$ is a maximal $\D{1}$-ideal.\label{thm:pid-prime-d1}
    \end{enumerate}
\end{theorem}

\begin{corollary}
    In a PID, every nonzero prime ideal is maximal.
\end{corollary}

\begin{proof}
    If $\Pp = (p)$ is a nonzero prime ideal, then $p$ is prime, so the result follows by Theorem \ref{thm:sato-pid}.\ref{thm:pid-prime-d1}.
\end{proof}

Here, we prove some more results concerning PIDs. Notice that in the proof of Theorem \ref{thm:pid-units}, $x_n$ is irreducible/prime iff $n \notin \zerojump$, hence we get a computable PID in which the primes are $\Ppi{1}$ complete. Thus, as a corollary, we get:

\begin{corollary}[($\RCA$)]
    The following are equivalent:
    \begin{enumerate}
        \item $\ACA$.
        
        \item For any PID $R$, the set of units of $R$ exists.
        
        \item For any PID $R$, the set of primes of $R$ exists.
    \end{enumerate}
\end{corollary}

Since the set of irreducibles in any ring is $\Ppi{2}$ in general, we would expect that this complexity for the primes is not optimal. Indeed, we can improve it to $\Ppi{2}$ complete:

\begin{theorem}\label{thm:pid-irr-pi2}
    There is a computable PID whose set of irreducibles is $\Ppi{2}$ complete.
\end{theorem}

\begin{proof}
    We will build a polynomial ring $R$ with variables $x_e$, $e \in \N$, such that $x_e$ is prime iff $W_e$ is infinite. Hence, we will have a many-one reduction from the primes of $R$ to $\Inf = \{ e: W_e \text{ is infinite} \}$, which is $\Ppi{2}$ complete.
    
    Start with $R_0 \defeq \Q[x_i, y_{i,0}, z_{i,0} : i \in \N] / (x_i - y_{i,0} z_{i,0} : i \in \N)$ localised at $$\left[ \bigcup (y_{i,0}) \cup \bigcup (z_{i,0}) \right]^\complement$$
    
    Now, to build $R$, we begin enumerating all $W_e$ in parallel. When a new element enters $W_e$ at stage $s$, we do the following:
    \begin{enumerate}
        \item Localise at $y_{e,t}$, where $t = \max \{ u : y_{e,u} \in R \text{ at stage } s \}$.
        
        \item Freely add elements $y_{e,s}$, $z_{e,s}$ to $R$, i.e.\ let $R_\text{new} = R_\text{old}[y_{e,s}, z_{e,s}]$.
        
        \item Set $x_e = y_{e,s} z_{e,s}$, i.e.\ quotient $R$ by $(x_e - y_{e,s} z_{e,s})$.
        
        \item Localise at $M = \big\{ p \in R : p \text{ contains } y_{e,s} \text{ or } z_{e,s} \text{ but } p \notin (y_{e,s}) \cup (z_{e,s}) \big\}$.
    \end{enumerate}
    
    Note that for any ring $S$ and element $s \in S$, the natural map $S \to S[y,z]/(s-yz)$ is injective. Hence, combining steps (ii) and (iii) above, we can consider this as a proper expansion of $R$. Since we only add elements, and never remove/quotient any, it follows that $R$ is c.e.. Thus, we may assume $R$ is computable by Theorem \ref{thm:ce-struct}.
    
    We claim $R$ is a PID. $R$ is a localisation of $\Q[x_e, y_{e,s}, z_{e,s} : \ang{e,s} \in K] / (x_e - y_{e,s} z_{e,s} : \ang{e,s} \in K)$ for some c.e.\ set $K$. By step (iv), the non-units of $R$ are all contained in $$\bigcup_{\ang{e,s} \in K} (y_{e,s}) \cup (z_{e,s})$$
    So, as in Proposition \ref{prop:pid-polynom}, every element of $R$ is (up to a unit) a product of the $y_{e,s}$ and $z_{e,s}$, so $R$ is a PID by the same argument.
    
    For each $x_e$, its only possible nontrivial splittings are $y_{e,s} z_{e,s}$ where $\ang{e,s} \in K$. If $W_e$ is finite, then the last $y_{e,s}$ we add will never be made a unit; hence $x_e$ is properly reducible into $y_{e,s} z_{e,s}$. Conversely, if $W_e$ is infinite, every $y_{e,s}$ we add will eventually be made a unit, so $x_e$ is irreducible.%
%
%
%
%
%
%
%
%
\end{proof}

There is a well-known characterisation of when the polynomial ring $R[x]$ is a PID: this is exactly when $R$ is a field. In fact, $R[x]$ is famously a Euclidean domain in this case. This characterisation is also provable in $\RCA$, but first we need the following technical lemma about polynomial division.

\begin{lemma}[($\RCA$)]\label{lem:poly-div}
    Suppose $R$ is a field, and fix $p, d \neq 0 \in R[x]$ with $\deg(p) \geq \deg(d)$. Then, there are polynomials $q, r \in R[x]$ such that
    \begin{enumerate}
        \item $p = dq + r$;
        \item $\deg(q) \leq \deg(p) - \deg(d)$;
        \item $\deg(r) < \deg(d)$.
    \end{enumerate}
%
\end{lemma}

\begin{proof}
    Write
    \begin{center}
    $\begin{array}{ccccccccccccc}
p & = & a_n x^n & + & \cdots & + & a_m x^m & + & \cdots & + & a_1 x & + & a_0 \\
d & = &         &   &        &   & b_m x^m & + & \cdots & + & b_1 x & + & b_0 \\
    \end{array}$
    \end{center}
    
    The claim is that there are
    \begin{center}
    $\begin{array}{ccccccccccc}
q & = & c_{n-m} x^{n-m} & + & \cdots          & + & \cdots & + & c_1 x & + & c_0 \\
r & = &                 &   & e_{m-1} x^{m-1} & + & \cdots & + & e_1 x & + & e_0 \\
    \end{array}$
    \end{center}
    with $p = dq + r$.
    
    Substituting the above expressions into equation (i), we get a system of $n+1$ linear equations in $n+1$ variables $e_i$, $c_i$, with coefficients in $a_i$, $b_i$ (see Figure \ref{fig:lin-eq1}). The corresponding $(n+1) \times (n+1)$ matrix $\mathbf{A}$ is upper-triangular (see Figure \ref{fig:lin-eq1}), and all entries on the main diagonal are nonzero (since $b_m \neq 0_R$). Thus, we can obtain a solution for $e_i$, $c_i$ by computing $\mathbf{A}^{-1} \mathbf{a}$, where $\mathbf{a}$ is the $(n+1)$-vector of coefficients of $p$. See Appendix \ref{app:lin-alg}.
\end{proof}

\begin{theorem}[($\RCA$)]\label{thm:rx-pid}
    For an integral domain $R$, the following are equivalent:
    \begin{enumerate}
        \item $R$ is a field.
        \item $R[x]$ is a Euclidean domain.
        \item $R[x]$ is a $\Si{1}$-PID.
    \end{enumerate}
\end{theorem}

\begin{proof}\
    \begin{description}
        \tfae{i}{ii} We claim the degree function $\deg$ is a Euclidean function on $R[x]$. Pick $p,p' \in R[x]$ with $p' \neq 0$. There are two cases:
        \begin{description}
            \item[$\deg(p) < \deg(p')$:] then $q \defeq 0$, $r \defeq p$ satisfy the definition of Euclidean.
            
            \item[$\deg(p) \geq \deg(p')$:] follows from Lemma \ref{lem:poly-div}.
        \end{description}
        
        \tfae{ii}{iii} Already proven (Corollary \ref{cor:eucl-is-pid}).
        
        \tfae{iii}{i} Pick $a \in R$ and consider the $\Si{1}$-ideal $(a,x)$ in $R[x]$. By assumption, there is $b \in R[x]$ such that $(b) = (a,x)$. We must have $b \in R$, else $(b) = (a,x)$ could not contain constants.
        
        Since $x \in (b)$, there is a polynomial $p \in R[x]$ such that $bp = x$. Now, $p$ must be linear, so $p = cx + d \implies (bc)x + bd = x$. Matching coefficients, we must have $bc = 1_R$, $d = 0$.
        
        Hence, $b$ is a unit and $(b) = (a,x) = R[x]$. In particular, $1_R \in (a,x)$, so there are polynomials $q,r \in R[x]$ such that $1_R = aq + xr$. Write $q = q'x + k$ for $k \in R$; then $1_R = ak + x(r + aq')$. Again matching coefficients, we get $r + aq' = 0$, $ak = 1_R$. Thus, $a$ is a unit.\qedhere
    \end{description}
\end{proof}

\subsection{Gauss' lemma}

Now, we return to our study of UFDs. An important theorem about UFDs is that $R$ is a UFD if and only if $R[x]$ is one. One direction is easy: if $R[x]$ is a UFD, then every $r \in R$ has a factorisation in $R[x]$, but by degree considerations, this must actually be a factorisation in $R$. The other direction is nontrivial, and we analyse its proof here.

\begin{definition}
    Let $R$ be a GCD domain, and $p \in R[x]$, $p \neq 0$. The \textit{content of $p$}, $\cont(p)$ is the gcd of its coefficients. We say $p$ is \textit{primitive} if $\cont(p) = 1$.
\end{definition}

For any nonzero $p \in R[x]$, we can take $c = \cont(p)$ and factor $p = cp'$ to get a primitive polynomial $p'$.

Before we prove the theorem, we need a series of lemmas. The name \textit{Gauss' lemma} is commonly given to any of these lemmas.

\begin{lemma}[($\RCA$; {\cite[Lem 12.1.8]{singh_basic_2011}})]\label{lem:gauss}
    Let $R$ be a UFD, and $K$ the field of fractions of $R$. Fix $p,q \in R[x]$.
    \begin{enumerate}
        \item If $p,q$ are primitive, then so is $pq$. \label{lem:prod-prim}
        
        \item $\cont(pq) \sim \cont(p) \cont(q)$.\label{lem:cont-mult}
        
        \item If $p$ is primitive and $q \div p$, then $q$ is primitive. \label{lem:prim-div}
    
        \item If $p$ is primitive, then the following are equivalent: \label{lem:gauss-prime-irr}
        \begin{enumerate}
            \item $p$ is irreducible in $R[x]$.
            \item $p$ is irreducible in $K[x]$.
            \item $p$ is prime in $R[x]$.
            \item $p$ is prime in $K[x]$.
        \end{enumerate}
    \end{enumerate}
\end{lemma}

\begin{proof}\
    \begin{enumerate}
        \item Fix a prime $r \in R$. Since $p,q$ are primitive, both must have some term not divisible by $r$. Picking the terms $a_n x^n$ in $p$ and $b_m x^m$ in $q$ of maximal degree with this property (by $\LPi{1}$), the $x^{n+m}$ term in $pq$ can't be divisible by $r$ either. As this is true for all primes $r$, we must have $\cont(pq) = 1$.

        \item Let $c_p = \cont(p)$ and $c_q = \cont(q)$, and write $p = c_p p'$, $q = c_q q'$. Then $pq = c_p c_q p'q'$, and since $p'q'$ is primitive by \ref{lem:prod-prim}, the result follows.

        \item From (ii), we get $\cont(q) \div \cont(p) \div 1$ as required.

        \item \begin{description}
            \tfae{a}{b} Suppose $p$ is primitive and irreducible in $R[x]$. We have $\cont(p) = 1$ but $p$ is not a unit, so $\deg(p) > 0$. By contradiction, suppose $p$ is reducible as $p = rs$ for non-units $r,s \in K[x]$ (i.e.\ $\deg(r), \deg(s) > 0$). By clearing denominators and dividing off contents, we can find $a,b,c,d \in R$ and primitive $r',s' \in R[x]$ so that $r = (a/b)r', s = (c/d)s'$.
            
            Then $bdp = acr's'$, so taking contents and using that $p, r', s'$ are primitive, \ref{lem:cont-mult} gives $bd \sim ac$, $p \sim r's'$. Now we've properly factored $p$ in $R[x]$, contradicting irreducibility of $p$.

            \tfae{b}{d} $K[x]$ is a PID by Theorem \ref{thm:rx-pid}, so AP by Theorem \ref{thm:sato-pid}.\ref{thm:pid-ap}.

            \tfae{d}{c} Suppose $p$ is primitive and prime in $K[x]$. Fix $r,s \in R[x]$ such that $p \div rs$ in $R[x]$. Then, $p \div rs$ in $K[x]$ too, so by primeness, $p \div r$ or $p \div s$ in $K[x]$. Without loss of generality, suppose $p \div r$ in $K[x]$, i.e.\ $r = ph$ for $h \in K[x]$. Write $r = cr'$ for $c = \cont(r)$, $r'$ primitive. Clearing denominators and taking contents, write $h = (a/b)h'$ for $a,b \in R$, $h' \in R[x]$ primitive.
            
            Then $bcr' = aph'$, so taking contents and using that $r', p, h'$ are primitive, \ref{lem:cont-mult} gives $bc \sim a$, $r' \sim ph'$. Now we see $r = cr' \sim p(ch')$ and $ch' \in R[x]$, so $p$ is prime in $R[x]$.

            \tfae{c}{a} By Proposition \ref{prop:basic-facts-prime}.\ref{prop:prime->irr}. \qedhere
        \end{description}
    \end{enumerate}
\end{proof}

\begin{theorem}[($\RCA+\ISi{3}$)]\label{thm:poly-ufd}
    If $R$ is a UFD, then $R[x]$ is a UFD.
\end{theorem}

\begin{proof}
    We proceed by strong induction on
    \begin{align*}
        \varphi(n)\ =\ (\forall p) \big( &\deg p = n\, \land\, p \text{ primitive} \\
        &\to\, (\exists u \text{ unit})(\exists q_1,\ldots,q_m \text{ prime})(p = u q_1 \cdots q_m) \big)
    \end{align*}
    
    By Lemma \ref{lem:gauss}, the $q_i$ must be primitive, so we can say they are prime/ irreducible in a $\Ppi{1}$ way: they can't be factored into polynomials of strictly lower degree. Hence, $\varphi$ is a $\Ppi{3}$ formula.
    
    If $\deg(p)=0$, then since $p$ is primitive, it is a unit, so we are done. So suppose $\deg(p) > 0$. If $p$ is irreducible, then it is prime by Lemma \ref{lem:gauss}.\ref{lem:gauss-prime-irr}. Otherwise, $p$ is properly reducible into non-units $q,r$. By Lemma \ref{lem:gauss}.\ref{lem:prim-div}, $q,r$ are primitive, so by induction, they have prime factorisations. The product of these factorisations gives the required factorisation of $p$.
    
    We have proven that every \textit{primitive} $p \in R[x]$ has a prime factorisation. This implies \textit{every} $p \in R[x]$ has a prime factorisation, since we can just factor out the content of $p$ as $p = cp'$, and combine the factorisations of $c$ and $p'$. Now the result follows from Theorem \ref{thm:ufd-defns}.
\end{proof}

One possible strategy to reduce the amount of induction required for Theorem \ref{thm:poly-ufd} would be to find a ``nice'' coding of $R[x]$, and do strong induction on the code for $p$. This would take the induction down to $\ISi{2}$. For this to work, we would need a coding $c\colon R[x] \cong \N$ such that whenever $p = qr$ for non-units $q,r$, we have $c(q), c(r) < c(p)$.

The problem is that, for general UFDs, this coding can't be computable. This is because it would allow us to compute the irreducibles of $R[x]$ from $\zerojump$, because we could simply check all $q < p$ to find factorisations of $p$, and ask $\zerojump$ if $q$ is a unit. \cite{dzhafarov_complexity_2018} constructed a UFD $R$ such that $\Irr(R)$ is $\Pi_2$-complete, and since $\Irr(R) = R \cap \Irr \big( R[x] \big)$, it follows that $\Irr \big( R[x] \big)$ is also is $\Pi_2$-complete. For this ring, such a coding must join $\zerojump$ above $\zerodouble$; in particular, it can't be computable (even from $\zerojump$).

More generally, it seems the only way to reduce the induction in Theorem \ref{thm:poly-ufd} is to bound the quantifiers on $p$ and $q_i$. However, if we could bound these, then we could compute the factorisation of any element by a finite search of all elements less than the bound. Essentially, reducing the complexity would require us to \textit{a priori} ``know'' the factorisation of $p$, which we cannot expect in general.

\comment{
\noindent Possible definitions of \textbf{\Bez:}
\begin{enumerate}
    \item Every finitely generated $\D{1}$-ideal is principal;

    \item Every finitely generated $\Si{1}$-ideal is principal;
    
    \item Every pair $a,b \in R$ has a gcd $d$, and there are $x,y \in R$ s.t. $ax + by = d$;
    
    \item Every pair $a,b \in R$ has a gcd $d$, and if $d=1$, there are $x,y \in R$ s.t. $ax + by = 1$;
    
    \item The sum of two principal $\D{1}$-ideals is a $\D{1}$-ideal;
    
    \item The sum of two principal $\D{1}$-ideals is a $\Si{1}$-ideal;
    
    \item The sum of two principal $\Si{1}$-ideals is a $\Si{1}$-ideal.
\end{enumerate}

Definition (vi) is {\color{red} (maybe?)} not as silly as it seems, as the sum of two

($\D{1}$-)ideals is $\Si{1}$ in general.

\vspace{10mm}

{\color{red} Look at theorem: \Bez\ is PID iff Noetherian iff a.c.c.p.}}
\chapter{Conclusion}
\label{sec:conc}

In this thesis, we used the tools of reverse mathematics to analyse several topics in ring theory, particularly radicals, Noetherian rings and integral domains. Every theorem we analysed turned out to be provable in $\ACA$, including Theorem \ref{thm:aca-prime-princ}, the usual proof of which uses Zorn's lemma. We also showed that some key results in commutative algebra actually require $\ACA$ (i.e.\ they are equivalent to $\ACA$ over $\RCA$):
\begin{itemize}
    \item The equivalence of several notions of Noetherian (Theorem \ref{thm:seq->d1n}).
    
    \item Every PID admits a Dedekind--Hasse norm (Corollary \ref{cor:aca-pid->dhd}).
\end{itemize}

Furthermore, we expect that most of the other results that were proven in $\ACA$ will turn out to be equivalent to $\ACA$. This suggests that $\ACA$ is the right axiom system in which to develop countable commutative algebra.

The majority of the results we studied were provable even in $\RCA$. Some were provable in $\RCA$ with extra induction:
\begin{itemize}
    \item The equivalence of weak and strict chain conditions for $\Si{1}$-ideals (Theorem \ref{thm:chain-s1-ideals}, $\ISi{2}$).
    
    \item If $R$ is a GCD domain, then every finite subset has a gcd (Proposition $\ref{prop:gcd-induct}$, $\ISi{3}$).
    
    \item If $R$ is a UFD, then so is $R[x]$ (Theorem \ref{thm:poly-ufd}, $\ISi{3}$).
\end{itemize}

Hence, if one had philosophical objections to working in a nonconstructive system such as $\ACA$, we can still develop the majority of countable commutative algebra ``computably'' in $\RCA \mathrel{(+}\ISi{3}$), including most basic facts about integral domains (\S\ref{sec:int-doms}). However, $\RCA$ still can't prove some important results such as the existence of irreducible factorisations under a.c.c.p.\ (Theorem \ref{thm:accp->irr-fact}), and the equivalence of different definitions of Noetherian (Theorems \ref{thm:seq->d1n} and \ref{thm:wkl-d1acc->seq}).

Finally, our work has given rise to many open questions. Two particularly interesting, related problems are reversing Corollary \ref{cor:pid->ufd} (PIDs are UFDs) and Theorem \ref{thm:aca-prime-princ} (all prime $\Si{1}$-ideals principal $\implies$ PID) in $\ACA$. The most obvious way is to construct a computable non-UFD in which every enumeration of a nonprincipal ideal computes $\zerojump$ (resp.\ a computable non-$\Si{1}$-PID in which every enumeration of a nonprincipal prime ideal computes $\zerojump$). To do this, we need to be able to force complexity on nonprincipal ideals, so that they require $\zerojump$ to be enumerated. We could do this, for example, by forcing the nonprincipal ideals to be $\Ppi{1}$-complete - it seems like novel techniques would be needed to do this.

Here are some of the other problems we'd particularly like to see solved:
\begin{itemize}
    \item Determine the exact reverse-mathematical strength of $\RAD$ (page \pageref{rad-principle}).
    
    \item Determine whether $\ISi{2}$ is necessary for Theorem \ref{thm:chain-s1-ideals}.
    
    \item Prove Conjecture \ref{thm:d1n->d1acc}.
    
    \item Use reverse mathematics to analyse important results about \Bez\ and GCD domains, such as:
    \begin{itemize}
        \item $R$ is \Bez\ iff it is a \textit{Pr\"ufer} GCD domain.
        
        \item The following are equivalent for a GCD domain: UFD, a.c.c.p, Noetherian, atomic.
        
        \item The following are equivalent for a \Bez\ domain: PID, Noetherian, UFD, a.c.c.p., atomic.
    \end{itemize}
    
    
    
    \item Determine whether $\ISi{3}$ is necessary for Theorem \ref{thm:poly-ufd}.
\end{itemize}

\appendix
\chapter{Linear algebra}
\label{app:lin-alg}


\renewcommand{\A}{\mathbf{A}}
\renewcommand{\B}{\mathbf{B}}
\renewcommand{\M}{\mathbf{M}}

This appendix is devoted to the proof of Lemma \ref{lem:poly-div}, and proving in $\RCA$ the necessary theorems of linear algebra.

We will only need to consider square matrices.

\begin{definition}
    Let $K$ be a field. An \textit{$n \times n$ matrix $\A$ over $K$} is an array of elements of $K$:
    $$\begin{bmatrix}
        A_{1,1} & A_{1,2} & \cdots & A_{1,n} \\
        A_{2,1} & A_{2,2} & \cdots & A_{2,n} \\
        \vdots  & \vdots  & \ddots & \vdots  \\
        A_{n,1} & A_{n,2} & \cdots & A_{n,n} \\
    \end{bmatrix}$$
\end{definition}

Notationally, we will use the same letter to refer to a matrix and its elements, but the matrix $\A$ will be in boldface, while its entries $A_{i,j}$ will be italicised. The definitions of matrix multiplication, identity matrix, invertible matrix are as usual. Furthermore, the usual proofs of associativity of matrix multiplication, uniqueness of inverses, etc.\ go through in $\RCA$.

\begin{definition}
    Let $\A$ be an $n \times n$ matrix, and $i,j \leq n$. The \textit{minor submatrix} $\A_{i,j}$ is the $(n-1) \times (n-1)$ matrix obtained from $\A$ by removing the $i$th row and $j$th column.
\end{definition}

\begin{definition}
    The \textit{determinant} of an $n \times n$ matrix $\A$ is given inductively on $n$. If $n=1$, then $\det(\A) = A_{1,1}$. If $n>1$, then $$\det(\A)\ \defeq\ \sum_{i=1}^n (-1)^{n+i} A_{n,i} \det(\A_{n,i})$$
\end{definition}

Although tedious, the usual proof of the Laplace expansion theorem goes through in $\RCA$: that is, we could equally well have done cofactor expansion along a different row/column to define the determinant.

\begin{lemma}[($\RCA$)]\label{lem:det-inv}
    If $\det(\A) \neq 0_K$, then $\A$ is invertible.
\end{lemma}

\begin{proof}
    As usual, we define the adjugate of $\A$ as the matrix $\B$ such that $B_{i,j} = (-1)^{i+j} \det(\A_{j,i})$. Then, we show that $\A\B = \B\A = \det(\A)\mathbf{I}$: the proof uses the aforementioned Laplace expansion theorem.
\end{proof}

\begin{definition}
    A matrix $\A$ is \textit{upper-triangular} if $A_{i,j} = 0$ for all $i>j$.
\end{definition}

\begin{lemma}[($\RCA$)]\label{lem:det-upt}
    The determinant of an upper-triangular $n \times n$ matrix $\A$ is $A_{1,1} \cdot A_{2,2} \cdots A_{n,n}$.
\end{lemma}

\begin{proof}
    Fix $\A$, and for all $k \leq n$, let $$
    \A^{(k)}\ \defeq\ \begin{bmatrix}
        A_{1,1} & A_{1,2} & \cdots & A_{1,k} \\
        A_{2,1} & A_{2,2} & \cdots & A_{2,k} \\
        \vdots  & \vdots  & \ddots & \vdots  \\
        A_{k,1} & A_{k,2} & \cdots & A_{k,k} \\
    \end{bmatrix}$$
    
    By $\Si{0}$ induction on $k$, we prove that for all $k \leq n$, $\det\left(\A^{(k)}\right) = A_{1,1} \cdot A_{2,2} \cdots A_{k,k}$. The result follows taking $k = n$.
    
    The $k=1$ case follows directly from the definition of determinant. Now, suppose $\det\left(\A^{(k)}\right) = A_{1,1} \cdot A_{2,2} \cdots A_{k,k}$. By definition, $$\det\left(\A^{(k+1)}\right)\ =\ A_{k+1,1} \det\left(\A^{(k+1)}_{k+1,1}\right)\ +\ \cdots\ +\ A_{k+1,k+1} \det\left(\A^{(k+1)}_{k+1,k+1}\right)$$ Since $\A$ is upper-triangular, all terms except the last are zero. Hence, $\det\left(\A^{(k+1)}\right) = A_{k+1,k+1} \det\left(\A^{(k+1)}_{k+1,k+1}\right)$. However, $\A^{(k+1)}_{k+1,k+1} = \A^{(k)}$, whence the claim follows.
\end{proof}

\begin{corollary}[($\RCA$)]\label{cor:upt-inv}
    Let $\A$ be an upper-triangular matrix. If every $A_{i,i} \neq 0_K$, then $\A$ is invertible.
\end{corollary}

\begin{proof}
    By Lemma \ref{lem:det-upt}, $\det(\A) \neq 0_K$, so the result follows from \ref{lem:det-inv}.
\end{proof}

Corollary \ref{cor:upt-inv} is enough to prove Lemma \ref{lem:poly-div}. The relevant system of equations and matrix are shown on the following page.

\pagebreak

\begin{landscape}
\hspace{-10mm} $\begin{array}{rrrrrrrrrrrrlrrrl}
    e_0 &     &        &         &           &        &         & + &         b_0\,c_0  &   &                   &   &        &   &                       & = & a_0               \\
        & e_1 &        &         &           &        &         & + &         b_1\,c_0  & + &         b_0\,c_1  &   &        &   &                       & = & a_1               \\
        &     & \ddots &         &           &        &         &   & \vdotswithin{c_0} &   & \vdotswithin{c_1} &   & \ddots &   &                       &   & \vdotswithin{a_0} \\
        &     &        & e_{n-m} &           &        &         & + &     b_{n-m}\,c_0  & + &   b_{n-m-1}\,c_1  & + & \cdots & + &         b_0\,c_{n-m}  & = & a_{n-m}           \\
        &     &        &         & e_{n-m+1} &        &         & + &   b_{n-m+1}\,c_0  & + &     b_{n-m}\,c_1  & + & \cdots & + &         b_1\,c_{n-m}  & = & a_{n-m+1}         \\
        &     &        &         &           & \ddots &         &   & \vdotswithin{c_0} &   & \vdotswithin{c_1} &   &        &   & \vdotswithin{c_{n-m}} &   & \vdotswithin{a_0} \\
        &     &        &         &           &        & e_{m-1} & + &     b_{m-1}\,c_0  & + &     b_{m-2}\,c_1  & + & \cdots & + &  b_{2m-n-1}\,c_{n-m}  & = & a_{m-1}           \\
        &     &        &         &           &        &         &   &         b_m\,c_0  & + &     b_{m-1}\,c_1  & + & \cdots & + &    b_{2m-n}\,c_{n-m}  & = & a_m               \\
        &     &        &         &           &        &         &   &                   &   &         b_m\,c_1  & + & \cdots & + &  b_{2m-n+1}\,c_{n-m}  & = & a_{m+1}           \\
        &     &        &         &           &        &         &   &                   &   &                   &   & \ddots &   & \vdotswithin{c_{n-m}} &   & \vdotswithin{a_0} \\
        &     &        &         &           &        &         &   &                   &   &                   &   &        &   &         b_m\,c_{n-m}  & = & a_n               \\
\end{array}$

\setcounter{MaxMatrixCols}{20}
\newcommand{\z}{{\color{lightgray} 0}}
\newcommand{\cdotg}{{\color{lightgray} \cdots}}
\newcommand{\vdotg}{{\color{lightgray} \vdots}}

\vspace{-20mm}\vfill
{
\hspace{-20mm}\begin{minipage}{0.55\linewidth}
$$\begin{bmatrix}
       1   &   \z   & \cdotg &   \z   &   \z   & \cdotg &   \z   &    b_0    &    \z     & \cdotg &     \z     \\
      \z   &    1   & \cdotg &   \z   &   \z   & \cdotg &   \z   &    b_1    &    b_0    & \cdotg &     \z     \\
    \vdotg & \vdotg & \ddots & \vdotg & \vdotg &        & \vdotg &  \vdots   &  \vdots   & \ddots &     \z     \\
      \z   &   \z   & \cdotg &    1   &   \z   & \cdotg &   \z   &  b_{n-m}  & b_{n-m-1} & \cdots &     b_0    \\
      \z   &   \z   & \cdotg &   \z   &    1   & \cdotg &   \z   & b_{n-m+1} &  b_{n-m}  & \cdots &     b_1    \\
    \vdotg & \vdotg &        & \vdotg & \vdotg & \ddots & \vdotg &  \vdots   &  \vdots   &        &   \vdots   \\
      \z   &   \z   & \cdotg &   \z   &   \z   & \cdotg &    1   &  b_{m-1}  &  b_{m-2}  & \cdots & b_{2m-n-1} \\
      \z   &   \z   & \cdotg &   \z   &   \z   & \cdotg &   \z   &    b_m    &  b_{m-1}  & \cdots &  b_{2m-n}  \\
      \z   &   \z   & \cdotg &   \z   &   \z   & \cdotg &   \z   &    \z     &    b_m    & \cdots & b_{2m-n+1} \\
    \vdotg & \vdotg &        & \vdotg & \vdotg &        & \vdotg &  \vdotg   &  \vdotg   & \ddots &   \vdots   \\
      \z   &   \z   & \cdotg &   \z   &   \z   & \cdotg &   \z   &    \z     &    \z     & \cdotg &     b_m    \\
\end{bmatrix}$$
\end{minipage}%
\hfill%
\begin{minipage}{0.4\linewidth}
    \vspace{10mm}
    \captionof{figure}{
    The system of linear equations and corresponding matrix obtained from the proof of Lemma \ref{lem:poly-div}, in the case when $2m > n$. The case $2m \leq n$ looks similar.
    }
    \label{fig:lin-eq1}
\end{minipage}
}
\end{landscape}

\pagebreak



\chapter{Zorn's lemma}
\label{app:zorn}

The motivation for this section originally came from Theorem \ref{thm:aca-prime-princ}. All of the proofs we could find used Zorn's lemma in an essential way; hence, we wondered if it was possible to formalise those arguments in second-order arithmetic. One might be tempted to say no, since Zorn's lemma for suborders of $\ang{\Pow(\N), \subseteq}$ seems to have an essential third-order quality. Our idea was to pull back the inclusion relation along an indexing of, say, the c.e.\ sets, reducing a third-order problem to a second-order one. As it turns out, this ``pull-back'' doesn't work, and we'll discuss why below.


However, supposing it did work, we would reduce the problem to a second-order version of Zorn's lemma. That is, given a set $A$ at some level $\Gamma$ of the arithmetical hierarchy ($\Gamma$ depends on the complexity of the original index set) and a partial order $\preccurlyeq$ at some level $\Theta$ of the arithmetical hierarchy (depending on the complexity of the sets being indexed), we want to show Zorn's lemma holds for $(A,\preccurlyeq)$.


We now proceed to the formal development in second-order arithmetic. The definition of partial orders (on subsets of $\N$) is as usual. For a partial ordering $\preccurlyeq$, $\prec$ will denote the corresponding strict relation. $\leq$ denotes the usual order relation on $\N$.

\begin{definition}
    Let $(A,\preccurlyeq)$ be a partial order.
    \begin{enumerate}
        \item A \textit{chain} in $(A,\preccurlyeq)$ is a function $f\colon \N \to A$ which is $\prec$-increasing: for all $n$, $f(n) \prec f(n+1)$.
        
        \item Given a chain $f$ in $(A,\preccurlyeq)$, an \textit{upper bound} for $f$ is an element $a \in A$ such that $f(n) \preccurlyeq a$ (equivalently, $f(n) \prec a$) for all $n$.
        
        \item A maximal element $a \in (A,\preccurlyeq)$ is one such that there is no $b \in A$ with $a \prec b$.
    \end{enumerate}
\end{definition}

\begin{definition}
    Let $\Gamma$ and $\Theta$ be classes
    of subsets of $\N$ and $\N^2$, respectively. $\ZL(\Gamma,\Theta)$ is the following statement: for every $\Gamma$ subset $A \subseteq \N$ and $\Theta$ relation $\preccurlyeq$ on $\N$, if $\preccurlyeq$ is a partial order, and every chain $f\colon \N \to (A,\preccurlyeq)$ has an upper bound, then $(A,\preccurlyeq)$ has a maximal element.
\end{definition}

Evidently, if $\Gamma \subseteq \Gamma'$ and $\Theta \subseteq \Theta'$, then $\ZL(\Gamma',\Theta') \models \ZL(\Gamma,\Theta)$. Here is the standard proof of arithmetical Zorn's lemma in $\ACA$.

\begin{theorem}\label{thm:aca->zl}
    $\ACA$ proves $\ZL(\D[1]{0},\D[1]{0})$.
\end{theorem}

\begin{proof}
    Suppose $A$ and $\preccurlyeq$ are $\D[1]{0}$ (arithmetical), and $\preccurlyeq$ is a partial order. By contradiction, suppose $(A,\preccurlyeq)$ has no maximal element. We will construct a chain $f\colon \N \to A$ with no upper bound.
    
    Fix an enumeration $\ang{a_0, a_1, \ldots}$ of $A$. We define $f$ by recursion, starting with $f(0) = a_0$. If $f(n) = a_k$ has been defined, then search for the next $\ell > k$ such that $a_\ell \succ a_k$, and set $f(n+1) = a_\ell$. Since $a_k$ is not maximal, we know we will always find $a_\ell \succ a_k$.
    
    Now, we claim $f$ has no upper bound in $A$. Suppose it did have an upper bound $a_k$ (i.e. $f(n) \prec a_k$ for all $n$). Then, we would have found $a_k$ at some stage of constructing $f$, and thus set $f(j) = a_k$ for some $j$. But $a_k \nprec a_k$, so this is a contradiction.
\end{proof}

Essentially the same proof shows that $\RCA \models \ZL(\Si{1}, \D{1})$. With some care, we can improve this to show:

\begin{theorem}\label{thm:rca-zl}
    $\RCA$ proves $\ZL(\Si{1},\Si{1})$.
\end{theorem}

\begin{proof}
    Suppose $A = \{ a_0, a_1, \ldots \}$ is c.e., and ${\preccurlyeq} = \bigcup_n {\preccurlyeq}_n$ is a c.e.\ partial order. We define $f$ by recursion, starting with $f(0) = a_0$. Then, for $i = 0, 1, \ldots$, we continue searching through pairs $\ang{n,k}$ to find one such that $a_n \succ_k f(i)$, and set $f(i+1) = a_n$. Since $f(i)$ is not maximal, we know there is $a_n \succ f(i)$, and this will be revealed at some finite stage $k$. Hence, $f$ is total. As in the proof of Theorem \ref{thm:aca->zl}, $f$ has no upper bound in $A$.
\end{proof}

Theorem \ref{thm:rca-zl} is optimal, in a sense:

\begin{theorem}[($\RCA$)]\
    \begin{enumerate}
        \item If the usual order relation $\leq$ is $\Theta$, then $\ZL(\Pi_1,\Theta)$ implies $\ACA$.
        
        \item If $\N$ is $\Gamma$, then $\ZL(\Gamma,\Pi_1)$ implies $\ACA$.
    \end{enumerate}
\end{theorem}

\begin{proof}\
    \begin{enumerate}
        \item Let $A$ be as in Lemma \ref{lem:ce-dense}. We consider the carrier set $A^\complement$ (which is $\Pi_1$) under the usual order relation $\leq$. Then, since $A^\complement$ is infinite, it has no maximal element, but any chain in $A^\complement$ computes $\zerojump$.
        
        \item We take $A=\N$ and build a $\Pi_1$ partial order $\preccurlyeq$ on $\N$, using a ``block merging strategy''. The blocks $B_n$ will be intervals in $(\N,\leq)$ such that all $b \in B_n$ are $\preccurlyeq$-incomparable, and $b \prec c$ for all $b \in B_n$, $c > \max(B_n)$. We begin with $B_n = {n}$, i.e.\ ${\preccurlyeq} = {\leq}$.
        
        To merge blocks $B_n$, $B_{n+1}$ means to remove from $\preccurlyeq$ all pairs $(b,c)$, where $b \in B_n$ and $c \in B_{n+1}$. We enumerate $\zerojump$, and if we see $n$ enter $\zerojump$ at stage $s$, and $s \in B_k$ for some $k > n$, then we merge the blocks $B_n,\ldots,B_k$. So, in the final partial order, we will have $a \preccurlyeq b$ iff $a$ is in a strictly earlier block than $b$.
%
        
        More formally, in $\RCA$ we can define markers $m_{n,s}$ by recursion on $s$, where $m_{n,s}$ marks the start of $B_n$ at stage $s$. To begin, $m_{n,0} = n$, and when we see $n$ enter $\zerojump$ at stage $s$, we find the least $j$ with $m_{j,s} > s$, and set $m_{n+1,s+1} = m_{j,s}$, $m_{n+2,s+1} = m_{j+1,s}$, etc. Then, let $a \preccurlyeq b \iff (\forall s)(\exists n<b)(a < m_{n,s} \leq b)$, which is $\Pi_1$.
        
        By bounded $\Si{1}$ comprehension, $\RCA$ can prove the existence of $\substr{\zerojump}{n}$ and $\substr{\varnothing'_s}{n}$ for every $s$ and $n$. Hence, the formula $$\varphi(n)\ \defeq\ (\exists s)(\forall m<n)(m \in \substr{\varnothing'_s}{n} \iff m \in \substr{\zerojump}{n})$$ is $\Si{1}$, so $\RCA$ can prove $(\forall n)\, \varphi(n)$ by induction. Now if $\substr{\zerojump}{n}$ has stabilised at stage $s$, it follows that $m_n$ will henceforth be fixed, so $\RCA$ proves all the $m_n$ stabilise and all the $B_n$ are finite.
        
        $(\N,\preccurlyeq)$ has no maximal element, since for every $n \in B_k$, we have $n \prec m$ for any $m \in B_{k+1}$. By $\ZL(\Gamma,\Pi_1)$, let $f$ be a chain in $(\N,\preccurlyeq)$. By $\Pi_1$ induction, we can prove that $f(n)$ is in block $B_n$ or higher, by inducting on $$\varphi(n) \defeq (\forall s)[f(n) > m_{n,s}]$$
        Hence, for all $n$, $f(n+1) \geq \mu_\zerojump(n)$, and so $f$ computes $\zerojump$.\qedhere
    \end{enumerate}
\end{proof}

\begin{corollary}[($\RCA$)]\label{cor:zl->aca}
    $\ZL(\Gamma,\Theta)$ is equivalent to $\ACA$ if
    \begin{enumerate}
        \item $\Gamma \supseteq \Pi_1$ and $\leq$ is $\Theta$, or
        
        \item $\N$ is $\Gamma$ and $\Theta \supseteq \Pi_1$.
    \end{enumerate}
\end{corollary}

\renewcommand{\M}{\mathcal{M}}
Now, say we are working in a model $\M$ of second-order arithmetic, and have a collection of $\Si{1}$-ideals which we are trying to apply Zorn's lemma to, e.g.: $$\K = \{ \I : \I \text{ is a nonprincipal }\Si{1}\text{-ideal} \}$$
Being $\Si{1}$ in $\M$, each element of $\K$ has an enumeration which exists in $\M$. So, the idea would be to index all the enumerations in $\M$, and pull back along the indexing to obtain a first-order partial ordering, to which $\ZL(\Gamma,\Theta)$ can be applied. The problem is that we cannot index all the possible enumerations in the model, as there may be uncountably many, e.g.\ when $\M = \Pow(\omega)$ is the full $\omega$-model.

We could attempt to fix this using an \textit{internal} notion of computability. The idea is we have a universal $\Si{1}$ formula $\varphi(e,n,X)$ such that for all $\Si{1}$ formulae $\psi(n,X)$, we can (in $\RCA$) find $e$ such that $$(\forall X)(\forall n)\big(\varphi(e,n,X) \leftrightarrow \psi(n,X) \big)$$
Defining $W_e^X \defeq \{ n: \varphi(e,n,X) \}$ (but not necessarily assuming this set exists), and given some nonprincipal $\Si{1}$-ideal $\I$, we can then look at the set $$K = \{ e : W_e^\I \text{ is a nonprincipal ideal} \}$$ which is $\Pi_3^\I$. We define a relation on $K$ by $e \preccurlyeq e' \iff W_e^\I \subseteq W_{e'}^\I$, which is $\Pi_2^\I$.

So, it seems that we have successfully reduced the problem to a second-order one. However, now another problem arises: since $\preccurlyeq$ is $\Pi_2^\I$, the chains in $K$ that we are trying to defeat are no longer $\I$-computable, but only $\I''$-computable. We can define an internal notion of $\I''$-computability, but the union of (internally) $\I$-c.e.\ sets indexed by an (internally) $\I''$-computable function may not be (internally) $\I$-c.e.\ itself---in general, it will only be (internally) $\Sigma_3^\I$.

So, to ensure closure under $\I''$-computable chains, we could instead look at indices for $\Sigma_3^\I$ nonprincipal ideals. However, now the inclusion relation is $\Pi_4^\I$, so the chains we need to defeat are $\I^{(4)}$-computable, so we would need to pass to $\Sigma_5^\I$ nonprincipal ideals to ensure closure. One can see that we will never be able to ``catch our tail''.

It is disappointing that the principles $\ZL(\Gamma,\Theta)$ don't seem to be applicable in reverse mathematics as we might have hoped. Nonetheless, we have left the results in this appendix, as we think they are interesting in their own right. It would be interesting to look at the statements $\ZL(\Gamma,\Theta)$ from the perspective of \textit{Weihrauch reducibility} \cite{brattka_weihrauch_2021}, where one could obtain a more fine-grained analysis than the crude classification we gave in Theorem \ref{thm:rca-zl} and Corollary \ref{cor:zl->aca}.

\phantomsection
\addcontentsline{toc}{chapter}{Bibliography}

\emergencystretch=2em
\printbibliography

\end{document}